\documentclass[11pt]{article}
\usepackage[T1]{fontenc}
\usepackage{lmodern}
\usepackage[margin=27mm]{geometry}
\usepackage{amsmath,amssymb,amsthm}
\usepackage{microtype}
\usepackage[colorlinks=true,linkcolor=blue,citecolor=blue,urlcolor=blue]{hyperref}
\usepackage{xcolor}
\usepackage{tocloft}
\usepackage{newtxtext}
\usepackage{newtxmath}

\numberwithin{equation}{section}
\newtheorem{theorem}{Theorem}[section]
\newtheorem{lemma}[theorem]{Lemma}
\newtheorem{proposition}[theorem]{Proposition}
\theoremstyle{remark}
\newtheorem{remark}[theorem]{Remark}

\newcommand{\T}{\mathbb T}
\newcommand{\R}{\mathbb R}
\newcommand{\Z}{\mathbb Z}

\theoremstyle{remark}
\newcommand{\eps}{\varepsilon}

\newcommand{\norm}[2]{\|#1\|_{#2}}
\newcommand{\abs}[1]{|#1|}

\DeclareMathOperator{\supp}{supp}
\title{Almost sharp global well-posedness and low-regularity
	blow-up for the two-dimensional periodic Zakharov system}
\author{Jiajun Wang}
\date{}
\begin{document}
\maketitle

\begin{abstract}
    We study global well-posedness of the Zakharov system
    on an arbitrary fixed rectangular two-dimensional torus
    for initial data in $H^s\times L^2\times H^{-1}$,
    where $s$ denotes the Sobolev regularity of the
    Schr\"odinger component. Kishimoto \cite{KishimotoLWP}
    established local well-posedness for $1/2\leq s\leq1$,
    whereas his global result \cite{KishimotoI} required
    $9/14<s<1$ under the condition that the Schr\"odinger
    mass is below the Euclidean cubic ground-state
    mass. We extend the global well-posedness result to every
    $1/2<s<1$,
    together with polynomial bounds for the Sobolev norms,
    under the same mass condition.
	
	The proof refines Kishimoto's resonantly modified energy
	by adding quartic and quintic corrections. The main
	difficulty is to cancel the interactions that obstruct
	the global iteration without introducing singular
	coefficients near resonance or losing control of the
	higher-order terms generated by the corrections.
	Our choices of coefficients exploit cancellations between
	successive nonlinear interactions: away from resonance,
	they cancel the problematic terms, while near resonance,
	algebraic rearrangements avoid division by small phases
	and retain enough frequency decay to control the
	remainders.
	
	For every $s<1/2$, we also construct solutions with
	arbitrarily small initial data in
	$H^s\times L^2\times H^{-1}$ whose wave component blows
	up in finite time while the Schr\"odinger component
	remains bounded in $H^s$. This construction combines
	explicit approximate solutions supported on three
	Fourier modes with stability at $H^{1/2}\times L^2$.
	The difficulty is to combine infinitely many such
	solutions while controlling their nonlinear interactions. A successive choice
	of increasing frequencies makes the accumulated errors
	summable, while the wave amplitudes produce blow-up.
	Together, the two results identify $s=1/2$ as the
	threshold separating global well-posedness from the
	existence of such blow-up solutions, leaving the
	endpoint unresolved.
\end{abstract}
\tableofcontents

\section{Introduction and statement of the main results}
Fix $\gamma_1,\gamma_2>0$ and equip
$\T^2=(\R/2\pi\gamma_1\Z)\times(\R/2\pi\gamma_2\Z)$
with Lebesgue measure. We consider the Cauchy problem for
the two-dimensional periodic Zakharov system:
\begin{equation}\label{eq:zakharov}
	\left\{
	\begin{array}{l}
		i\partial_t u+\Delta u=nu,\\[2pt]
		\partial_t^2 n-\Delta n=\Delta |u|^2,\\[2pt]
		(u,n,\partial_t n)|_{t=0}=(u_0,n_0,n_1).
	\end{array}
	\right.
\end{equation}
Here $u$ is complex-valued and $n$ is real-valued.
The system was introduced by Zakharov \cite{Zakharov1972}
as a model for Langmuir waves in a plasma: $u$ describes
the slowly varying envelope of the electric field, while
$n$ describes the variation of the ion density.

Write $D=\sqrt{-\Delta}$. We assume that $n_1$ has zero
spatial mean, so that $\partial_t n$ has zero spatial
mean throughout the evolution. The operator $D^{-1}$
is defined on zero-mean distributions by the Fourier
multiplier $|k|^{-1}$ for $k\neq0$, with zero Fourier
coefficient at $k=0$. On the fixed torus, the zero-mean
subspace of $H^{-1}$ is precisely $DL^2$, and
\begin{equation*}
	\|D^{-1}f\|_{L^2}\sim\|f\|_{H^{-1}},
	\qquad \widehat f(0)=0.
\end{equation*}
For sufficiently regular solutions, \eqref{eq:zakharov}
conserves the Schr\"odinger mass
\begin{equation*}
	M(u(t))=M(u_0),
	\qquad
	M(f)=\int_{\T^2}|f(x)|^2\,dx,
\end{equation*}
and the Hamiltonian
\begin{equation*}
	\mathcal H(u,n,\partial_t n)
	=
	\|\nabla u\|_{L^2}^2
	+\frac12\|n\|_{L^2}^2
	+\frac12\|D^{-1}\partial_t n\|_{L^2}^2
	+\int_{\T^2}n|u|^2\,dx.
\end{equation*}
The associated energy space is
$H^1\times L^2\times H^{-1}$, with the zero-mean
condition on the third component.

The Cauchy problem on $\R^d$ has an extensive
well-posedness theory. Ozawa and Tsutsumi
\cite{OzawaTsutsumi1992} established local well-posedness
in $H^2\times H^1\times L^2$ for $d\leq3$.
Bourgain and Colliander \cite{BourgainColliander1996}
subsequently obtained local well-posedness at energy regularity
in dimensions two and three, together with global
results under suitable smallness assumptions.
Ginibre, Tsutsumi, and Velo \cite{GinibreTsutsumiVelo1997}
developed a systematic local theory in
$H^s\times H^\ell\times H^{\ell-1}$ in arbitrary
dimension. Further improvements relied on more precise
estimates for the interaction between the Schr\"odinger
and wave evolutions. In two dimensions, Bejenaru,
Herr, Holmer, and Tataru \cite{BHHT2009} reached
\begin{equation*}
	L^2(\R^2)\times H^{-1/2}(\R^2)\times H^{-3/2}(\R^2).
\end{equation*}
In three dimensions, Bejenaru and Herr
\cite{BejenaruHerr2011} obtained local well-posedness
in the full subcritical range on the line
$\ell=s-1/2$, namely $s>0$.
Chen and Wu \cite{ChenWu2021} and Sanwal
\cite{Sanwal2022} further enlarged the admissible
Sobolev ranges in dimensions two and three.
For $d\geq4$, Candy, Herr, and Nakanishi
\cite{CandyHerrNakanishi2023} determined the sharp
Sobolev range for local well-posedness with
real-analytic dependence on the initial data.

Global existence at low regularity depends strongly
on the dimension. On $\R$, Colliander, Holmer, and
Tzirakis \cite{CollianderHolmerTzirakis2008} proved
global well-posedness for arbitrary data in
$L^2\times H^{-1/2}\times H^{-3/2}$ by combining mass
conservation with control of the wave component.
In three dimensions, Hani, Pusateri, and Shatah
\cite{HaniPusateriShatah2013} established global
existence and scattering for small localized data.
For radial data, Guo, Nakanishi, and Wang
\cite{GuoNakanishiWang2013} proved a scattering and
grow-up dichotomy below the ground-state energy.
Kato and Kinoshita \cite{KatoKinoshita2023}
established small-data global well-posedness and
scattering at the critical regularities
$(s,\ell)=(0,-1/2)$ under radial symmetry.
More recently, Guo and Nakanishi
\cite{GuoNakanishi2026}, in a preprint, classified
the radial dynamics at energies slightly above the
ground-state energy, distinguishing scattering,
grow-up, and solutions trapped near the ground state.

The four-dimensional problem is energy-critical.
Bejenaru, Guo, Herr, and Nakanishi \cite{BGHN2015}
proved small-data global well-posedness and scattering
in the energy space. Kato and Tsugawa
\cite{KatoTsugawa2017} obtained small-data global
well-posedness and scattering at critical regularity
in dimensions $d\geq4$.
Guo and Nakanishi \cite{GuoNakanishi2021}
established the radial dynamics below the
four-dimensional ground state, while Candy, Herr,
and Nakanishi \cite{CandyHerrNakanishi2021}
proved global well-posedness in the corresponding
potential well without radial symmetry.
Following the concentration-compactness reduction
of Candy \cite{Candy2024}, the recent preprint of
Candy and Nakanishi \cite{CandyNakanishi2026}
establishes scattering throughout this potential
well without a symmetry assumption.
In the opposite direction, Krieger and Schmid
\cite{KriegerSchmid2024I,KriegerSchmid2024II}
constructed finite-time type II blow-up solutions
for the four-dimensional system, subject to the
numerical nondegeneracy assumptions stated in
their work.

Ill-posedness results provide complementary
restrictions on the local theory. These restrictions
must be distinguished according to the property that
fails: failure of twice differentiability of the
solution map, norm inflation, and finite-time blow-up
are different conclusions.
For the Euclidean problem, Grube \cite{Grube2025}
proved norm inflation in several Sobolev ranges
in arbitrary dimension, establishing the sharpness,
up to endpoints, of a number of local
well-posedness results.

For the periodic case, the available
dispersive estimates and the geometry of the
frequency lattice lead to a different regularity
theory. Bourgain \cite{Bourgain1994} initiated the
study of the periodic Cauchy problem and invariant
measures. Takaoka \cite{Takaoka1999} obtained
low-regularity local well-posedness in one dimension
and identified the dependence of the regularity
threshold on the period.
For the one-dimensional periodic problem, Compaan
\cite{Compaan2019} combined the $I$-method with a
high-low argument to obtain global existence below
energy regularity, including Schr\"odinger data in
$H^{1/2+\varepsilon}$ and wave data in $L^2$.

In two dimensions, Kishimoto
\cite[Theorem 1.1]{KishimotoLWP} established local
well-posedness on arbitrary rectangular tori in
$H^s\times H^\ell\times H^{\ell-1}$ whenever
\begin{equation*}
	0\leq s-\ell\leq1,
	\qquad
	2s\geq\ell+1\geq1.
\end{equation*}
In particular, for $\ell=0$ this includes
\begin{equation*}
	H^s(\T^2)\times L^2(\T^2)\times H^{-1}(\T^2),
	\qquad \frac12\leq s\leq1.
\end{equation*}
The endpoint $s=1/2$ is therefore included in the
local theory. The counterexamples in
\cite[Theorem 1.2]{KishimotoLWP} show that, with
$\ell=0$, the solution map cannot be twice
differentiable at the origin when $s<1/2$.
This gives a precise obstruction to a smooth
iteration theory below that threshold.
Unconditional uniqueness and related questions
for periodic Zakharov systems were subsequently
studied in \cite{Kishimoto2022}.

In higher dimensions, Kinoshita, Nakamura, and
Sanwal \cite{KNS2023} improved the periodic local
theory by establishing a decoupling inequality for
the paraboloid under a shell restriction and a
trilinear estimate involving the paraboloid and
the cone. On the line $\ell=s-1/2$, their result
gives local well-posedness for
\begin{equation*}
	s>\frac12\quad(d=3),\qquad
	s>\frac34\quad(d=4),\qquad
	s>\frac{d-3}{2}\quad(d\geq5).
\end{equation*}
These conclusions also hold on rectangular tori.
For $d=3$ and $d\geq5$, the thresholds are sharp
up to the endpoint within the iteration framework.
The estimates in \cite{KNS2023} illustrate the
importance of exploiting the geometry of the
interacting frequencies in the periodic problem.

We now return to global existence in two dimensions.
Let $Q$ be the positive radial ground state satisfying
\begin{equation}\label{eq:groundstate}
	-\Delta Q+Q-Q^3=0
	\quad\hbox{on }\R^2,
	\qquad
	M(Q)=\int_{\R^2}|Q(x)|^2\,dx.
\end{equation}
The number $M(Q)$ is the natural mass threshold
for the energy-space theory. The sharp
Gagliardo--Nirenberg inequality, together with mass
and energy conservation, controls the energy norm
under the strict condition $M(u_0)<M(Q)$.
On $\R^2$, Glangetas and Merle
\cite{GlangetasMerle1994I,GlangetasMerle1994II}
constructed finite-time blow-up solutions with
Schr\"odinger mass arbitrarily close to $M(Q)$
from above and studied concentration at blow-up.
For the periodic problem, Kishimoto and Maeda
\cite{KishimotoMaeda} constructed finite-time
blow-up solutions with mass arbitrarily close to
the same threshold from above. They also proved
global existence in the energy space when
$M(u_0)\leq M(Q)$, including the equality case,
which requires more than the coercive estimate
available under a strict mass gap.

Below energy regularity, the Hamiltonian need not
be finite, and conservation of mass alone does not
control the Sobolev norms required to continue the
local solution. The $I$-method, introduced by
Colliander, Keel, Staffilani, Takaoka, and Tao
\cite{CKSTT2002}, addresses this difficulty through
an almost conserved energy.
For the two-dimensional Euclidean Zakharov system,
Fang, Pecher, and Zhong \cite{FangPecherZhong2009}
used this method to obtain global well-posedness
for $3/4<s<1$ with wave data at energy regularity,
under a smallness condition on the Schr\"odinger
mass. Kishimoto \cite{KishimotoI} refined the
modified energy by separating resonant and
nonresonant interactions. In particular,
\cite[Theorem 1.1]{KishimotoI} proves global
well-posedness on every rectangular torus for
\begin{equation*}
	(u_0,n_0,n_1)\in
	H^s\times L^2\times H^{-1},
	\qquad
	\frac9{14}<s<1,
	\qquad
	M(u_0)<M(Q).
\end{equation*}
The same work also improves the Euclidean global
theory and allows negative Sobolev regularity
for the wave component.

The purpose of the present work is to extend the
two-dimensional periodic global theory to every
\begin{equation*}
	\frac12<s<1
\end{equation*}
under the same ground-state mass condition.
We retain the modified local theory and the
resonant cubic correction of \cite{KishimotoI},
and introduce additional quartic and quintic
corrections to the energy. Their coefficients
are chosen to cancel the interactions responsible
for the more restrictive regularity requirement
in the previous almost conservation argument.
The resulting energy increments can be summed
over successive local intervals for every
$s>1/2$. This yields global well-posedness
arbitrarily close to the endpoint of the local
theory, without a smallness assumption on the
wave data or on the $H^s$ norm of $u_0$.
We state the result as follows.

\begin{theorem}\label{thm:main}
	Let $1/2<s<1$. Suppose that $u_0\in H^s(\T^2)$,
	$n_0\in L^2(\T^2;\R)$ and $n_1\in H^{-1}(\T^2;\R)$ satisfy
	\begin{equation}\label{eq:massassumption}
		M(u_0)<M(Q).
	\end{equation}
	Then \eqref{eq:zakharov}, with data $(u_0,n_0,n_1)$, has a global
	solution in $C(\R;H^s\times L^2\times H^{-1})$. The solution is unique
	in the space specified in Lemma~\ref{lem:local},
	and depends continuously on the data in these Sobolev spaces.
	For every $\eps>0$ there is a constant $C_\eps$ such that
	\begin{equation}\label{eq:growth}
		\sup_{|t|\le T}\bigl(\norm{u(t)}{H^s}+\norm{n(t)}2
		+\norm{D^{-1}n_t(t)}2\bigr)
		\le C_\eps(1+T)^{(1-s)/(2s-1)+\eps},\qquad T\ge0.
	\end{equation}
	The constant may depend on the data, $s$, $\gamma$ and the strict
	gap in \eqref{eq:massassumption}.
\end{theorem}

We next consider regularities below $s=1/2$.
The following theorem shows that arbitrarily small initial
data in $H^s\times L^2\times H^{-1}$ can generate solutions
whose wave component blows up in finite time, while the
Schr\"odinger component remains bounded in $H^s$.
Thus the regularity assumption in Theorem~\ref{thm:main}
cannot be lowered below $s=1/2$ to obtain global
well-posedness for all small data in these spaces.

\begin{theorem}\label{thm:blowup:main}
	Let $s<1/2$, $T>0$, and $\eta>0$.
	There exist initial data
	\begin{equation*}
		(u_0,n_0,n_1)
		\in
		H^s(\T^2)
		\times L^2(\T^2)
		\times H^{-1}(\T^2)
	\end{equation*}
	with $n_0,n_1$ real-valued,
	\begin{equation*}
		\widehat n_0(0)=\widehat n_1(0)=0,
		\qquad
		u_0\notin H^{1/2}(\T^2),
	\end{equation*}
	and
	\begin{equation}\label{eq:blowup:smallness}
		\|u_0\|_{H^s}
		+\|n_0\|_{L^2}
		+\|n_1\|_{H^{-1}}
		<\eta,
		\qquad
		0<\|u_0\|_{L^2}<\eta,
	\end{equation}
	for which the Zakharov system has a distributional solution
	$(u,n)$ on $[0,T)$ with these initial data.
	
	Moreover, for any prescribed
	\begin{equation}\label{eq:blowup:s0}
		\max(s,0)<s_0<1/2,
	\end{equation}
	the construction can be carried out so that
	\begin{equation}\label{eq:blowup:regularity}
		u\in C([0,T];H^{s_0}(\T^2)),
		\qquad
		n\in C([0,T);L^2(\T^2))
		\cap C^1([0,T);H^{-1}(\T^2)).
	\end{equation}
	The Schr\"odinger mass is conserved:
	\begin{equation}\label{eq:blowup:mass}
		\|u(t)\|_{L^2}=\|u_0\|_{L^2},
		\qquad 0\leq t<T.
	\end{equation}
	There exist constants $c,C>0$ and $t_*<T$ such that
	\begin{equation}\label{eq:blowup:n-rate}
		\frac{c}{T-t}
		\leq \|n(t)\|_{L^2}
		\leq \frac{C}{T-t},
		\qquad t_*<t<T,
	\end{equation}
	and
	\begin{equation}\label{eq:blowup:nt-rate}
		\frac{c}{T-t}
		\leq \|\partial_t n(t)\|_{H^{-1}}
		\leq \frac{C}{T-t},
		\qquad t_*<t<T.
	\end{equation}
	
\end{theorem}

We prove Theorem~\ref{thm:blowup:main} by constructing smooth
solutions successively. At each step, we add an explicit approximate solution
whose $u$ component has Fourier support at two frequencies
and whose $n_+$ component has Fourier support at one frequency. We choose these frequencies so that
the error tends to zero in the local solution spaces, while
the $H^{s_0}$ norm of the added Schr\"odinger component also
tends to zero. A perturbation argument then gives an exact
smooth solution with the prescribed new initial data.

We first establish the estimates needed for this perturbation
argument. We then construct the frequencies on
$\T^2$ and verify the approximate solution
directly. Finally, we choose the successive errors to be
summable and prove
\eqref{eq:blowup:regularity}--\eqref{eq:blowup:nt-rate}
for the limit.

\vspace{10pt}
\noindent
\textbf{Disclosure of AI assistance.}
The author used ChatGPT-6 Astra during the preparation and
revision of this paper to assist with developing partial
arguments (especially the division of frequency space into regions used to define
the correction terms  and the construction of blow-up), checking calculations, and revising the exposition
and \LaTeX{} code.
The author takes full responsibility for the mathematical
claims, proofs, and references presented in the paper.

\vspace{10pt}
\noindent
\textbf{Organization.}
The paper consists of two parts.
Part~I proves the global well-posedness result in
Theorem~\ref{thm:main}.
We first recall the local theory and the estimates used
in the modified-energy argument.
Sections~\ref{sec:quartic} and~\ref{sec:quintic} construct
the quartic and quintic corrections and estimate the
corresponding remainders.
We then establish fixed-time bounds for the corrections.
The resulting almost conservation law is combined with
the local theory to prove global well-posedness and
polynomial bounds for the Sobolev norms.
Part~II proves the blow-up result in
Theorem~\ref{thm:blowup:main}.
We first establish the Duhamel estimates
and stability result needed for the construction.
We then select suitable lattice frequencies, construct
explicit approximate solutions, and estimate their
interactions with a fixed smooth solution.
Finally, we construct a sequence of smooth exact
solutions, pass to a limit, and verify the finite-time blow-up of the wave component.

\vspace{10pt}
\noindent
\textbf{Notation.}
\begin{itemize}
	\item By $u\in C([0,T];B)$ (or $L^p([0,T];B)$) for a
	Banach space $B$, we mean that $u$ is a continuous
	(or Bochner $L^p$) map from $[0,T]$ to $B$;
	see \cite[page 301]{evans2022partial}.
	
	\item By $A\lesssim B$ (resp. $A\sim B$), we mean that
	there exists a positive constant $C$ such that $A\le CB$
	(resp. $C^{-1}B\le A\le CB$).
	If the constant $C$ depends on $p$, then we write
	$A\lesssim_p B$ (resp. $A\sim_p B$).
	Dependence on the fixed torus and the fixed regularity
	exponents is suppressed.
	
	\item By $A\ll B$, we mean that $A/B$ is sufficiently small.
	
	\item By $a+$, we mean $a+\varepsilon$ for an arbitrarily
	small $\varepsilon>0$. An estimate with exponent $a+$
	holds for every $\varepsilon>0$, with an implicit constant
	that may depend on $\varepsilon$.
	The notation $a-$ is defined similarly.
	
	\item We write $\langle k\rangle=(1+|k|^2)^{1/2}$.
	Spatial Fourier frequencies belong to
	$\gamma_1^{-1}\Z\times\gamma_2^{-1}\Z$.
	
	\item Dyadic spatial frequencies are denoted by $N,N_j$,
	and dyadic modulations by $L,L_j$.
	All dyadic sums range over $1,2,4,\ldots$.
	At the lowest spatial scale, the condition
	$|k|\sim1$ is understood as $|k|\lesssim1$.
\end{itemize}

\vspace{15pt}
\part{Global well-posedness for $s>1/2$}\label{qwe}
\section{Function spaces and preliminary estimates}

\subsection{Fourier conventions and the first-order equations}

The dual lattice is
$\gamma_1^{-1}\Z\times\gamma_2^{-1}\Z$. We set
\begin{equation}\label{eq:fourier}
	\widehat f(k)=\abs{\T^2}^{-1}\int_{\T^2}f(x)e^{-ik\cdot x}\,dx,
	\qquad \norm f2^2=\abs{\T^2}\sum_k|\widehat f(k)|^2.
\end{equation}
In particular, products have the usual discrete convolution formula
without an additional factor. Constants may depend on $\gamma$. For convenience, we only consider the case $|\T^2|=1$ in Part \ref{qwe}.

The mean of $n$ is constant because $n_1$ has mean zero. If it equals
$\mu$, the substitution $u=e^{-i\mu t}\widetilde u$ and
$n=\widetilde n+\mu$ reduces \eqref{eq:zakharov} to the zero-mean
case. We make that reduction during the proof. It does not change
the mass or the Sobolev norm of the Schr\"odinger component.
Put 
\begin{equation*}
	n_+=n+iD^{-1}n_t,\qquad n_-=\overline{n_+}.
\end{equation*}
Then $n=(n_++n_-)/2$ and the equation (\ref{eq:zakharov}) is transformed into the first order system:
\begin{equation}\label{eq:firstorder}
	iu_t+\Delta u=\frac12(n_++n_-)u,\qquad
	i\partial_tn_\pm=\pm Dn_\pm\pm D|u|^2.
\end{equation}
Smooth solutions conserve their mass and the Hamiltonian
\begin{equation}\label{eq:hamiltonian}
	\mathcal H(u,n_+)=\norm{\nabla u}2^2+\frac12\norm{n_+}2^2
	+\frac12\int_{\T^2}(n_++n_-)|u|^2\,dx.
\end{equation}
All later differentiation is performed first for smooth solutions.
The passage to rough local solutions is justified in
Remark~\ref{rem:extension}.

\vspace{8pt}
\subsection{Bourgain norms and time cutoffs}
Let $P_N$ denote a smooth dyadic spatial projection, and let
$Q_L^S,Q_L^{W_\pm}$ denote smooth dyadic projections to
$\langle\tau+|k|^2\rangle\sim L$ and
$\langle\tau\pm|k|\rangle\sim L$, respectively. All dyadic
indices below range over $1,2,4,\ldots$. We use
\begin{equation}\label{eq:resolution}
	\norm f{X_*^{r,b,1}}
	=\left(\sum_N N^{2r}
	\left(\sum_L L^b\norm{P_NQ_L^*f}{L^2_{t,x}}\right)^2
	\right)^{1/2},\qquad *=S,W_\pm.
\end{equation}
The norm on an interval is the infimum over full-line extensions.
Recall the classical embedding
\begin{equation}\label{time uniform}
	X_*^{s,1/2,1}(J)\hookrightarrow C(J;H^s),
	\qquad *=S,W_\pm.
\end{equation}
This follows from Fourier inversion in time,
Cauchy--Schwarz on each modulation interval, and
the $\ell^1$ summation over dyadic modulations.

Here and below $\widetilde f(\tau,k)$ denotes the spacetime Fourier transform.

\begin{lemma}\label{lem:inputs}
	Let $J=[t_0,t_0+h]$, where $0<h\le1$, and let $r\in\mathbb R$.
	Fix $\psi\in C_c^\infty(\mathbb R)$ such that $\psi=1$ on
	$[-1,1]$, and define
	\begin{equation*}
		\psi_h(t)=\psi\!\left(\frac{t-t_0}{h}\right).
	\end{equation*}
	Thus $\psi_h=1$ on $J$. For $0<b\le1/2$ and
	$*=S,W_\pm$, we have
	\begin{equation}\label{eq:smoothcutoff}
		\norm{\psi_h f}{X_*^{r,b,1}}
		\lesssim_b h^{1/2-b}\norm f{X_*^{r,1/2,1}}.
	\end{equation}
	The characteristic function $1_J$ satisfies the corresponding
	estimate with the strict restriction $0<b<1/2$:
	\begin{equation}\label{eq:intervalcutoff}
		\norm{1_J f}{X_*^{r,b,1}}
		\lesssim_b h^{1/2-b}\norm f{X_*^{r,1/2,1}}.
	\end{equation}
	
	If $\widetilde f$ is supported where $|k|\sim N_1$ and
	$\langle\tau\pm|k|\rangle\sim L_1$, then
	\begin{equation}\label{eq:sourcewave}
		\norm f{L^4_{t,x}}
		\lesssim N_1^{3/8}L_1^{3/8}\norm f{L^2_{t,x}}.
	\end{equation}
	If $\supp\widetilde u\subset\{|k|\sim N_1\}$ and
	$\supp\widetilde n\subset\{|k|\sim N_2\}$, then, for small
	$\eps>0$,
	\begin{equation}\label{eq:sourcepair}
		\norm{un}{L^2_{t,x}}\lesssim_\eps
		\norm u{X_S^{2\eps,1/2,1}}
		\norm n{X_{W_\pm}^{0,1/2-\eps,1}}
		+\norm u{X_S^{1/2+\eps,1/2,1}}
		\norm n{X_{W_\pm}^{0,0,1}}.
	\end{equation}
	Both choices of sign are allowed, and the constants are
	independent of the dyadic spatial and modulation frequencies.
\end{lemma}

Estimate \eqref{eq:smoothcutoff} is
\cite[Lemma 2.3]{KishimotoI}.
Estimate \eqref{eq:intervalcutoff} is explicitly stated in
\cite[(4.5)]{KishimotoI} for
$\chi_h=1_{[0,h]}$, with $0<b<1/2$.
The estimates for $J=[t_0,t_0+h]$ follow by translating time:
\begin{equation*}
	\mathcal F_{t,x}[f(t+t_0,x)](\tau,k)
	=e^{it_0\tau}\widetilde f(\tau,k),
	\qquad
	\norm{f(t+t_0,x)}{X_*^{r,b,1}}
	=\norm f{X_*^{r,b,1}}.
\end{equation*}
Finally, \eqref{eq:sourcewave} and \eqref{eq:sourcepair}
are \cite[Lemma 2.11]{KishimotoI} and
\cite[Lemma 4.2]{KishimotoI}, respectively.

\vspace{10pt}
\section{Modified energy and resonant decomposition}

\subsection{The multiplier \texorpdfstring{\(I\)}{I}}

Fix $1/2<s<1$. For $N\geq1$, define the Fourier multiplier
operator $I$ by
\begin{equation*}
	\widehat{If}(k)=m(k)\widehat f(k),
	\qquad
	m(k)=
	\left\{
	\begin{array}{ll}
		1, & |k|\leq N,\\[2pt]
		\displaystyle\left(\frac{N}{|k|}\right)^{1-s},
		& |k|\geq2N.
	\end{array}
	\right.
\end{equation*}
For $N<|k|<2N$, choose a smooth radial interpolation
so that $m$ is non-increasing in $|k|$, $0<m\leq1$,
and $g(k)=|k|^2(1-m(k)^2)$ satisfies
\begin{equation}\label{eq:gproperties}
	0\leq g(k)\leq |k|^2,
	\qquad
	|g(p)-g(l)|
	\lesssim \bigl||p|^2-|l|^2\bigr|.
\end{equation}

\vspace{8pt}
\subsection{The first modified energy and the resonant decomposition}

We first display the quantity whose coercivity will control the solution:
\begin{equation}\label{eq:firstmodified}
	\mathcal H(Iu,n_+)=\norm{\nabla Iu}2^2+\frac12\norm{n_+}2^2
	+\frac12\int_{\T^2}(n_++n_-)|Iu|^2\,dx.
\end{equation}
The operator $I$ acts only on the Schr\"odinger component.  Directly from the definition of $m$,
\begin{equation}\label{eq:I-map}
	\norm u{H^s}\lesssim\norm{Iu}{H^1}
	\lesssim N^{1-s}\norm u{H^s},\qquad
	\norm{Iu}2\leq\norm u2.
\end{equation}
The two inequalities
in \eqref{eq:I-map} will be used to choose the initial energy scale and
recover the original Sobolev norm in the global iteration.

We modify the cubic part of the energy (\ref{eq:firstmodified}) to cancel the cubic terms
in its time derivative. Indeed, we seek  multipliers
$\sigma_\pm(p,l)$ and define
\begin{equation}\label{eq:Htilde}
	\widetilde{\mathcal H}(u,n_+)
	=\norm{\nabla Iu}2^2+\frac12\norm{n_+}2^2
	+\frac12\sum_{p,l}\sum_\pm
	\sigma_\pm(p,l)\overline{\widehat u(p)}\widehat u(l)
	\widehat{n_\pm}(p-l).
\end{equation}

For smooth solutions, differentiating the quadratic part of
\eqref{eq:Htilde} gives the cubic contribution
\begin{equation}\label{eq:quadratic-cubic}
	-\frac{i}{2}\sum_{p,l}\sum_\pm
	\bigl(|p|^2m(p)^2-|l|^2m(l)^2\mp|p-l|\bigr)
	\overline{\widehat u(p)}\widehat u(l)
	\widehat{n_\pm}(p-l).
\end{equation}
Differentiating the cubic part also gives the cubic contribution
\begin{equation}\label{eq:cubic-linear}
	\frac{i}{2}\sum_{p,l}\sum_\pm
	\sigma_\pm(p,l)
	\bigl(|p|^2-|l|^2\mp|p-l|\bigr)
	\overline{\widehat u(p)}\widehat u(l)
	\widehat{n_\pm}(p-l).
\end{equation}
To cancel \eqref{eq:quadratic-cubic} against
\eqref{eq:cubic-linear}, we therefore choose
\begin{equation*}
	\sigma_\pm(p,l)=
	\frac{|p|^2m(p)^2-|l|^2m(l)^2\mp|p-l|}
	{|p|^2-|l|^2\mp|p-l|}
	\quad\text{if }\bigl||p|^2-|l|^2\bigr|>2|p-l|.
\end{equation*}
In the complementary (resonant) region, we instead define
\begin{equation*}
	\sigma_\pm(p,l)=
	\frac{|p|^2m(p)^2-|l|^2m(l)^2}{|p|^2-|l|^2}
	\quad\text{if }\bigl||p|^2-|l|^2\bigr|\leq2|p-l|.
\end{equation*}
For the subsequent correction terms, write
\begin{equation}\label{eq:b1-def}
	b_{1,\pm}(p,l)=1-\sigma_\pm(p,l).
\end{equation}
Using $g$ from \eqref{eq:gproperties}, we obtain
\begin{equation}\label{eq:b1non}
	b_{1,\pm}(p,l)
	=\frac{g(p)-g(l)}{|p|^2-|l|^2\mp|p-l|}
	\quad\text{if }\bigl||p|^2-|l|^2\bigr|>2|p-l|,
\end{equation}
and
\begin{equation}\label{eq:b1res}
	b_{1,\pm}(p,l)
	=\frac{g(p)-g(l)}{|p|^2-|l|^2}
	\quad\text{if }\bigl||p|^2-|l|^2\bigr|\leq2|p-l|.
\end{equation}
Thus, we have $|b_{1,\pm}(p,l)|\lesssim 1$.
By \eqref{eq:b1-def}, the sum of the coefficients in
\eqref{eq:quadratic-cubic} and \eqref{eq:cubic-linear}, after
removing the common factor $i/2$, is
\begin{equation*}
	g(p)-g(l)
	-\bigl(|p|^2-|l|^2\mp|p-l|\bigr)b_{1,\pm}(p,l).
\end{equation*}
By \eqref{eq:b1non}, this coefficient vanishes in the
nonresonant region. By \eqref{eq:b1res}, in the resonant
region it equals
\begin{equation*}
	\pm|p-l|\,b_{1,\pm}(p,l).
\end{equation*}
Consequently,
\begin{equation}\label{eq:b1residual}
	\left|
	g(p)-g(l)
	-\bigl(|p|^2-|l|^2\mp|p-l|\bigr)b_{1,\pm}(p,l)
	\right|
	\lesssim |p-l|.
\end{equation}

Based on the above choices, Kishimoto controlled the following fixed-time difference between $\widetilde{\mathcal H}(u,n_+)$ and $\mathcal H(Iu,n_+)$  (see \cite[Proposition 3.2]{KishimotoI}):
\begin{proposition}\label{prop:cubiccomparison}
	Let $1/2<s<1$, $u\in H^s(\T^2)$, and $n_+\in L^2(\T^2)$.
	Then, for every sufficiently small $\eps>0$,
	\begin{equation}\label{eq:cubiccomparison}
		|\widetilde{\mathcal H}(u,n_+)-\mathcal H(Iu,n_+)|
		\lesssim_\eps N^{-1+\eps}\norm{Iu}{H^1}^2\norm{n_+}2.
	\end{equation}
\end{proposition}

\vspace{8pt}
\subsection{The additional energy corrections}

Starting from \eqref{eq:firstmodified} and its modification
\eqref{eq:Htilde}, we will subtract a quartic term $\mathcal C_4$ and a quintic
term $\mathcal C_5$ from \eqref{eq:Htilde}. In fact, the energy used in the iteration is
\begin{equation}\label{eq:energydef}
	\mathcal E_N(u,n_+)
	=\widetilde{\mathcal H}(u,n_+)-\mathcal C_4(u,n_+)
	-\mathcal C_5(u,n_+),
\end{equation}
where 
\begin{equation}\label{eq:C4-definition}
	\mathcal C_4(u,n_+)
	=\frac14\sum_{p,r,l}\sum_{\pm_1,\pm_2}
	b_{2,\pm_1,\pm_2}(p,r,l)\,
	\overline{\widehat u(p)}\widehat u(l)
	\widehat{n_{\pm_1}}(p-r)
	\widehat{n_{\pm_2}}(r-l)
\end{equation}
and
\begin{equation}\label{eq:C5-definition}
	\mathcal C_5(u,n_+)
	=\frac18\sum_{p,r,z,l}\sum_{\pm_1,\pm_2,\pm_3}
	b_{3,\pm_1,\pm_2,\pm_3}(p,r,z,l)\,
	\overline{\widehat u(p)}\widehat u(l)
	\widehat{n_{\pm_1}}(p-r)
	\widehat{n_{\pm_2}}(r-z)
	\widehat{n_{\pm_3}}(z-l).
\end{equation}
The coefficients $b_2$ and $b_3$ in
\eqref{eq:C4-definition} and \eqref{eq:C5-definition}
are specified in Sections~\ref{sec:quartic}
and~\ref{sec:quintic}, respectively.

 The estimates for
the two additional terms and their derivatives are the new
parts of the argument. Proposition~\ref{prop:comparison} bounds
their fixed-time size, and Proposition~\ref{prop:increment} gives
the resulting almost conservation law. 

\vspace{10pt}
\section{Global solutions for the periodic case}

\subsection{Fixed-time estimate and almost conservation law}

\begin{proposition}\label{prop:comparison}
	Let $1/2<s<1$, and let $(u,n_+)$ be smooth functions.
	Set
	\begin{equation*}
		a=1+\norm{Iu}{H^1}+\norm{n_+}2.
	\end{equation*}
	Then the energy \eqref{eq:energydef} satisfies
	\begin{equation}\label{eq:comparison}
		|\mathcal E_N(u,n_+)-\mathcal H(Iu,n_+)|
		\lesssim_\eps
		N^{-1+\eps}a^3+N^{-2+\eps}a^4+N^{-3+\eps}a^5.
	\end{equation}
	The estimate holds for every sufficiently small $\eps>0$,
	with constants independent of $N,u,n_+$.
\end{proposition}

For $0<h\leq1$ and a fixed small $\eps>0$, put
\begin{equation}\label{eq:d3}
	d_3(N,h)=N^{-1+\eps}h^{1/2-\eps},
\end{equation}
\begin{equation}\label{eq:d4}
	d_4(N,h)=N^{-2+\eps}+N^{-1+\eps}h^{1-\eps},
\end{equation}
\begin{equation}\label{eq:d5}
	d_5(N,h)=N^{-3+\eps}+N^{-2+\eps}h^{1-\eps},
\end{equation}
and
\begin{equation}\label{eq:d6}
	d_6(N,h)=N^{-5/2+\eps}h^{1/2-\eps}.
\end{equation}

\begin{proposition}\label{prop:increment}
	Let $1/2<s<1$, and let $(u,n_+)$ be a smooth solution on
	$J=[t_0,t_1]$, where $0<h=t_1-t_0\leq1$. Set
	\begin{equation}\label{eq:A}
		A=1+\norm{Iu}{X_S^{1,1/2,1}(J)}
		+\norm{n_+}{X_{W_+}^{0,1/2,1}(J)}.
	\end{equation}
	Then the energy \eqref{eq:energydef} satisfies
	\begin{equation}\label{eq:increment}
		|\mathcal E_N(t_1)-\mathcal E_N(t_0)|
		\lesssim_\eps\sum_{j=3}^6d_j(N,h)A^j.
	\end{equation}
	The estimate holds for every sufficiently small $\eps>0$,
	with constants independent of $N,h$ and the solution.
\end{proposition}

 The subsequent sections
complete the proofs of Propositions~\ref{prop:comparison}
and~\ref{prop:increment}.

\vspace{8pt}
\subsection{Consequences of the preliminary estimates}

We now derive the estimates that will be used to bound the
quartic, quintic, and sixth-order expressions in the derivative
of the modified energy. We first estimate each spatial frequency
separately and retain the dependence on the spatial frequencies
so that the resulting sums can later be estimated by
Cauchy--Schwarz.

\paragraph{1. Choice of functions on the whole time line.}

The norms in \eqref{eq:A} are defined on $J$, whereas the estimates
in Lemma~\ref{lem:inputs} are stated for functions on
$\mathbb R\times\mathbb T^2$. By the definition of the norms on $J$,
we may choose functions on $\mathbb R\times\mathbb T^2$ that agree
with $u$ and $n_+$ on $J$ and satisfy
\begin{equation*}
	1+\norm{Iu}{X_S^{1,1/2,1}}
	+\norm{n_+}{X_{W_+}^{0,1/2,1}}
	\leq 2A.
\end{equation*}
We use the same letters for these functions. The equations are
used only on $J$; no equation is required outside $J$.
We also set $n_-=\overline{n_+}$ on the whole time line $\R$.

For dyadic spatial frequencies $N_1,N_2\geq1$, write
\begin{equation*}
	u_{N_1}=P_{N_1}u,\qquad
	n_{\pm,N_2}=P_{N_2}n_\pm,
\end{equation*}
and set
\begin{equation*}
	A_{N_1}=\norm{Iu_{N_1}}{X_S^{1,1/2,1}},
	\qquad
	B_{N_2}=\norm{n_{+,N_2}}{X_{W_+}^{0,1/2,1}}.
\end{equation*}
Consequently, we have
\begin{equation*}
	\left(\sum_{N_1}A_{N_1}^2\right)^{1/2}
	+\left(\sum_{N_2}B_{N_2}^2\right)^{1/2}
	\lesssim A.
\end{equation*}

Conjugation changes $(\tau,k)$ into $(-\tau,-k)$ and changes
the modulation $\langle\tau+|k|\rangle$ into
$\langle\tau-|k|\rangle$. Consequently,
\begin{equation*}
	\norm{n_{-,N_2}}{X_{W_-}^{0,1/2,1}}
	=\norm{n_{+,N_2}}{X_{W_+}^{0,1/2,1}}
	=B_{N_2}.
\end{equation*}
Thus the same quantity $B_{N_2}$ controls either choice of sign.

\paragraph{2. Estimates at each fixed time.}

For any of the modulation projections
$Q_{L_1}^*$, Fourier inversion in time and Minkowski's inequality give
\begin{equation*}
	\norm{P_{N_1}Q_{L_1}^*f}{L_t^\infty L_x^2}
	\lesssim
	\left(
	\sum_k
	\left(
	\int_{\mathbb R}
	\left|\widetilde{P_{N_1}Q_{L_1}^*f}(\tau,k)\right|\,d\tau
	\right)^2
	\right)^{1/2}.
\end{equation*}
For each fixed $k$, the integrand is supported on a set of
$\tau$ of measure $O(L_1)$. Cauchy--Schwarz in $\tau$, followed
by Plancherel's identity, therefore yields
\begin{equation*}
	\norm{P_{N_1}Q_{L_1}^*f}{L_t^\infty L_x^2}
	\lesssim L_1^{1/2}\norm{P_{N_1}Q_{L_1}^*f}{L^2_{t,x}}.
\end{equation*}
Summing in $L_1$ proves
\begin{equation*}
	\norm{P_{N_1}f}{L_t^\infty L_x^2}
	\lesssim
	\sum_{L_1}L_1^{1/2}\norm{P_{N_1}Q_{L_1}^*f}{L^2_{t,x}}.
\end{equation*}

On the support of $P_{N_1}$, we have
$m(k)\sim m(N_1)$ and $\langle k\rangle \sim N_1$.
Applying the preceding estimate to $Iu_{N_1}$ gives
\begin{equation*}
	\norm{u_{N_1}}{L_t^\infty L_x^2}
	\lesssim m(N_1)^{-1}\norm{Iu_{N_1}}{L_t^\infty L_x^2}
	\lesssim N_1^{-1}m(N_1)^{-1}A_{N_1}.
\end{equation*}
The two-dimensional Bernstein inequality gives
\begin{equation*}
	\norm{u_{N_1}}{L_{t,x}^\infty}
	\lesssim N_1\norm{u_{N_1}}{L_t^\infty L_x^2}.
\end{equation*}
We have therefore proved
\begin{equation}\label{eq:spatialu}
	\norm{u_{N_1}}{L_t^\infty L_x^2}
	\lesssim N_1^{-1}m(N_1)^{-1}A_{N_1},
	\qquad \norm{u_{N_1}}{L_{t,x}^\infty}
	\lesssim m(N_1)^{-1}A_{N_1}.
\end{equation}

Applying the same estimate to $n_{\pm,N_2}$ gives
\begin{equation}\label{eq:spatialwave}
	\norm{n_{\pm,N_2}}{L_t^\infty L_x^2}\lesssim B_{N_2},
	\qquad \norm{n_{\pm,N_2}}{L_{t,x}^\infty}\lesssim N_2 B_{N_2}.
\end{equation}
The estimates \eqref{eq:spatialu} and \eqref{eq:spatialwave}
will be used when spatial H\"older inequalities suffice.

\paragraph{3. An estimate for the product of $u_{N_1}$ and
	$n_{\pm,N_2}$.}

Choose the smooth function $\psi_h$ from
Lemma~\ref{lem:inputs}, so that $\psi_h=1$ on $J$.
Then
\begin{equation*}
	(\psi_hu_{N_1})(1_Jn_{\pm,N_2})=1_Ju_{N_1}n_{\pm,N_2}.
\end{equation*}

Apply \eqref{eq:sourcepair} with $\eps$ replaced by $\eps/2$. We obtain
\begin{equation}\label{eq:pair-before-cutoffs}
	\begin{split}
		\norm{(\psi_hu_{N_1})(1_Jn_{\pm,N_2})}{L^2_{t,x}}
		\lesssim_\eps{}&
		\norm{\psi_hu_{N_1}}{X_S^{\eps,1/2,1}}
		\norm{1_Jn_{\pm,N_2}}{X_{W_\pm}^{0,1/2-\eps/2,1}}\\
		&+
		\norm{\psi_hu_{N_1}}{X_S^{1/2+\eps/2,1/2,1}}
		\norm{1_Jn_{\pm,N_2}}{X_{W_\pm}^{0,0,1}}.
	\end{split}
\end{equation}
Spatial frequency localization and \eqref{eq:smoothcutoff}
at $b=1/2$ imply, for either spatial exponent used above,
\begin{equation*}
	\norm{\psi_hu_{N_1}}{X_S^{r,1/2,1}}
	\lesssim
	N_1^{r-1}m(N_1)^{-1}
	\norm{\psi_hIu_{N_1}}{X_S^{1,1/2,1}}
	\lesssim N_1^{r-1}m(N_1)^{-1}A_{N_1}.
\end{equation*}
For the first occurrence of $n_{\pm,N_2}$ in
\eqref{eq:pair-before-cutoffs}, the estimate with $1_J$
at $b=1/2-\eps/2$ gives
\begin{equation*}
	\norm{1_Jn_{\pm,N_2}}{X_{W_\pm}^{0,1/2-\eps/2,1}}
	\lesssim_\eps h^{\eps/2}B_{N_2}.
\end{equation*}
For the second occurrence, we can similarly obtain
\begin{equation*}
	\norm{1_Jn_{\pm,N_2}}{X_{W_\pm}^{0,0,1}}
	\leq
	\norm{1_Jn_{\pm,N_2}}{X_{W_\pm}^{0,\eps/2,1}}
	\lesssim_\eps h^{1/2-\eps/2}B_{N_2}.
\end{equation*}
Substitution into \eqref{eq:pair-before-cutoffs} yields
\begin{equation}\label{eq:pair}
	\norm{(\psi_hu_{N_1})(1_Jn_{\pm,N_2})}{L^2_{t,x}}
	\lesssim_\eps m(N_1)^{-1}
	\left(
	N_1^{-1+\eps}
	+h^{1/2-\eps}N_1^{-1/2+\eps}
	\right)A_{N_1}B_{N_2}.
\end{equation}
The same estimate holds with $\overline{u_{N_1}}$ in place
of $u_{N_1}$: conjugating the whole product replaces
$\overline{u_{N_1}}n_{\pm,N_2}$ by $u_{N_1}n_{\mp,N_2}$ and leaves
its $L^2$ norm unchanged.

\paragraph{4. An $L^4$ estimate for $n_{\pm,N_2}$.}

Decomposing $1_Jn_{\pm,N_2}$ in modulation, \eqref{eq:sourcewave} gives
\begin{equation*}
		\norm{1_Jn_{\pm,N_2}}{L^4_{t,x}}
		\leq
		\sum_{L_1}
		\norm{Q_{L_1}^{W_\pm}(1_Jn_{\pm,N_2})}{L^4_{t,x}}\\
		\lesssim
		N_2^{3/8}\sum_{L_1}L_1^{3/8}
		\norm{Q_{L_1}^{W_\pm}(1_Jn_{\pm,N_2})}{L^2_{t,x}}\\
		\lesssim
		N_2^{3/8}
		\norm{1_Jn_{\pm,N_2}}{X_{W_\pm}^{0,3/8,1}}.
\end{equation*}
Since $0<3/8<1/2$, \eqref{eq:intervalcutoff} gives
\begin{equation*}
	\norm{1_Jn_{\pm,N_2}}{X_{W_\pm}^{0,3/8,1}}
	\lesssim h^{1/2-3/8}B_{N_2}=h^{1/8}B_{N_2}.
\end{equation*}
Thus, one can obtain the following $L^4$ estimate
\begin{equation}\label{eq:wavefour}
	\norm{1_Jn_{\pm,N_2}}{L^4_{t,x}}\lesssim N_2^{3/8}h^{1/8}B_{N_2}.
\end{equation}
 We will use \eqref{eq:wavefour} when H\"older's
inequality assigns exponent $4$ to $n_{\pm,N_2}$.

\begin{remark}\label{absolute}
	
	Let $f_1,\ldots,f_j$ denote the functions in one summand,
	and let $M(k_1,\ldots,k_j)$ be its spatial coefficient.
	At each fixed time, for $1\leq p_a\leq\infty$ with
	$\sum_{a=1}^j1/p_a=1$, we have
	\begin{equation*}
		\left|
		\sum_{k_1+\cdots+k_j=0}
		M(k_1,\ldots,k_j)
		\prod_{a=1}^j\widehat f_a(t,k_a)
		\right|
		\leq
		\sup|M|
		\sum_{k_1+\cdots+k_j=0}
		\prod_{a=1}^j|\widehat f_a(t,k_a)|
	\end{equation*}
	\begin{equation*}
		=
		\sup|M|
		\int_{\T^2}
		\prod_{a=1}^j
		\left(\mathcal F_x^{-1}
		|\widehat f_a(t,\cdot)|\right)(x)\,dx
		\lesssim
		\sup|M|
		\prod_{a=1}^j
		\left\|
		\mathcal F_x^{-1}|\widehat f_a(t,\cdot)|
		\right\|_{L_x^{p_a}}.
	\end{equation*}
	The resulting product of norms can then be integrated in
	time and estimated by H\"older's inequality.
	
	Taking absolute values of the spatial Fourier coefficients
	preserves their support and the spatial $L^2$ norm.
	Consequently, \eqref{eq:spatialu} and \eqref{eq:spatialwave}
	remain valid with $u_{N_1}$ and $n_{\pm,N_2}$ replaced by
	$\mathcal F_x^{-1}|\widehat{u_{N_1}}|$ and
	$\mathcal F_x^{-1}|\widehat{n_{\pm,N_2}}|$, respectively.
	Likewise, taking absolute values of the spacetime Fourier
	transforms preserves their supports and the
	$X_*^{r,b,1}$ norms. Thus \eqref{eq:pair} and
	\eqref{eq:wavefour} remain valid with $\psi_hu_{N_1}$
	and $1_Jn_{\pm,N_2}$ replaced by
	$\mathcal F_{t,x}^{-1}|\widetilde{\psi_hu_{N_1}}|$ and
	$\mathcal F_{t,x}^{-1}|\widetilde{1_Jn_{\pm,N_2}}|$,
	respectively.
\end{remark}

\vspace{8pt}
\subsection{Fourier expressions and their time differentiation}

We first record the equations (\ref{eq:firstorder}) under Fourier transform:
\begin{equation}\label{eq:fourier-time-derivatives}
	\begin{split}
		\partial_t\widehat u(p)
		&=-i|p|^2\widehat u(p)
		-\frac{i}{2}\sum_r\sum_\pm
		\widehat{n_\pm}(p-r)\widehat u(r),\\
		\partial_t\overline{\widehat u(p)}
		&=i|p|^2\overline{\widehat u(p)}
		+\frac{i}{2}\sum_r\sum_\pm
		\widehat{n_\pm}(r-p)\overline{\widehat u(r)},\\
		\partial_t\widehat{n_\pm}(q)
		&=\mp i|q|\widehat{n_\pm}(q)
		\mp i|q|\sum_k\widehat u(k)\overline{\widehat u(k-q)}.
	\end{split}
\end{equation}
The second identity follows from the first using $	\overline{\widehat{n_\pm}(q)}
=\widehat{n_\mp}(-q)$.

\paragraph{1. A common notation for the correction terms.}
We introduce a common notation for the quartic and quintic corrections defined earlier, so that their time
derivatives can be computed by a single formula.

For $1\leq j\leq4$, define
\begin{equation*}
	\Omega_{j,\pm_1,\ldots,\pm_j}(p_0,\ldots,p_j)
	=
	|p_0|^2-|p_j|^2
	-\sum_{a=1}^j\bigl(\pm_a|p_{a-1}-p_a|\bigr).
\end{equation*}
 For a coefficient
$c_{\pm_1,\ldots,\pm_j}(p_0,\ldots,p_j)$, define
\begin{equation}\label{eq:form}
		\mathcal C_{j+2}[c](u,n_+)
		=\frac{1}{2^j}
		\sum_{p_0,\ldots,p_j}
		\sum_{\pm_1,\ldots,\pm_j}
		c_{\pm_1,\ldots,\pm_j}(p_0,\ldots,p_j)
		\overline{\widehat u(p_0)}\widehat u(p_j)
		\prod_{a=1}^j
		\widehat{n_{\pm_a}}(p_{a-1}-p_a).
\end{equation}
With the notation of \eqref{eq:form}, the quartic and quintic
corrections defined in \eqref{eq:C4-definition} and
\eqref{eq:C5-definition} can be written as
\begin{equation*}
	\mathcal C_4=\mathcal C_4[b_2],
	\qquad
	\mathcal C_5=\mathcal C_5[b_3].
\end{equation*}
Likewise, by \eqref{eq:b1-def}, the cubic part of
\eqref{eq:Htilde} is
\begin{equation*}
	\mathcal C_3[1-b_1]
	=\mathcal C_3[1]-\mathcal C_3[b_1].
\end{equation*}
We write $\mathcal C_3=\mathcal C_3[b_1]$.

The coefficients $b_1,b_2,b_3$ appearing in the preceding
definitions are real and satisfy
\begin{equation}\label{eq:coefficient-reversal}
	b_{j,\pm_1,\ldots,\pm_j}(p_0,\ldots,p_j)
	=
	b_{j,\mp_j,\ldots,\mp_1}(p_j,\ldots,p_0),
	\qquad 1\leq j\leq3.
\end{equation}
For $b_1$, this follows directly from
\eqref{eq:b1non} and \eqref{eq:b1res}.
For $b_2$ and $b_3$, the identity follows from their
definitions in Sections~\ref{sec:quartic} and
\ref{sec:quintic}.

By \eqref{eq:coefficient-reversal} and
$\overline{\widehat{n_\pm}(q)}=\widehat{n_\mp}(-q)$,
complex conjugation followed by reversing the frequencies
and signs gives
$\overline{\mathcal C_{j+2}[b_j]}=\mathcal C_{j+2}[b_j]$.
Thus the correction terms are real-valued.

To write the time derivative uniformly for
$1\leq j\leq4$, set
\begin{equation*}
	b_0(p)=g(p),\qquad b_4=0.
\end{equation*}
These conventions do not change $b_1,b_2,b_3$.
For $1\leq j\leq4$, define
\begin{equation*}
	R_{j,\pm_1,\ldots,\pm_j}(p_0,\ldots,p_j)
	=
	\Omega_{j,\pm_1,\ldots,\pm_j}(p_0,\ldots,p_j)
	b_{j,\pm_1,\ldots,\pm_j}(p_0,\ldots,p_j)
\end{equation*}
\begin{equation}\label{eq:Rj}
	+b_{j-1,\pm_2,\ldots,\pm_j}(p_1,\ldots,p_j)
	-b_{j-1,\pm_1,\ldots,\pm_{j-1}}(p_0,\ldots,p_{j-1}).
\end{equation}
For $j=1$ and $j=4$, respectively, this reads
\begin{equation*}
	R_{1,\pm}(p_0,p_1)
	=
	\bigl(|p_0|^2-|p_1|^2\mp|p_0-p_1|\bigr)b_{1,\pm}(p_0,p_1)
	+g(p_1)-g(p_0)
\end{equation*}
and
\begin{equation*}
		R_{4,\pm_1,\pm_2,\pm_3,\pm_4}(p_0,p_1,p_2,p_3,p_4)
		=b_{3,\pm_2,\pm_3,\pm_4}(p_1,p_2,p_3,p_4)
		-b_{3,\pm_1,\pm_2,\pm_3}(p_0,p_1,p_2,p_3).
\end{equation*}
Using \eqref{eq:form}, we write
\begin{equation*}
	\mathcal R_j=\mathcal C_{j+2}[R_j],
	\qquad 1\leq j\leq4.
\end{equation*}

\paragraph{2. Differentiation of the correction terms.}

We claim that, for $1\leq j\leq3$,
\begin{equation}\label{eq:correction-derivative}
	\frac{d}{dt}\mathcal C_{j+2}[b_j]
	=
	i\mathcal C_{j+2}[\Omega_jb_j]
	+i\mathcal C_{j+3}[R_{j+1}-\Omega_{j+1}b_{j+1}]
	+\mathcal F_j,
\end{equation}
where 
\begin{equation}\label{eq:wave-replacement}
	\begin{split}
		\mathcal F_j
		={}&-\frac{i}{2^j}
		\sum_{p_0,\ldots,p_j}
		\sum_{\pm_1,\ldots,\pm_j}
		b_{j,\pm_1,\ldots,\pm_j}(p_0,\ldots,p_j)
		\overline{\widehat u(p_0)}\widehat u(p_j)\\
		&\quad{}\times
		\sum_{a=1}^j
		\bigl(\pm_a|p_{a-1}-p_a|\bigr)
		\left(
		\sum_k\widehat u(k)
		\overline{\widehat u(k-p_{a-1}+p_a)}
		\right)
		\prod_{\substack{1\leq d\leq j\\d\neq a}}
		\widehat{n_{\pm_d}}(p_{d-1}-p_d).
	\end{split}
\end{equation}

Differentiating a summand of \eqref{eq:form} in time and
substituting the linear parts of
\eqref{eq:fourier-time-derivatives}, we obtain
\begin{equation*}
	\begin{split}
		&\left(
		i|p_0|^2-i|p_j|^2
		-i\sum_{a=1}^j\bigl(\pm_a|p_{a-1}-p_a|\bigr)
		\right)
		\overline{\widehat u(p_0)}\widehat u(p_j)
		\prod_{a=1}^j\widehat{n_{\pm_a}}(p_{a-1}-p_a)\\
		&\qquad=
		i\Omega_{j,\pm_1,\ldots,\pm_j}(p_0,\ldots,p_j)
		\overline{\widehat u(p_0)}\widehat u(p_j)
		\prod_{a=1}^j\widehat{n_{\pm_a}}(p_{a-1}-p_a).
	\end{split}
\end{equation*}
After multiplication by $b_j/2^j$ and summation, this
equals $i\mathcal C_{j+2}[\Omega_jb_j]$.

The nonlinear part of $\overline{\hat{u}(p_0)}$ in \eqref{eq:fourier-time-derivatives}
contributes
\begin{equation*}
		\frac{i}{2^{j+1}}
		\sum_{p_0,\ldots,p_j,r}
		\sum_{\pm,\pm_1,\ldots,\pm_j}
		b_{j,\pm_1,\ldots,\pm_j}(p_0,\ldots,p_j)
		\widehat{n_\pm}(r-p_0)
		\overline{\widehat u(r)}\widehat u(p_j)
		\prod_{a=1}^j\widehat{n_{\pm_a}}(p_{a-1}-p_a).
\end{equation*}
Relabel
\begin{equation*}
	(r,p_0,\ldots,p_j)\longmapsto(p_0,p_1,\ldots,p_{j+1}),
	\qquad
	(\pm,\pm_1,\ldots,\pm_j)
	\longmapsto(\pm_1,\pm_2,\ldots,\pm_{j+1}).
\end{equation*}
The preceding contribution becomes
\begin{equation}\label{eq:left-u-derivative}
		\frac{i}{2^{j+1}}
		\sum_{p_0,\ldots,p_{j+1}}
		\sum_{\pm_1,\ldots,\pm_{j+1}}
		b_{j,\pm_2,\ldots,\pm_{j+1}}(p_1,\ldots,p_{j+1})
		\overline{\widehat u(p_0)}\widehat u(p_{j+1})
		\prod_{a=1}^{j+1}
		\widehat{n_{\pm_a}}(p_{a-1}-p_a).
\end{equation}

Similarly, the nonlinear part of $\widehat u(p_j)$ in \eqref{eq:fourier-time-derivatives} contributes
\begin{equation}\label{eq:right-u-derivative}
		-\frac{i}{2^{j+1}}
		\sum_{p_0,\ldots,p_{j+1}}
		\sum_{\pm_1,\ldots,\pm_{j+1}}
		b_{j,\pm_1,\ldots,\pm_j}(p_0,\ldots,p_j)
		\overline{\widehat u(p_0)}\widehat u(p_{j+1})
		\prod_{a=1}^{j+1}
		\widehat{n_{\pm_a}}(p_{a-1}-p_a).
\end{equation}
By \eqref{eq:Rj},
\begin{equation*}
	b_{j,\pm_2,\ldots,\pm_{j+1}}(p_1,\ldots,p_{j+1})
	-b_{j,\pm_1,\ldots,\pm_j}(p_0,\ldots,p_j)
\end{equation*}
\begin{equation*}
		=R_{j+1,\pm_1,\ldots,\pm_{j+1}}(p_0,\ldots,p_{j+1})-
		\Omega_{j+1,\pm_1,\ldots,\pm_{j+1}}(p_0,\ldots,p_{j+1})
		b_{j+1,\pm_1,\ldots,\pm_{j+1}}(p_0,\ldots,p_{j+1}).
\end{equation*}
Thus \eqref{eq:left-u-derivative} plus
\eqref{eq:right-u-derivative} equals
\begin{equation*}
	i\mathcal C_{j+3}[R_{j+1}-\Omega_{j+1}b_{j+1}].
\end{equation*}
It remains to compute the contribution from the nonlinear
parts of the equations for $n_\pm$. By
\eqref{eq:fourier-time-derivatives}, these parts are
\begin{equation*}
	\partial_t\widehat{n_{\pm_a}}(p_{a-1}-p_a)
	\;\longrightarrow\;
	-i\bigl(\pm_a|p_{a-1}-p_a|\bigr)
	\sum_k\widehat u(k)
	\overline{\widehat u(k-p_{a-1}+p_a)}.
\end{equation*}
The product rule gives
\begin{equation*}
	\partial_t
	\left(
	\prod_{a=1}^j
	\widehat{n_{\pm_a}}(p_{a-1}-p_a)
	\right)
	=
	\sum_{a=1}^j
	\left(
	\partial_t\widehat{n_{\pm_a}}(p_{a-1}-p_a)
	\right)
	\prod_{\substack{1\leq d\leq j\\d\neq a}}
	\widehat{n_{\pm_d}}(p_{d-1}-p_d).
\end{equation*}
Substituting the nonlinear parts displayed above, their
contribution to $\frac{d}{dt}\mathcal C_{j+2}[b_j]$ is
\begin{equation*}
	-\frac{i}{2^j}
	\sum_{p_0,\ldots,p_j}
	\sum_{\pm_1,\ldots,\pm_j}
	b_{j,\pm_1,\ldots,\pm_j}(p_0,\ldots,p_j)
	\overline{\widehat u(p_0)}\widehat u(p_j)
\end{equation*}
\begin{equation*}
	\qquad{}\times
	\sum_{a=1}^j
	\bigl(\pm_a|p_{a-1}-p_a|\bigr)
	\left(
	\sum_k\widehat u(k)
	\overline{\widehat u(k-p_{a-1}+p_a)}
	\right)
	\prod_{\substack{1\leq d\leq j\\d\neq a}}
	\widehat{n_{\pm_d}}(p_{d-1}-p_d).
\end{equation*}
This is exactly $\mathcal F_j$ as defined in
\eqref{eq:wave-replacement}. Adding it to the contributions
already computed proves \eqref{eq:correction-derivative}.

\paragraph{3. Differentiation of $\widetilde{\mathcal H}$.}

We first compute the time derivatives of its quadratic terms.
Using \eqref{eq:fourier-time-derivatives},
\begin{equation*}
	\frac{d}{dt}\norm{\nabla Iu}2^2
	=
	2\operatorname{Re}
	\sum_p |p|^2m(p)^2
	\overline{\widehat u(p)}\,\partial_t\widehat u(p)
\end{equation*}
\begin{equation*}
	=
    \operatorname{Re}
	\left(
	-i\sum_{p,l}\sum_\pm
	|p|^2m(p)^2
	\overline{\widehat u(p)}\widehat u(l)
	\widehat{n_\pm}(p-l)
	\right)
	=
	-\frac{i}{2}\sum_{p,l}\sum_\pm
	\bigl(|p|^2m(p)^2-|l|^2m(l)^2\bigr)
	\overline{\widehat u(p)}\widehat u(l)
	\widehat{n_\pm}(p-l).
\end{equation*}

Similarly,
\begin{equation*}
	\frac{d}{dt}\frac12\norm{n_+}2^2
	=
	\operatorname{Re}
	\left(
	-i\sum_q|q|\overline{\widehat{n_+}(q)}
	\sum_k\widehat u(k)\overline{\widehat u(k-q)}
	\right)
\end{equation*}
\begin{equation*}
	=
	\operatorname{Re}
	\left(
	-i\sum_{p,l}|p-l|
	\overline{\widehat u(p)}\widehat u(l)
	\widehat{n_-}(p-l)
	\right)
	=
	\frac{i}{2}\sum_{p,l}|p-l|
	\overline{\widehat u(p)}\widehat u(l)
	\bigl(\widehat{n_+}(p-l)-\widehat{n_-}(p-l)\bigr).
\end{equation*}
Consequently,
\begin{equation}\label{eq:quadratic-energy-derivative}
		\frac{d}{dt}
		\left(\norm{\nabla Iu}2^2+\frac12\norm{n_+}2^2\right)
		=-\frac{i}{2}\sum_{p,l}\sum_\pm
		\bigl(|p|^2m(p)^2-|l|^2m(l)^2\mp|p-l|\bigr)
		\overline{\widehat u(p)}\widehat u(l)
		\widehat{n_\pm}(p-l).
\end{equation}
Since $\sigma_\pm=1-b_{1,\pm}$, we have
\begin{equation*}
	\widetilde{\mathcal H}
	=
	\norm{\nabla Iu}2^2+\frac12\norm{n_+}2^2
	+\mathcal C_3[1]-\mathcal C_3[b_1].
\end{equation*}
We first differentiate $\mathcal C_3[1]$. Substituting
\eqref{eq:fourier-time-derivatives} and applying the product
rule gives
\begin{equation*}
	\frac{d}{dt}\mathcal C_3[1]
	=
	i\mathcal C_3[\Omega_1]
	+i\mathcal C_4[1-1]
	-\frac{i}{2}
	\sum_{p,l,k}
	\overline{\widehat u(p)}\widehat u(l)
	\widehat u(k)\overline{\widehat u(k-p+l)}
	\sum_\pm(\pm|p-l|).
\end{equation*}
Here the term $i\mathcal C_3[\Omega_1]$ comes from the
linear parts of the equations. The term
$i\mathcal C_4[1-1]$ comes from the nonlinear parts of
the equations for $u$ and $\overline u$, as computed in
\eqref{eq:left-u-derivative} and
\eqref{eq:right-u-derivative}. The final sum comes from
the nonlinear part of the equation for $n_\pm$.
Since
\begin{equation*}
	1-1=0,
	\qquad
	\sum_\pm(\pm|p-l|)=0,
\end{equation*}
we obtain
\begin{equation}\label{eq:constant-cubic-derivative}
	\frac{d}{dt}\mathcal C_3[1]
	=i\mathcal C_3[\Omega_1].
\end{equation}

Taking $j=1$ in \eqref{eq:correction-derivative} gives
\begin{equation*}
	\frac{d}{dt}\mathcal C_3[b_1]
	=
	i\mathcal C_3[\Omega_1b_1]
	+i\mathcal C_4[R_2-\Omega_2b_2]
	+\mathcal F_1.
\end{equation*}
Subtracting this identity from
\eqref{eq:constant-cubic-derivative}, we find
\begin{equation}\label{eq:cubic-part-derivative}
	\frac{d}{dt}\mathcal C_3[1-b_1]
	=
	i\mathcal C_3[\Omega_1(1-b_1)]
	-i\mathcal C_4[R_2-\Omega_2b_2]
	-\mathcal F_1.
\end{equation}

We next combine the first term on the right of
\eqref{eq:cubic-part-derivative} with
\eqref{eq:quadratic-energy-derivative}.
Using $g(k)=|k|^2(1-m(k)^2)$, we have
\begin{equation*}
	-\bigl(|p|^2m(p)^2-|l|^2m(l)^2\mp|p-l|\bigr)
	=
	g(p)-g(l)-\Omega_{1,\pm}(p,l).
\end{equation*}
Consequently,
\begin{equation*}
	\frac{d}{dt}
	\left(\norm{\nabla Iu}2^2+\frac12\norm{n_+}2^2\right)
	+i\mathcal C_3[\Omega_1(1-b_1)]
	=
	i\mathcal C_3[g(p)-g(l)-\Omega_1]
	+i\mathcal C_3[\Omega_1(1-b_1)]
\end{equation*}
\begin{equation*}
	=
	i\mathcal C_3[g(p)-g(l)-\Omega_1b_1]
	=-i\mathcal C_3[R_1]
	=-i\mathcal R_1,
\end{equation*}
Adding the remaining terms in
\eqref{eq:cubic-part-derivative} therefore yields
\begin{equation}\label{eq:Htilde-derivative}
	\frac{d}{dt}\widetilde{\mathcal H}
	=
	-i\mathcal R_1
	-i\mathcal C_4[R_2-\Omega_2b_2]
	-\mathcal F_1.
\end{equation}

\paragraph{4. Derivative of $\mathcal E_N$.}

Taking $j=2$ and $j=3$ in
\eqref{eq:correction-derivative}, respectively, gives
\begin{equation*}
	\frac{d}{dt}\mathcal C_4
	=
	i\mathcal C_4[\Omega_2b_2]
	+i\mathcal C_5[R_3-\Omega_3b_3]
	+\mathcal F_2
\end{equation*}
and, since $b_4=0$,
\begin{equation*}
	\frac{d}{dt}\mathcal C_5
	=
	i\mathcal C_5[\Omega_3b_3]
	+i\mathcal C_6[R_4]
	+\mathcal F_3.
\end{equation*}
Subtracting the first identity from
\eqref{eq:Htilde-derivative}, we obtain
\begin{equation*}
	\frac{d}{dt}
	\bigl(\widetilde{\mathcal H}-\mathcal C_4\bigr)
	=
	-i\mathcal R_1
	-i\mathcal C_4[R_2-\Omega_2b_2]
	-i\mathcal C_4[\Omega_2b_2]
	-i\mathcal C_5[R_3-\Omega_3b_3]
	-\mathcal F_1-\mathcal F_2
\end{equation*}
\begin{equation*}
	=
	-i\mathcal R_1-i\mathcal R_2
	-i\mathcal C_5[R_3-\Omega_3b_3]
	-\mathcal F_1-\mathcal F_2.
\end{equation*}
Subtracting the derivative of $\mathcal C_5$ next gives
\begin{equation*}
	\frac{d}{dt}\mathcal E_N
	=
	-i\mathcal R_1-i\mathcal R_2
	-i\mathcal C_5[R_3-\Omega_3b_3]
	-i\mathcal C_5[\Omega_3b_3]
	-i\mathcal C_6[R_4]
	-\mathcal F_1-\mathcal F_2-\mathcal F_3.
\end{equation*}
Hence
\begin{equation}\label{eq:energyidentity}
	\frac{d}{dt}\mathcal E_N
	=-i\sum_{j=1}^4\mathcal R_j
	-\sum_{j=1}^3\mathcal F_j.
\end{equation}
\begin{remark}
	The additional corrections are designed to improve the
	quartic term in \cite[(4.2)]{KishimotoI}. In our notation, this
	term is
	\begin{equation*}
		-i\mathcal C_4\bigl[
		b_{1,\pm_2}(p_1,p_2)-b_{1,\pm_1}(p_0,p_1)
		\bigr],
	\end{equation*}
	and contains two occurrences of $u$ or $\overline u$ and two
	occurrences of $n_\pm$. Its estimate in Kishimoto's argument
	includes the bound
	$N^{-5/4+\eps}h^{1/4-\eps}A^4$, which limits the resulting
	global iteration to $s>9/14$.
	
	Subtracting $\mathcal C_4[b_2]$ preserves this quartic structure
	but changes its coefficient to
	\begin{equation*}
		R_{2,\pm_1,\pm_2}(p_0,p_1,p_2)
		=
		b_{1,\pm_2}(p_1,p_2)-b_{1,\pm_1}(p_0,p_1)
		+\Omega_{2,\pm_1,\pm_2}(p_0,p_1,p_2)
		b_{2,\pm_1,\pm_2}(p_0,p_1,p_2).
	\end{equation*}
	Thus the natural choice, whenever it is well-defined, is
	\begin{equation*}
		b_{2,\pm_1,\pm_2}(p_0,p_1,p_2)
		=
		\frac{
			b_{1,\pm_1}(p_0,p_1)-b_{1,\pm_2}(p_1,p_2)
		}{
			\Omega_{2,\pm_1,\pm_2}(p_0,p_1,p_2)
		}.
	\end{equation*}
	This makes $R_2$ vanish in the selected region. 
	
	Differentiating $\mathcal C_4[b_2]$ also produces a
	quintic term containing two occurrences of $u$ or
	$\overline u$ and three occurrences of $n_\pm$, with coefficient
	\begin{equation*}
		b_{2,\pm_2,\pm_3}(p_1,p_2,p_3)
		-b_{2,\pm_1,\pm_2}(p_0,p_1,p_2).
	\end{equation*}
	The quintic correction $\mathcal C_5[b_3]$ treats this difference in the same
	way. Its natural coefficient is
	\begin{equation*}
		b_{3,\pm_1,\pm_2,\pm_3}(p_0,p_1,p_2,p_3)
		=
		\frac{
			b_{2,\pm_1,\pm_2}(p_0,p_1,p_2)
			-b_{2,\pm_2,\pm_3}(p_1,p_2,p_3)
		}{
			\Omega_{3,\pm_1,\pm_2,\pm_3}(p_0,p_1,p_2,p_3)
		},
	\end{equation*}
	again only where this division is justified. 
	
	The terms $\mathcal R_1$ and 
	$\mathcal F_1$ are unchanged by the additional corrections.
	They can therefore be estimated as in the treatment of
	\cite[Proposition 4.1]{KishimotoI},
	using the trilinear estimates for $\mathcal R_1$ and the
	frequency decomposition and H\"older estimates for
	$\mathcal F_1$. For $1/2<s<1$, these give
	\begin{equation}\label{eq:oldcubic}
		\bigg|\int_J\mathcal R_1(t)\,dt\bigg|
		\lesssim_\eps N^{-1+\eps}h^{1/2-\eps}A^3
	\end{equation}
	and
	\begin{equation}\label{eq:oldquartic}
		\bigg|\int_J\mathcal F_1(t)\,dt\bigg|
		\lesssim_\eps
		(N^{-2+\eps}+N^{-1+\eps}h^{1-\eps})A^4.
	\end{equation}

	It turns out that these bounds impose no regularity restriction beyond
	$s>1/2$ in the global iteration. 
	
	In summary, the intended improvement therefore comes from cancelling
	the difficult part of the quartic term before
	estimating it, rather than applying a stronger estimate
	directly to the same term. 
\end{remark}

\vspace{10pt}
\section{The quartic correction}\label{sec:quartic}

We now specify the coefficient $b_2$ in
\eqref{eq:C4-definition}. Recall that
\begin{equation*}
	R_{2,\pm_1,\pm_2}(p,r,l)
	=
	\Omega_{2,\pm_1,\pm_2}(p,r,l)b_{2,\pm_1,\pm_2}(p,r,l)
	+b_{1,\pm_2}(r,l)-b_{1,\pm_1}(p,r),
\end{equation*}
where
\begin{equation*}
	\Omega_{2,\pm_1,\pm_2}(p,r,l)
	=
	|p|^2-|l|^2\mp_1|p-r|\mp_2|r-l|.
\end{equation*}
We choose $b_2$ to cancel this coefficient in the regions
specified below. Fix $c=1/100$ and take $N\geq10^4$.

\paragraph{1. Definition of $b_2$.}

If $\max(|p|,|l|)\leq N/4$ and $g(r)\neq0$, define
\begin{equation}\label{eq:b2low}
	b_{2,\pm_1,\pm_2}(p,r,l)
	=
	-\frac{g(r)}
	{\bigl(|p|^2-|r|^2\mp_1|p-r|\bigr)
		\bigl(|r|^2-|l|^2\mp_2|r-l|\bigr)}.
\end{equation}
Note that $g(r)\neq0$ implies $|r|>N$ and $|p|,|l|\leq N/4$,
so both denominator factors have absolute value comparable to
$|r|^2$ for sufficiently large $N$, and hence are nonzero.
If $\max(|p|,|l|)\leq N/4$ and $g(r)=0$, set
$b_{2,\pm_1,\pm_2}(p,r,l)=0$.

In the region
\begin{equation*}
	\max(|p|,|l|)>N/4,\qquad
	\min(|p|,|l|)\leq N/4,\qquad \max(|p|,|l|)\ge 4\min(|p|,|l|),
\end{equation*}
subject also to
\begin{equation*}
	\max(|p|,|r|,|l|)
	\leq c\max(|p|,|l|)^2,
\end{equation*}
define
\begin{equation}\label{eq:b2highlow}
	b_{2,\pm_1,\pm_2}(p,r,l)
	=
	\frac{b_{1,\pm_1}(p,r)-b_{1,\pm_2}(r,l)}
	{\Omega_{2,\pm_1,\pm_2}(p,r,l)}.
\end{equation}
In all remaining regions, set
$b_{2,\pm_1,\pm_2}(p,r,l)=0$.

We first check the denominator in (\ref{eq:b2highlow}) is non-vanishing. The assumptions give
\begin{equation*}
	\bigl||p|^2-|l|^2\bigr|
	\geq\frac{15}{16}\max(|p|,|l|)^2,
\end{equation*}
whereas
\begin{equation*}
	|p-r|+|r-l|
	\leq4\max(|p|,|r|,|l|)
	\leq4c\max(|p|,|l|)^2.
\end{equation*}
Hence
\begin{equation}\label{eq:b2phase-lower}
	|\Omega_{2,\pm_1,\pm_2}(p,r,l)|
	\geq
	\left(\frac{15}{16}-4c\right)\max(|p|,|l|)^2>0.
\end{equation}

The defining regions are invariant under
$(p,r,l;\pm_1,\pm_2)\mapsto(l,r,p;\mp_2,\mp_1)$.
Moreover,
\begin{equation*}
	\Omega_{2,\mp_2,\mp_1}(l,r,p)
	=-\Omega_{2,\pm_1,\pm_2}(p,r,l).
\end{equation*}
 Thus, one can check that
\begin{equation}\label{sym}
	b_{2,\pm_1,\pm_2}(p,r,l)
	=b_{2,\mp_2,\mp_1}(l,r,p),
\end{equation}
which satisfies the assumption (\ref{eq:coefficient-reversal}).

\begin{remark}
	The choice of $b_2$ is guided by the cancellation condition
	\begin{equation*}
		\Omega_{2,\pm_1,\pm_2}(p,r,l)b_{2,\pm_1,\pm_2}(p,r,l)
		=
		b_{1,\pm_1}(p,r)-b_{1,\pm_2}(r,l).
	\end{equation*}
	In \eqref{eq:b2highlow}, the phase $\Omega_2$ has a lower bound,
	so we define $b_2$ by direct division.
	In \eqref{eq:b2low}, such a lower bound is unavailable,
	but the numerator contains the same phase as a factor.
	Cancelling this factor gives \eqref{eq:b2low} without
	dividing by a potentially vanishing quantity.
	Thus the two formulas actually arise from the same cancellation condition. See \eqref{eq:b2identity} for the detailed identities.
\end{remark}

We next estimate the size of $b_2$ and the quartic term
$\mathcal R_2$ resulting from this choice.
\begin{lemma}\label{lem:b2}
	If $b_{2,\pm_1,\pm_2}(p,r,l)\neq0$, then
	\begin{equation}\label{as}
		\max(|p|,|r|,|l|)>N,
		\qquad
		\min(|p|,|l|)\leq N/4.
	\end{equation}
	Furthermore,
	\begin{equation}\label{eq:b2bound}
		|b_{2,\pm_1,\pm_2}(p,r,l)|
		\lesssim \max(|p|,|r|,|l|)^{-2}
	\end{equation}
	on its nonzero support. For smooth solutions on $J$,
	\begin{equation}\label{eq:quarticremainder}
		\bigg|\int_J\mathcal R_2(t)\,dt\bigg|
		\lesssim_{\eps}
		\bigl(N^{-2+\eps}+N^{-1+\eps}h^{1-\eps}\bigr)A^4.
	\end{equation}
\end{lemma}
\begin{proof}
	
	\textbf{The bound for $b_2$.}
	
	In \eqref{eq:b2low}, the restriction directly gives
	\begin{equation*}
		|b_{2,\pm_1,\pm_2}(p,r,l)|
		\lesssim \frac{g(r)}{|r|^4}
		\lesssim|r|^{-2}.
	\end{equation*}
	Since $|p|,|l|\leq|r|/4$, this proves
	\eqref{eq:b2bound} in that region.
	
	Consider \eqref{eq:b2highlow}. If
	$|r|\leq4\max(|p|,|l|)$, then
	\begin{equation*}
		|b_{2,\pm_1,\pm_2}(p,r,l)|
		\lesssim \max(|p|,|l|)^{-2}
		\lesssim\max(|p|,|r|,|l|)^{-2}
	\end{equation*}
	by \eqref{eq:b2phase-lower} and the boundedness of $b_1$.
	
	It remains to consider
	$|r|>4\max(|p|,|l|)$.
	By the symmetry (\ref{sym}), we may assume
	\begin{equation*}
		|p|=\max(|p|,|l|)>N/4,
		\qquad |l|\leq N/4,
		\qquad g(l)=0.
	\end{equation*}
	By \eqref{eq:b1non},
	\begin{equation*}
		b_{1,\pm_1}(p,r)
		=
		\frac{g(r)-g(p)}
		{|r|^2-|p|^2\pm_1|p-r|},
		\qquad
		b_{1,\pm_2}(r,l)
		=
		\frac{g(r)}
		{|r|^2-|l|^2\mp_2|r-l|}.
	\end{equation*}
	Subtracting these fractions gives
	\begin{equation}\label{eq:b1difference}
		b_{1,\pm_1}(p,r)-b_{1,\pm_2}(r,l)
		=
		\frac{g(r)\Omega_{2,\pm_1,\pm_2}(p,r,l)}
		{\bigl(|r|^2-|p|^2\pm_1|p-r|\bigr)
			\bigl(|r|^2-|l|^2\mp_2|r-l|\bigr)}
		-\frac{g(p)}
		{|r|^2-|p|^2\pm_1|p-r|}.
	\end{equation}
	Each denominator factor in \eqref{eq:b1difference} has
	absolute value at least $|r|^2/2$. Using \eqref{eq:b2phase-lower}, we obtain
	\begin{equation*}
		|b_{2,\pm_1,\pm_2}(p,r,l)|
		\lesssim
		\frac{g(r)}{|r|^4}
		+\frac{g(p)}
		{|r|^2|\Omega_{2,\pm_1,\pm_2}(p,r,l)|}
		\lesssim |r|^{-2}.
	\end{equation*}
	This proves \eqref{eq:b2bound}.
	
	Finally, if $\max(|p|,|r|,|l|)\leq N$, then
	\begin{equation*}
		g(p)=g(r)=g(l)=0,
		\qquad
		b_{1,\pm_1}(p,r)=b_{1,\pm_2}(r,l)=0.
	\end{equation*}
	Both definitions of $b_2$ therefore give zero.
	The condition $\min(|p|,|l|)\leq N/4$ follows directly
	from the defining regions. This concludes (\ref{as}).
	
	\medskip
	\noindent\textbf{The regions where $R_2=0$.}
	
	In the region of \eqref{eq:b2highlow}, substitution gives
	\begin{equation*}
		R_{2,\pm_1,\pm_2}(p,r,l)
		=
		b_{1,\pm_1}(p,r)-b_{1,\pm_2}(r,l)
		+b_{1,\pm_2}(r,l)-b_{1,\pm_1}(p,r)
		=0.
	\end{equation*}
	
	Suppose next that $\max(|p|,|l|)\leq N/4$.
	If $g(r)=0$, all three frequencies have vanishing $g$,
	and both $b_1$ coefficients vanish.
	Thus $R_2=0$.
	
	If $g(r)\neq0$, then $|r|>N$ and
	\begin{equation*}
		b_{1,\pm_1}(p,r)
		=-\frac{g(r)}{|p|^2-|r|^2\mp_1|p-r|},
		\qquad
		b_{1,\pm_2}(r,l)
		=\frac{g(r)}{|r|^2-|l|^2\mp_2|r-l|}.
	\end{equation*}
	Adding the two denominators gives $\Omega_2$, so
	\begin{equation}\label{eq:b2identity}
		b_{1,\pm_1}(p,r)-b_{1,\pm_2}(r,l)
		=
		-\frac{g(r)\Omega_{2,\pm_1,\pm_2}(p,r,l)}
		{\bigl(|p|^2-|r|^2\mp_1|p-r|\bigr)
			\bigl(|r|^2-|l|^2\mp_2|r-l|\bigr)}.
	\end{equation}
	By \eqref{eq:b2low},
	\begin{equation*}
		b_{1,\pm_1}(p,r)-b_{1,\pm_2}(r,l)
		=\Omega_{2,\pm_1,\pm_2}(p,r,l)
		b_{2,\pm_1,\pm_2}(p,r,l),
	\end{equation*}
	and hence $R_2=0$.
	
	Outside the above regions, $b_2=0$ and
	\begin{equation}\label{eq:R2-outside}
		R_{2,\pm_1,\pm_2}(p,r,l)
		=b_{1,\pm_2}(r,l)-b_{1,\pm_1}(p,r).
	\end{equation}
	In particular, $|R_2|\lesssim1$ everywhere.
	
\medskip
\noindent\textbf{ The regions where $R_2$ may be nonzero.}

We first identify the frequencies that remain after the
cancellations proved above. If $R_{2,\pm_1,\pm_2}(p,r,l)\neq0$,
then at least one of the following three conditions holds:
\begin{equation}\label{eq:R2-region-one}
	|p|\geq N/16,\qquad |l|\geq N/16;
\end{equation}
\begin{equation}\label{eq:R2-region-two}
	|l|<N/16,\qquad |p|>N/4,\qquad |r|>c|p|^2;
\end{equation}
\begin{equation}\label{eq:R2-region-three}
	|p|<N/16,\qquad |l|>N/4,\qquad |r|>c|l|^2.
\end{equation}

To verify this inclusion, suppose that
\eqref{eq:R2-region-one} fails. After exchanging $p,l$
and replacing $(\pm_1,\pm_2)$ by $(\mp_2,\mp_1)$ if necessary,
we may assume $|l|<N/16$. The cancellation in
\eqref{eq:b2low} implies
\begin{equation*}
	R_{2,\pm_1,\pm_2}(p,r,l)\neq0
	\quad\Longrightarrow\quad
	|p|>N/4>4|l|.
\end{equation*}
The cancellation in \eqref{eq:b2highlow} then implies
\begin{equation*}
	\max(|p|,|r|,|l|)>c|p|^2.
\end{equation*}
Since $N\geq10^4$ and $c=1/100$, we can deduce that
\begin{equation*}
	|r|>c|p|^2,
\end{equation*}
which is \eqref{eq:R2-region-two}.

For either region, decompose
\begin{equation*}
	u=\sum_{N_1}u_{N_1},
	\qquad
	u=\sum_{N_2}u_{N_2},
	\qquad
	n_{\pm_1}=\sum_{N_3}n_{\pm_1,N_3},
	\qquad
	n_{\pm_2}=\sum_{N_4}n_{\pm_2,N_4},
\end{equation*}
where $u_{N_a}=P_{N_a}u$ and
$n_{\pm,N_a}=P_{N_a}n_\pm$. In every nonzero summand,
\begin{equation}\label{eq:R2-dyadic-support}
	|p|\sim N_1,\qquad |l|\sim N_2,\qquad
	|p-r|\sim N_3,\qquad |r-l|\sim N_4.
\end{equation}
For $N_a=1$, the corresponding condition is $|k|\lesssim1$.

\medskip
\noindent\textbf{ The region \eqref{eq:R2-region-one}.}

For fixed signs and dyadic indices, we estimate
\begin{equation*}
	\left|
	\int_J
	\sum_{\substack{p,r,l\\|p|,|l|\geq N/16}}
	R_{2,\pm_1,\pm_2}(p,r,l)
	\overline{\widehat{u_{N_1}}(t,p)}
	\widehat{u_{N_2}}(t,l)
	\widehat{n_{\pm_1,N_3}}(t,p-r)
	\widehat{n_{\pm_2,N_4}}(t,r-l)\,dt
	\right|.
\end{equation*}
The support conditions imply $N_1,N_2\geq N/32$.

Insert $\psi_h$ on both occurrences of $u$ and $1_J$
on both occurrences of $n$. Since $\psi_h^2 1_J^2=1_J$,
the integral is unchanged. From the discussion in Remark \ref{absolute} and $|R_2|\lesssim1$, we apply
Cauchy--Schwarz and obtain
\begin{equation*}
	\left|
	\int_J
	\sum_{\substack{p,r,l\\|p|,|l|\geq N/16}}
	R_{2,\pm_1,\pm_2}(p,r,l)
	\overline{\widehat{u_{N_1}}(t,p)}
	\widehat{u_{N_2}}(t,l)
	\widehat{n_{\pm_1,N_3}}(t,p-r)
	\widehat{n_{\pm_2,N_4}}(t,r-l)\,dt
	\right|
\end{equation*}
\begin{equation*}
	\lesssim
	\left\|
	\left(\mathcal F_{t,x}^{-1}
	\left|\widetilde{\psi_h\overline{u_{N_1}}}\right|\right)
	\left(\mathcal F_{t,x}^{-1}
	\left|\widetilde{1_Jn_{\pm_1,N_3}}\right|\right)
	\right\|_{L^2_{t,x}}
	\left\|
	\left(\mathcal F_{t,x}^{-1}
	\left|\widetilde{\psi_hu_{N_2}}\right|\right)
	\left(\mathcal F_{t,x}^{-1}
	\left|\widetilde{1_Jn_{\pm_2,N_4}}\right|\right)
	\right\|_{L^2_{t,x}}.
\end{equation*}
Applying \eqref{eq:pair} with $\eps/8$ in place of $\eps$
bounds the last product by
\begin{equation*}
	C_{\eps}
	m(N_1)^{-1}m(N_2)^{-1}
	A_{N_1}A_{N_2}B_{N_3}B_{N_4}
\end{equation*}
\begin{equation*}
	\qquad{}\times
	\bigl(N_1^{-1+\eps/8}
	+h^{1/2-\eps/8}N_1^{-1/2+\eps/8}\bigr)
	\bigl(N_2^{-1+\eps/8}
	+h^{1/2-\eps/8}N_2^{-1/2+\eps/8}\bigr).
\end{equation*}

We now sum in $N_3,N_4$. \eqref{eq:R2-dyadic-support} implies
\begin{equation*}
	\bigl||p-r|-|r-l|\bigr|
	\leq |p-l|
	\leq2N_1+2N_2.
\end{equation*}
Thus, we have
\begin{equation*}
	\max(N_3,N_4)\lesssim \max(N_1,N_2)
	\quad\text{or}\quad
	N_3\sim N_4.
\end{equation*}
For the first possibility,
\begin{equation*}
	\sum_{N_3,N_4\lesssim\max(N_1,N_2)}
	B_{N_3}B_{N_4}
	=
	\left(
	\sum_{N_3\lesssim \max(N_1,N_2)}B_{N_3}
	\right)^2
	\leq
	\#\{N_3\text{ dyadic}:N_3\lesssim\max(N_1,N_2)\}
	\sum_{N_3}B_{N_3}^2
\end{equation*}
\begin{equation*}
	\lesssim
	\bigl(1+\log(2+\max(N_1,N_2))\bigr)A^2
	\lesssim_\eps (N_1N_2)^{\eps/8}A^2.
\end{equation*}
For the second possibility,
\begin{equation*}
	\sum_{N_3\sim N_4}B_{N_3}B_{N_4}
	\lesssim\sum_{N_3}B_{N_3}^2
	\lesssim A^2.
\end{equation*}
Using $A_{N_1},A_{N_2}\lesssim A$ and summing the
preceding estimates, we obtain
\begin{equation*}
	\left|
	\int_J
	\sum_{\substack{p,r,l\\|p|,|l|\geq N/16}}
	\sum_{\pm_1,\pm_2}
	R_{2,\pm_1,\pm_2}(p,r,l)
	\overline{\widehat u(t,p)}\widehat u(t,l)
	\widehat{n_{\pm_1}}(t,p-r)
	\widehat{n_{\pm_2}}(t,r-l)\,dt
	\right|
\end{equation*}
\begin{equation*}
	\lesssim_{\eps}
	A^4
	\left[
	\sum_{N_1\geq N/32}
	m(N_1)^{-1}
	\bigl(
	N_1^{-1+\eps/4}
	+h^{1/2-\eps/8}N_1^{-1/2+\eps/4}
	\bigr)
	\right]^2.
\end{equation*}

For $N_1\ge N/32$, we have
$m(N_1)^{-1}\sim N^{s-1}N_1^{1-s}$. Choose
$0<\eps<\min(1,4(s-1/2))$. Then
\begin{equation*}
	\sum_{N_1\ge N/32}
	m(N_1)^{-1}N_1^{-1+\eps/4}
	\sim
	N^{s-1}\sum_{N_1\ge N/32}N_1^{-s+\eps/4}
	\lesssim_{\eps}N^{-1+\eps/4},
\end{equation*}
and
\begin{equation*}
	\sum_{N_1\ge N/32}
	m(N_1)^{-1}N_1^{-1/2+\eps/4}
	\sim
	N^{s-1}\sum_{N_1\ge N/32}N_1^{1/2-s+\eps/4}
	\lesssim_{s,\eps}N^{-1/2+\eps/4}.
\end{equation*}
Thus
\begin{equation*}
	\left|
	\int_J
	\sum_{\substack{p,r,l\\|p|,|l|\geq N/16}}
	\sum_{\pm_1,\pm_2}
	R_{2,\pm_1,\pm_2}(p,r,l)
	\overline{\widehat u(t,p)}\widehat u(t,l)
	\widehat{n_{\pm_1}}(t,p-r)
	\widehat{n_{\pm_2}}(t,r-l)\,dt
	\right|
\end{equation*}
\begin{equation*}
	\lesssim_{\eps}
	\bigl(N^{-2+\eps/2}
	+h^{1-\eps/4}N^{-1+\eps/2}\bigr)A^4
	\leq
	\bigl(N^{-2+\eps}
	+h^{1-\eps}N^{-1+\eps}\bigr)A^4.
\end{equation*}

\medskip
\noindent\textbf{ The region \eqref{eq:R2-region-two}.}

We now have the restrictions:
\begin{equation*}
	|l|<N/16,\qquad |p|>N/4,\qquad |r|>c|p|^2.
\end{equation*}
Consequently, 
\begin{equation*}
	|p-r|\sim |r-l|\sim |r|.
\end{equation*}
In this region $b_2=0$. By \eqref{eq:R2-outside} and
\eqref{eq:b1difference}, one can obtain
\begin{equation*}
	|R_{2,\pm_1,\pm_2}(p,r,l)|
	\lesssim
	\frac{g(r)|\Omega_{2,\pm_1,\pm_2}(p,r,l)|}{|r|^4}
	+\frac{g(p)}{|r|^2}.
\end{equation*}
Since
\begin{equation*}
	|\Omega_{2,\pm_1,\pm_2}(p,r,l)|
	\leq |p|^2+|l|^2+|p-r|+|r-l|
	\lesssim |p|^2+|r|,
\end{equation*}
we obtain
\begin{equation*}
	|R_{2,\pm_1,\pm_2}(p,r,l)|
	\lesssim \frac1{|r|}+
	\frac{|p|^2}{|r|^2}
	\lesssim |r|^{-1}.
\end{equation*}

Keep the dyadic notation \eqref{eq:R2-dyadic-support}.
The present restrictions imply
\begin{equation*}
	N_1\geq N/8,\qquad
	N_2\leq N/8,\qquad
	N_3\gtrsim N_1^2,\qquad
	N_4\sim N_3.
\end{equation*}
In particular, $|r|\sim N_3$ and $m(N_2)=1$.

For fixed signs and dyadic indices, one can similarly derive
\begin{equation*}
	\left|
	\int_J
	\sum_{\substack{p,r,l\\
			|l|<N/16,\ |p|>N/4\\
			|r|>c|p|^2}}
	R_{2,\pm_1,\pm_2}(p,r,l)
	\overline{\widehat{u_{N_1}}(t,p)}
	\widehat{u_{N_2}}(t,l)
	\widehat{n_{\pm_1,N_3}}(t,p-r)
	\widehat{n_{\pm_2,N_4}}(t,r-l)\,dt
	\right|
\end{equation*}
\begin{equation*}
	\lesssim
	N_3^{-1}\int_J
	\left\|\mathcal F_x^{-1}
	|\widehat{u_{N_1}}(t,\cdot)|\right\|_{L_x^\infty}
	\left\|\mathcal F_x^{-1}
	|\widehat{u_{N_2}}(t,\cdot)|\right\|_{L_x^\infty}
	\|n_{\pm_1,N_3}(t)\|_{L_x^2}
	\|n_{\pm_2,N_4}(t)\|_{L_x^2}\,dt
\end{equation*}
\begin{equation*}
	\lesssim
	hN_3^{-1}m(N_1)^{-1}
	A_{N_1}A_{N_2}B_{N_3}B_{N_4},
\end{equation*}
where we used \eqref{eq:spatialu},
\eqref{eq:spatialwave}, and $m(N_2)=1$.

Summing these estimates yields
\begin{equation}\label{eq:quarticfar}
	\left|
	\int_J
	\sum_{\substack{p,r,l\\
			|l|<N/16,\ |p|>N/4\\
			|r|>c|p|^2}}
	\sum_{\pm_1,\pm_2}
	R_{2,\pm_1,\pm_2}(p,r,l)
	\overline{\widehat u(t,p)}\widehat u(t,l)
	\widehat{n_{\pm_1}}(t,p-r)
	\widehat{n_{\pm_2}}(t,r-l)\,dt
	\right|
\end{equation}
\begin{equation*}
	\lesssim
	h\sum_{N_1\geq N/8}
	m(N_1)^{-1}A_{N_1}
	\sum_{N_2\leq N/8}A_{N_2}
	\sum_{N_3\gtrsim N_1^2}N_3^{-1}
	\sum_{N_4\sim N_3}B_{N_3}B_{N_4}.
\end{equation*}
By Cauchy--Schwarz,
\begin{equation*}
	\sum_{N_2\leq N/8}A_{N_2}
	\lesssim
	(1+\log(2+N))^{1/2}
	\left(\sum_{N_2}A_{N_2}^2\right)^{1/2}
	\lesssim_\eps N^\eps A.
\end{equation*}
Also,
\begin{equation*}
	A_{N_1}\lesssim A,
	\qquad
	\sum_{N_4\sim N_3}B_{N_3}B_{N_4}\lesssim A^2, \qquad 	\sum_{N_3\gtrsim N_1^2}N_3^{-1}
	\lesssim N_1^{-2}.
\end{equation*}
Thus the right-hand side of \eqref{eq:quarticfar} can be bounded by
\begin{equation*}
	 h N^\eps A^4\sum_{N_1\geq N/8}N_1^{-2}m(N_1)^{-1}
	\lesssim h N^\eps A^4
	\left(N^{s-1}\sum_{N_1\ge N/8}N_1^{-1-s}\right)
	\lesssim  h N^{-2+\eps}A^4.
\end{equation*}
Consequently,
\begin{equation*}
	\left|
	\int_J
	\sum_{\substack{p,r,l\\
			|l|<N/16,\ |p|>N/4\\
			|r|>c|p|^2}}
	\sum_{\pm_1,\pm_2}
	R_{2,\pm_1,\pm_2}(p,r,l)
	\overline{\widehat u(t,p)}\widehat u(t,l)
	\widehat{n_{\pm_1}}(t,p-r)
	\widehat{n_{\pm_2}}(t,r-l)\,dt
	\right|
\end{equation*}
\begin{equation*}
	\lesssim_{\eps}hN^{-2+\eps}A^4.
\end{equation*}
Exchanging $p,l$ and reversing the signs gives the same
bound for \eqref{eq:R2-region-three}.
Since $R_2$ vanishes outside the three regions listed above,
the definition of $\mathcal R_2$ now gives
\begin{equation*}
	\bigg|\int_J\mathcal R_2(t)\,dt\bigg|
	\lesssim_{\eps}
	\bigl(N^{-2+\eps}+N^{-1+\eps}h^{1-\eps}
	+hN^{-2+\eps}\bigr)A^4
	\lesssim
	\bigl(N^{-2+\eps}+N^{-1+\eps}h^{1-\eps}\bigr)A^4,
\end{equation*}
which proves \eqref{eq:quarticremainder}. 

\end{proof}

\vspace{10pt}
\section{The quintic correction}\label{sec:quintic}

We now specify the coefficient $b_3$ in
\eqref{eq:C5-definition}. Recall that
\begin{equation*}
	R_{3,\pm_1,\pm_2,\pm_3}(p,r,z,l)
	=
	\Omega_{3,\pm_1,\pm_2,\pm_3}(p,r,z,l)
	b_{3,\pm_1,\pm_2,\pm_3}(p,r,z,l)
\end{equation*}
\begin{equation*}
	\qquad{}
	+b_{2,\pm_2,\pm_3}(r,z,l)
	-b_{2,\pm_1,\pm_2}(p,r,z).
\end{equation*}
Recall that
\begin{equation*}
	\Omega_{1,\pm_1}(p,r)=|p|^2-|r|^2\mp_1|p-r|,
	\qquad
	\Omega_{1,\pm_2}(r,z)=|r|^2-|z|^2\mp_2|r-z|,
\end{equation*}
\begin{equation*}
	\Omega_{1,\pm_3}(z,l)=|z|^2-|l|^2\mp_3|z-l|,
\end{equation*}
\begin{equation}\label{2}
	\Omega_{2,\pm_1,\pm_2}(p,r,z)=\Omega_{1,\pm_1}(p,r)+\Omega_{1,\pm_2}(r,z),
	\qquad
	\Omega_{2,\pm_2,\pm_3}(r,z,l)=\Omega_{1,\pm_2}(r,z)+\Omega_{1,\pm_3}(z,l),
\end{equation}
and
\begin{equation}\label{eq:quintic-phase-identities}
	\Omega_{3,\pm_1,\pm_2,\pm_3}(p,r,z,l)
	=\Omega_{1,\pm_1}(p,r)+\Omega_{2,\pm_2,\pm_3}(r,z,l)
	=\Omega_{2,\pm_1,\pm_2}(p,r,z)+\Omega_{1,\pm_3}(z,l).
\end{equation}
We retain $c=1/100$ and $N\geq10^4$.

\paragraph{1. Definition of $b_3$.}

We define $b_3$ on the five disjoint regions below and set
$b_3=0$ elsewhere. In \textup{(i)--(iv)}, we assume
\begin{equation}\label{eq:b3-small-exterior}
	\max(|p|,|l|)\leq N/16.
\end{equation}

\medskip
\noindent\textbf{(i) $|r|>N/4$ and $|z|\leq N/4$.}

If $g(r)\neq0$, define
\begin{equation}\label{eq:b3-single-r}
	b_{3,\pm_1,\pm_2,\pm_3}(p,r,z,l)
	=
	-\frac{g(r)}
	{\Omega_{1,\pm_1}(p,r)
		\Omega_{1,\pm_2}(r,z)
		\Omega_{2,\pm_2,\pm_3}(r,z,l)}.
\end{equation}
If $g(r)=0$, set $b_{3,\pm_1,\pm_2,\pm_3}(p,r,z,l)=0$.

\medskip
\noindent\textbf{(ii) $|z|>N/4$ and $|r|\leq N/4$.}

If $g(z)\neq0$, define
\begin{equation}\label{eq:b3-single-z}
	b_{3,\pm_1,\pm_2,\pm_3}(p,r,z,l)
	=
	\frac{g(z)}
	{\Omega_{2,\pm_1,\pm_2}(p,r,z)
		\Omega_{1,\pm_2}(r,z)
		\Omega_{1,\pm_3}(z,l)}.
\end{equation}
If $g(z)=0$, set $b_{3,\pm_1,\pm_2,\pm_3}(p,r,z,l)=0$.

\medskip
\noindent\textbf{(iii) $\min(|r|,|z|)>N/4$ and
	$\max(|r|,|z|)>4\min(|r|,|z|)$.}

We additionally require
\begin{equation}\label{eq:b3-internal-cutoff}
	\max(|r|,|z|)\leq c\min(|r|,|z|)^2.
\end{equation}
On this region, define
\begin{equation}\label{eq:b3separated}
	b_{3,\pm_1,\pm_2,\pm_3}(p,r,z,l)
	=
	-\frac{g(r)}
	{\Omega_{1,\pm_1}(p,r)
		\Omega_{1,\pm_2}(r,z)
		\Omega_{2,\pm_2,\pm_3}(r,z,l)}
\end{equation}
\begin{equation*}
	\qquad{}
	+\frac{g(z)}
	{\Omega_{2,\pm_1,\pm_2}(p,r,z)
		\Omega_{1,\pm_2}(r,z)
		\Omega_{1,\pm_3}(z,l)}.
\end{equation*}

\medskip
\noindent\textbf{(iv) $\min(|r|,|z|)>N/4$ and
	$\max(|r|,|z|)\leq4\min(|r|,|z|)$.}

We again require \eqref{eq:b3-internal-cutoff}.
On this region, define
\begin{equation}\label{eq:b3comparable}
	b_{3,\pm_1,\pm_2,\pm_3}(p,r,z,l)
	=
	-\frac{b_{1,\pm_2}(r,z)}2
	\left(
	\frac1{\Omega_{2,\pm_2,\pm_3}(r,z,l)
		\Omega_{1,\pm_1}(p,r)}
	+\frac1{\Omega_{2,\pm_1,\pm_2}(p,r,z)
		\Omega_{1,\pm_3}(z,l)}
	\right)
\end{equation}
\begin{equation*}
	\qquad{}
	-\frac{g(r)+g(z)}
	{2\Omega_{1,\pm_1}(p,r)\Omega_{1,\pm_3}(z,l)}
	\left(
	\frac1{\Omega_{2,\pm_1,\pm_2}(p,r,z)}
	-\frac1{\Omega_{2,\pm_2,\pm_3}(r,z,l)}
	\right).
\end{equation*}

\medskip
\noindent\textbf{(v) $\max(|p|,|l|)>N/16$.}

In this case, we require
\begin{equation}\label{eq:b3-large-exterior}
	\max(|p|,|l|)>N/16,
	\qquad
	\min(|p|,|l|)\leq N/64,
\end{equation}
and
\begin{equation}\label{eq:b3-exterior-cutoff}
	\max(|p|,|r|,|z|,|l|)
	\leq c\max(|p|,|l|)^2.
\end{equation}
On this region, define
\begin{equation}\label{eq:b3highlow}
	b_{3,\pm_1,\pm_2,\pm_3}(p,r,z,l)
	=
	\frac{
		b_{2,\pm_1,\pm_2}(p,r,z)
		-b_{2,\pm_2,\pm_3}(r,z,l)
	}{
		\Omega_{3,\pm_1,\pm_2,\pm_3}(p,r,z,l)
	}.
\end{equation}

We also set $b_3=0$ whenever none of
\textup{(i)--(v)}, including all the conditions specified
in the corresponding case, applies.

\begin{remark}
	The definitions of $b_3$ are guided by the condition $R_3=0$.
	In \eqref{eq:b3highlow}, the phase $\Omega_3$ has a lower bound,
	so we divide the adjacent $b_2$ difference directly by
	$\Omega_3$.
	For \eqref{eq:b3-single-r} and \eqref{eq:b3-single-z},
	the numerator contains $\Omega_3$ as a factor, which
	can be cancelled; see (\ref{eq:b3-single-r-factorization})
	and the corresponding identity obtained by exchanging
	the frequencies and reversing the signs.
	Similar factorization gives \eqref{eq:b3separated};
	see (\ref{eq:b3-separated-r-factorization}) and
	(\ref{eq:b3-separated-z-factorization}).
	
	For \eqref{eq:b3comparable}, we first use the relation
	between $g(r)-g(z)$ and
	$\Omega_{1,\pm_2}(r,z)b_{1,\pm_2}(r,z)$ to extract $\Omega_3$, and average the two expressions obtained
	by eliminating $g(r)$ or $g(z)$.
	This gives a symmetric formula with a factor
	$\Omega_{3}$;
	see (\ref{eq:centralresidual}).
	Thus the construction gives exact cancellation in the nonresonant region and a controlled
	remainder otherwise.
\end{remark}

\paragraph{2. The denominators and the symmetry of $b_3$.}
We first check that all denominators defined above are non-vanishing.
In \eqref{eq:b3-single-r}, a nonzero numerator implies
$|r|>N$. Hence
\begin{equation*}
	|p|,|z|,|l|\leq |r|/4.
\end{equation*}
It follows that
\begin{equation}\label{1}
	|\Omega_{1,\pm_1}(p,r)|,\ |\Omega_{1,\pm_2}(r,z)|, \ |\Omega_{2,\pm_2,\pm_3}(r,z,l)|
	\gtrsim |r|^2,
\end{equation}
The corresponding bounds for \eqref{eq:b3-single-z}
follow by exchanging $p,l$ and $r,z$ and reversing the signs.

Under \eqref{eq:b3-small-exterior},
$\min(|r|,|z|)>N/4$, and \eqref{eq:b3-internal-cutoff},
we have
\begin{equation*}
	|p|,|l|\leq\frac14\min(|r|,|z|),
\end{equation*}
and
\begin{equation*}
	|p-r|+|r-z|+|z-l|
	\leq6\max(|r|,|z|)
	\leq6c\min(|r|,|z|)^2.
\end{equation*}
Consequently,
\begin{equation}\label{eq:b3-denominator-r}
	|\Omega_{1,\pm_1}(p,r)|\sim |r|^2,
	\qquad
	|\Omega_{2,\pm_2,\pm_3}(r,z,l)|\sim |r|^2,
\end{equation}
and
\begin{equation}\label{eq:b3-denominator-z}
	|\Omega_{2,\pm_1,\pm_2}(p,r,z)|\sim |z|^2,
	\qquad
	|\Omega_{1,\pm_3}(z,l)|\sim |z|^2.
\end{equation}

When $\max(|r|,|z|)>4\min(|r|,|z|)$, we also have
\begin{equation*}
	\bigl||r|^2-|z|^2\bigr|
	\geq\frac{15}{16}\max(|r|,|z|)^2
	>2|r-z|,
\end{equation*}
and hence
\begin{equation}\label{eq:b3-central-denominator}
	|\Omega_{1,\pm_2}(r,z)|\sim\max(|r|,|z|)^2.
\end{equation}
Thus all denominators in \eqref{eq:b3separated}
and \eqref{eq:b3comparable} are bounded away from zero
on their defining regions.

In \eqref{eq:b3highlow}, we have
\begin{equation*}
	\max(|p|,|l|)>4\min(|p|,|l|).
\end{equation*}
Together with \eqref{eq:b3-exterior-cutoff}, this gives
\begin{equation}\label{eq:b3-exterior-phase}
	|\Omega_{3,\pm_1,\pm_2,\pm_3}(p,r,z,l)|
	\gtrsim\max(|p|,|l|)^2.
\end{equation}

Under the substitution
\begin{equation*}
	(p,r,z,l;\pm_1,\pm_2,\pm_3)
	\longmapsto
	(l,z,r,p;\mp_3,\mp_2,\mp_1),
\end{equation*}
we have
\begin{equation*}
	\Omega_{1,\pm_1}(p,r)\longmapsto-\Omega_{1,\pm_3}(z,l),
	\qquad
	\Omega_{1,\pm_2}(r,z)\longmapsto-\Omega_{1,\pm_2}(r,z),
\end{equation*}
\begin{equation*}
	\Omega_{2,\pm_1,\pm_2}(p,r,z)\longmapsto-\Omega_{2,\pm_2,\pm_3}(r,z,l),
	\qquad
	\Omega_{3,\pm_1,\pm_2,\pm_3}(p,r,z,l)\longmapsto-\Omega_{3,\pm_1,\pm_2,\pm_3}(p,r,z,l).
\end{equation*}
Then one can directly check that
\begin{equation}\label{eq:b3-symmetry}
	b_{3,\pm_1,\pm_2,\pm_3}(p,r,z,l)
	=
	b_{3,\mp_3,\mp_2,\mp_1}(l,z,r,p),
\end{equation}
which satisfies the assumption (\ref{eq:coefficient-reversal}).

We next estimate $b_3$ and the resulting remainder
$\mathcal R_3$.

\begin{lemma}\label{lem:b3geometry}\label{lem:quinticremainder}
	Let $1/2<s<1$. If
	$b_{3,\pm_1,\pm_2,\pm_3}(p,r,z,l)\neq0$, then
	\begin{equation}\label{eq:b3-support}
		\max(|p|,|r|,|z|,|l|)>N,
		\qquad
		\min(|p|,|l|)\leq N/16.
	\end{equation}
	Furthermore,
	\begin{equation}\label{eq:b3bound}
		|b_{3,\pm_1,\pm_2,\pm_3}(p,r,z,l)|
		\lesssim
		N^{-2}\max(|p|,|r|,|z|,|l|)^{-2}
	\end{equation}
	on its nonzero support. For smooth solutions on $J$,
	\begin{equation}\label{eq:quinticremainder}
		\bigg|\int_J\mathcal R_3(t)\,dt\bigg|
		\lesssim_{\eps}
		\bigl(N^{-3+\eps}+N^{-2+\eps}h^{1-\eps}\bigr)A^5.
	\end{equation}
\end{lemma}

\begin{proof}
	\textbf{The bound for $b_3$.}
	
	In \eqref{eq:b3-single-r}, a nonzero coefficient requires
	$g(r)\neq0$, and hence $|r|>N$. The denominator estimates (\ref{1}) give
	\begin{equation*}
		|b_{3,\pm_1,\pm_2,\pm_3}(p,r,z,l)|
		\lesssim
		\frac{g(r)}{|r|^6}
		\lesssim |r|^{-4}.
	\end{equation*}
	Since
	\begin{equation*}
		\max(|p|,|r|,|z|,|l|)=|r|>N,
	\end{equation*}
	we obtain
	\begin{equation*}
		|b_{3,\pm_1,\pm_2,\pm_3}(p,r,z,l)|
		\lesssim
		N^{-2}\max(|p|,|r|,|z|,|l|)^{-2}.
	\end{equation*}
	By symmetry, the same estimate holds
	in the region of \eqref{eq:b3-single-z}.
	
	Consider \eqref{eq:b3separated}. By
	\eqref{eq:b3-denominator-r},
	\eqref{eq:b3-denominator-z}, and
	\eqref{eq:b3-central-denominator},
	\begin{equation*}
		|b_{3,\pm_1,\pm_2,\pm_3}(p,r,z,l)|
		\lesssim
		\frac{g(r)}{|r|^4\max(|r|,|z|)^2}
		+\frac{g(z)}{|z|^4\max(|r|,|z|)^2}
	\end{equation*}
	\begin{equation*}
		\lesssim
		\min(|r|,|z|)^{-2}\max(|r|,|z|)^{-2}.
	\end{equation*}
	In this region,
	\begin{equation*}
		\min(|r|,|z|)>N/4,
		\qquad
		\max(|p|,|r|,|z|,|l|)=\max(|r|,|z|).
	\end{equation*}
	Thus \eqref{eq:b3bound} follows.
	
	In \eqref{eq:b3comparable}, we have
	\begin{equation*}
		|r|\sim|z|,
		\qquad
		\min(|r|,|z|)>N/4.
	\end{equation*}
	Using the boundedness of $b_1$ and the denominator
	estimates, we obtain
	\begin{equation*}
		|b_{3,\pm_1,\pm_2,\pm_3}(p,r,z,l)|
		\lesssim
		\max(|r|,|z|)^{-4}
		+\bigl(g(r)+g(z)\bigr)\max(|r|,|z|)^{-6}
	\end{equation*}
	\begin{equation*}
		\lesssim
		\max(|r|,|z|)^{-4}
		\lesssim
		N^{-2}\max(|p|,|r|,|z|,|l|)^{-2}.
	\end{equation*}
	
	Finally, consider \eqref{eq:b3highlow}.
	For the fraction containing $b_{2,\pm_1,\pm_2}(p,r,z)$,
	\eqref{eq:b2bound} and \eqref{eq:b3-exterior-phase} give
	\begin{equation*}
		\left|
		\frac{b_{2,\pm_1,\pm_2}(p,r,z)}
		{\Omega_{3,\pm_1,\pm_2,\pm_3}(p,r,z,l)}
		\right|
		\lesssim
		\frac1{
			\max(|p|,|l|)^2\max(|p|,|r|,|z|)^2}.
	\end{equation*}
	If this fraction is nonzero, then
	\begin{equation*}
		\max(|p|,|l|)>N/16,
		\qquad
		\max(|p|,|r|,|z|)>N.
	\end{equation*}
	Since
	\begin{equation*}
		\max\bigl(\max(|p|,|l|),\max(|p|,|r|,|z|)\bigr)
		=\max(|p|,|r|,|z|,|l|),
	\end{equation*}
	we conclude that
	\begin{equation*}
		\left|
		\frac{b_{2,\pm_1,\pm_2}(p,r,z)}
		{\Omega_{3,\pm_1,\pm_2,\pm_3}(p,r,z,l)}
		\right|
		\lesssim
		N^{-2}\max(|p|,|r|,|z|,|l|)^{-2}.
	\end{equation*}
	For the other fraction, the same argument gives
	\begin{equation*}
		\left|
		\frac{b_{2,\pm_2,\pm_3}(r,z,l)}
		{\Omega_{3,\pm_1,\pm_2,\pm_3}(p,r,z,l)}
		\right|
		\lesssim
		\frac1{
			\max(|p|,|l|)^2\max(|r|,|z|,|l|)^2}
		\lesssim
		N^{-2}\max(|p|,|r|,|z|,|l|)^{-2}.
	\end{equation*}
	Adding these two bounds proves \eqref{eq:b3bound}.
	
	It remains to verify \eqref{eq:b3-support}.
	If $\max(|p|,|r|,|z|,|l|)\leq N$, then
	\begin{equation*}
		g(p)=g(r)=g(z)=g(l)=0,
	\end{equation*}
	and all the relevant $b_1,b_2$ coefficients vanish.
	Every defining formula for $b_3$ therefore gives zero.
	The condition $\min(|p|,|l|)\leq N/16$ follows directly
	from \eqref{eq:b3-small-exterior} and
	\eqref{eq:b3-large-exterior}.
	This proves \eqref{eq:b3-support}.

\medskip
\noindent\textbf{The regions where $R_3=0$.}

In the region of \eqref{eq:b3highlow}, substitution gives
\begin{equation*}
	R_{3,\pm_1,\pm_2,\pm_3}(p,r,z,l)
	=
	b_{2,\pm_1,\pm_2}(p,r,z)
	-b_{2,\pm_2,\pm_3}(r,z,l)
	+b_{2,\pm_2,\pm_3}(r,z,l)
	-b_{2,\pm_1,\pm_2}(p,r,z)
	=0.
\end{equation*}

Suppose next that $\max(|p|,|l|)\leq N/16$.
If $\max(|r|,|z|)\leq N/4$, then
\begin{equation*}
	b_{3,\pm_1,\pm_2,\pm_3}(p,r,z,l)
	=b_{2,\pm_1,\pm_2}(p,r,z)
	=b_{2,\pm_2,\pm_3}(r,z,l)=0,
\end{equation*}
so $R_3=0$.

Consider $|r|>N/4$ and $|z|\leq N/4$.
If $g(r)=0$, both $b_2$, $b_3$ coefficients vanish,
and hence $R_3=0$.
If $g(r)\neq0$, then $|r|>N$ and
\begin{equation*}
	b_{2,\pm_1,\pm_2}(p,r,z)
	=
	-\frac{g(r)}
	{\Omega_{1,\pm_1}(p,r)\Omega_{1,\pm_2}(r,z)}.
\end{equation*}
Moreover,
\begin{equation*}
	|r|>4|l|,
	\qquad
	\max(|r|,|z|,|l|)=|r|\leq c|r|^2.
\end{equation*}
Thus \eqref{eq:b2highlow} applies to $(r,z,l)$, giving
\begin{equation*}
	b_{2,\pm_2,\pm_3}(r,z,l)
	=
	\frac{b_{1,\pm_2}(r,z)-b_{1,\pm_3}(z,l)}
	{\Omega_{2,\pm_2,\pm_3}(r,z,l)}
	=
	\frac{g(r)}
	{\Omega_{1,\pm_2}(r,z)\Omega_{2,\pm_2,\pm_3}(r,z,l)}.
\end{equation*}
Consequently,
\begin{equation}\label{eq:b3-single-r-factorization}
	b_{2,\pm_2,\pm_3}(r,z,l)
	-b_{2,\pm_1,\pm_2}(p,r,z)
	=
	\frac{g(r)}{\Omega_{1,\pm_2}(r,z)}
	\left(
	\frac1{\Omega_{2,\pm_2,\pm_3}(r,z,l)}
	+\frac1{\Omega_{1,\pm_1}(p,r)}
	\right)
\end{equation}
\begin{equation*}
	=
	\frac{g(r)\Omega_{3,\pm_1,\pm_2,\pm_3}(p,r,z,l)}
	{\Omega_{1,\pm_1}(p,r)
		\Omega_{1,\pm_2}(r,z)
		\Omega_{2,\pm_2,\pm_3}(r,z,l)}.
\end{equation*}
By \eqref{eq:b3-single-r},
\begin{equation*}
	b_{2,\pm_2,\pm_3}(r,z,l)
	-b_{2,\pm_1,\pm_2}(p,r,z)
	=
	-\Omega_{3,\pm_1,\pm_2,\pm_3}(p,r,z,l)
	b_{3,\pm_1,\pm_2,\pm_3}(p,r,z,l),
\end{equation*}
and hence $R_3=0$.
Exchanging $p,l$ and $r,z$ and reversing the signs
proves the same assertion when $|z|>N/4$ and $|r|\leq N/4$.

We next consider 
\begin{equation*}
	\max(|p|,|l|)\leq N/16,
	\qquad
	\min(|r|,|z|)>N/4,
	\qquad
	\max(|r|,|z|)\leq c\min(|r|,|z|)^2.
\end{equation*}
Since
\begin{equation*}
	|p|\leq N/16,\qquad |z|>N/4,\qquad
	|z|>4|p|,\qquad
	\max(|p|,|r|,|z|)\leq c|z|^2,
\end{equation*}
\eqref{eq:b2highlow} gives
\begin{equation*}
	b_{2,\pm_1,\pm_2}(p,r,z)
	=
	\frac{b_{1,\pm_1}(p,r)-b_{1,\pm_2}(r,z)}
	{\Omega_{2,\pm_1,\pm_2}(p,r,z)}.
\end{equation*}
Similarly, since
\begin{equation*}
	|l|\leq N/16,\qquad |r|>N/4,\qquad
	|r|>4|l|,\qquad
	\max(|r|,|z|,|l|)\leq c|r|^2,
\end{equation*}
\eqref{eq:b2highlow} gives
\begin{equation*}
	b_{2,\pm_2,\pm_3}(r,z,l)
	=
	\frac{b_{1,\pm_2}(r,z)-b_{1,\pm_3}(z,l)}
	{\Omega_{2,\pm_2,\pm_3}(r,z,l)}.
\end{equation*}
Using
\begin{equation*}
	b_{1,\pm_1}(p,r)
	=-\frac{g(r)}{\Omega_{1,\pm_1}(p,r)},
	\qquad
	b_{1,\pm_3}(z,l)
	=\frac{g(z)}{\Omega_{1,\pm_3}(z,l)},
\end{equation*}
we obtain
\begin{equation}\label{eq:adjacentdifference}
	b_{2,\pm_2,\pm_3}(r,z,l)
	-b_{2,\pm_1,\pm_2}(p,r,z)
	=
	b_{1,\pm_2}(r,z)
	\left(
	\frac1{\Omega_{2,\pm_2,\pm_3}(r,z,l)}
	+\frac1{\Omega_{2,\pm_1,\pm_2}(p,r,z)}
	\right)
\end{equation}
\begin{equation*}
	\qquad{}
	+\frac{g(r)}
	{\Omega_{1,\pm_1}(p,r)\Omega_{2,\pm_1,\pm_2}(p,r,z)}
	-\frac{g(z)}
	{\Omega_{1,\pm_3}(z,l)\Omega_{2,\pm_2,\pm_3}(r,z,l)}.
\end{equation*}

On the region of \eqref{eq:b3separated},
\begin{equation*}
	\bigl||r|^2-|z|^2\bigr|>2|r-z|,
	\qquad
	b_{1,\pm_2}(r,z)
	=\frac{g(r)-g(z)}{\Omega_{1,\pm_2}(r,z)}.
\end{equation*}
The identities (\ref{2}) and (\ref{eq:quintic-phase-identities}) imply
\begin{equation}\label{eq:b3-separated-r-factorization}
	\frac1{\Omega_{1,\pm_2}(r,z)}
	\left(
	\frac1{\Omega_{2,\pm_2,\pm_3}(r,z,l)}
	+\frac1{\Omega_{2,\pm_1,\pm_2}(p,r,z)}
	\right)
	+\frac1{\Omega_{1,\pm_1}(p,r)\Omega_{2,\pm_1,\pm_2}(p,r,z)}
\end{equation}
\begin{equation*}
	=
	\frac{\Omega_{3,\pm_1,\pm_2,\pm_3}(p,r,z,l)}
	{\Omega_{1,\pm_1}(p,r)
		\Omega_{1,\pm_2}(r,z)
		\Omega_{2,\pm_2,\pm_3}(r,z,l)},
\end{equation*}
and
\begin{equation}\label{eq:b3-separated-z-factorization}
	-\frac1{\Omega_{1,\pm_2}(r,z)}
	\left(
	\frac1{\Omega_{2,\pm_2,\pm_3}(r,z,l)}
	+\frac1{\Omega_{2,\pm_1,\pm_2}(p,r,z)}
	\right)
	-\frac1{\Omega_{1,\pm_3}(z,l)\Omega_{2,\pm_2,\pm_3}(r,z,l)}
\end{equation}
\begin{equation*}
	=
	-\frac{\Omega_{3,\pm_1,\pm_2,\pm_3}(p,r,z,l)}
	{\Omega_{2,\pm_1,\pm_2}(p,r,z)
		\Omega_{1,\pm_2}(r,z)
		\Omega_{1,\pm_3}(z,l)}.
\end{equation*}
Multiplying these identities by $g(r)$ and $g(z)$,
respectively, and substituting into
\eqref{eq:adjacentdifference}, we obtain
\begin{equation*}
	b_{2,\pm_2,\pm_3}(r,z,l)
	-b_{2,\pm_1,\pm_2}(p,r,z)
	=
	-\Omega_{3,\pm_1,\pm_2,\pm_3}(p,r,z,l)
	b_{3,\pm_1,\pm_2,\pm_3}(p,r,z,l).
\end{equation*}
Thus $R_3=0$ on the region of \eqref{eq:b3separated}.

For \eqref{eq:b3comparable}, the same identities give
\begin{equation*}
	\frac1{\Omega_{1,\pm_1}(p,r)\Omega_{2,\pm_1,\pm_2}(p,r,z)}
	-\frac{\Omega_{3,\pm_1,\pm_2,\pm_3}(p,r,z,l)}
	{2\Omega_{1,\pm_1}(p,r)\Omega_{1,\pm_3}(z,l)}
	\left(
	\frac1{\Omega_{2,\pm_1,\pm_2}(p,r,z)}
	-\frac1{\Omega_{2,\pm_2,\pm_3}(r,z,l)}
	\right)
\end{equation*}
\begin{equation*}
	=
	\frac12
	\left(
	\frac1{\Omega_{1,\pm_1}(p,r)\Omega_{2,\pm_1,\pm_2}(p,r,z)}
	+\frac1{\Omega_{1,\pm_3}(z,l)\Omega_{2,\pm_2,\pm_3}(r,z,l)}
	\right),
\end{equation*}
\begin{equation*}
	-\frac1{\Omega_{1,\pm_3}(z,l)\Omega_{2,\pm_2,\pm_3}(r,z,l)}
	-\frac{\Omega_{3,\pm_1,\pm_2,\pm_3}(p,r,z,l)}
	{2\Omega_{1,\pm_1}(p,r)\Omega_{1,\pm_3}(z,l)}
	\left(
	\frac1{\Omega_{2,\pm_1,\pm_2}(p,r,z)}
	-\frac1{\Omega_{2,\pm_2,\pm_3}(r,z,l)}
	\right)
\end{equation*}
\begin{equation*}
	=
	-\frac12
	\left(
	\frac1{\Omega_{1,\pm_1}(p,r)\Omega_{2,\pm_1,\pm_2}(p,r,z)}
	+\frac1{\Omega_{1,\pm_3}(z,l)\Omega_{2,\pm_2,\pm_3}(r,z,l)}
	\right),
\end{equation*}
and
\begin{equation*}
	\frac1{\Omega_{2,\pm_2,\pm_3}(r,z,l)}
	+\frac1{\Omega_{2,\pm_1,\pm_2}(p,r,z)}
	-\frac{\Omega_{3,\pm_1,\pm_2,\pm_3}(p,r,z,l)}
	{2\Omega_{2,\pm_2,\pm_3}(r,z,l)\Omega_{1,\pm_1}(p,r)}
	-\frac{\Omega_{3,\pm_1,\pm_2,\pm_3}(p,r,z,l)}
	{2\Omega_{2,\pm_1,\pm_2}(p,r,z)\Omega_{1,\pm_3}(z,l)}
\end{equation*}
\begin{equation*}
	=
	-\frac{\Omega_{1,\pm_2}(r,z)}2
	\left(
	\frac1{\Omega_{1,\pm_1}(p,r)\Omega_{2,\pm_1,\pm_2}(p,r,z)}
	+\frac1{\Omega_{1,\pm_3}(z,l)\Omega_{2,\pm_2,\pm_3}(r,z,l)}
	\right).
\end{equation*}
Multiplying these identities by $g(r)$, $g(z)$, and
$b_{1,\pm_2}(r,z)$, respectively, and adding yields
\begin{equation}\label{eq:centralresidual}
	R_{3,\pm_1,\pm_2,\pm_3}(p,r,z,l)
	=
	\frac{g(r)-g(z)
		-\Omega_{1,\pm_2}(r,z)b_{1,\pm_2}(r,z)}2
\end{equation}
\begin{equation*}
	\qquad{}\times
	\left(
	\frac1{\Omega_{1,\pm_1}(p,r)\Omega_{2,\pm_1,\pm_2}(p,r,z)}
	+\frac1{\Omega_{1,\pm_3}(z,l)\Omega_{2,\pm_2,\pm_3}(r,z,l)}
	\right).
\end{equation*}
Then, if 
\begin{equation*}
	\bigl||r|^2-|z|^2\bigr|>2|r-z|
\end{equation*}
we have $R_{3,\pm_1,\pm_2,\pm_3}(p,r,z,l)=0$
on the region of \eqref{eq:b3comparable}.

If
\begin{equation*}
	\bigl||r|^2-|z|^2\bigr|\leq2|r-z|,
\end{equation*}
then \eqref{eq:b1residual} gives
\begin{equation*}
	|g(r)-g(z)
	-\Omega_{1,\pm_2}(r,z)b_{1,\pm_2}(r,z)|
	\lesssim |r-z|.
\end{equation*}
Since $|r|\sim|z|$, \eqref{eq:centralresidual} implies
\begin{equation}\label{eq:quinticcentralpointwise}
	|R_{3,\pm_1,\pm_2,\pm_3}(p,r,z,l)|
	\lesssim |r-z|\max(|r|,|z|)^{-4}.
\end{equation}

Where $b_3=0$, we have
\begin{equation}\label{eq:R3-outside}
	R_{3,\pm_1,\pm_2,\pm_3}(p,r,z,l)
	=
	b_{2,\pm_2,\pm_3}(r,z,l)
	-b_{2,\pm_1,\pm_2}(p,r,z).
\end{equation}
	
	\medskip
	\noindent\textbf{The regions where $R_3$ may be nonzero.}
	
	We first identify the frequencies that remain after the
	cancellations proved above.
	If $R_{3,\pm_1,\pm_2,\pm_3}(p,r,z,l)\neq0$, then one
	of the following four conditions holds:
	\begin{equation}\label{eq:R3-region-one}
		\max(|p|,|l|)\leq N/16,\qquad
		\min(|r|,|z|)>N/4,\qquad
		\max(|r|,|z|)>c\min(|r|,|z|)^2;
	\end{equation}
	\begin{equation}\label{eq:R3-region-two}
		\max(|p|,|l|)\leq N/16,\qquad
		\min(|r|,|z|)>N/4,\qquad
		\max(|r|,|z|)\leq c\min(|r|,|z|)^2,
	\end{equation}
	\begin{equation*}
		\max(|r|,|z|)\leq4\min(|r|,|z|),
		\qquad
		\bigl||r|^2-|z|^2\bigr|\leq2|r-z|;
	\end{equation*}
	\begin{equation}\label{eq:R3-region-three}
		\max(|p|,|l|)>N/16,\qquad
		\min(|p|,|l|)\leq N/64,
	\end{equation}
	\begin{equation*}
		\max(|p|,|r|,|z|,|l|)
		>c\max(|p|,|l|)^2;
	\end{equation*}
	\begin{equation}\label{eq:R3-region-four}
		\max(|p|,|l|)>N/16,\qquad
		\min(|p|,|l|)>N/64.
	\end{equation}
	
	To verify this inclusion, suppose first that
	$\max(|p|,|l|)\leq N/16$.
	The cancellations in \eqref{eq:b3-single-r},
	\eqref{eq:b3-single-z}, and the case
	$\max(|r|,|z|)\leq N/4$ imply
	\begin{equation*}
		R_{3,\pm_1,\pm_2,\pm_3}(p,r,z,l)\neq0
		\quad\Longrightarrow\quad
		\min(|r|,|z|)>N/4.
	\end{equation*}
	If \eqref{eq:b3-internal-cutoff} fails, we obtain
	\eqref{eq:R3-region-one}.
	If it holds, the cancellations in \eqref{eq:b3separated} and \eqref{eq:b3comparable} imply
	\begin{equation*}
		\max(|r|,|z|)\leq4\min(|r|,|z|),
		\qquad
		\bigl||r|^2-|z|^2\bigr|\leq2|r-z|,
	\end{equation*}
	which gives \eqref{eq:R3-region-two}.
	
	If $\max(|p|,|l|)>N/16$ and $\min(|p|,|l|)\leq N/64$,
	the cancellation in \eqref{eq:b3highlow} implies
	\eqref{eq:R3-region-three}.
	The remaining possibility is \eqref{eq:R3-region-four}.
	
	For all four regions, decompose
	\begin{equation*}
		u=\sum_{N_1}u_{N_1},
		\qquad
		u=\sum_{N_2}u_{N_2},
	\end{equation*}
	\begin{equation*}
		n_{\pm_1}=\sum_{N_3}n_{\pm_1,N_3},
		\qquad
		n_{\pm_2}=\sum_{N_4}n_{\pm_2,N_4},
		\qquad
		n_{\pm_3}=\sum_{N_5}n_{\pm_3,N_5},
	\end{equation*}
	where $u_{N_a}=P_{N_a}u$ and
	$n_{\pm,N_a}=P_{N_a}n_\pm$.
	In every nonzero summand,
	\begin{equation*}
		|p|\sim N_1,\qquad |l|\sim N_2,\qquad
		|p-r|\sim N_3,\qquad |r-z|\sim N_4,\qquad
		|z-l|\sim N_5.
	\end{equation*}
	For $N_a=1$, the corresponding condition is $|k|\lesssim1$.

	\medskip
	\noindent\textbf{The region \eqref{eq:R3-region-one}.}
	
	Exchanging $p,l$ and $r,z$ and reversing the signs,
	we may assume
	\begin{equation*}
		|r|>c|z|^2,\qquad |z|>N/4.
	\end{equation*}
	Since $N\geq10^4$ and $c=1/100$,
	\begin{equation*}
		|r|>4|z|,
		\qquad
		|p-r|\sim|r-z|\sim|r|,
		\qquad
		|z-l|\sim|z|.
	\end{equation*}
	Here $b_3=0$, so \eqref{eq:R3-outside} and
	\eqref{eq:b2bound} give
	\begin{equation*}
		|R_{3,\pm_1,\pm_2,\pm_3}(p,r,z,l)|
		\leq
		|b_{2,\pm_2,\pm_3}(r,z,l)|
		+|b_{2,\pm_1,\pm_2}(p,r,z)|
		\lesssim |r|^{-2}.
	\end{equation*}
	The support conditions imply
	\begin{equation*}
		N_1,N_2\leq N/8,\qquad
		N_5\geq N/8,\qquad
		N_3\gtrsim N_5^2,\qquad
		N_4\sim N_3.
	\end{equation*}
	In particular, $m(N_1)=m(N_2)=1$.
	
	For fixed signs and dyadic indices, Remark~\ref{absolute}
	and spatial H\"older's inequality give
	\begin{equation*}
		\left|
		\int_J
		\sum_{\substack{p,r,z,l\\
				\max(|p|,|l|)\leq N/16\\
				|z|>N/4,\ |r|>c|z|^2}}
		R_{3,\pm_1,\pm_2,\pm_3}(p,r,z,l)
		\overline{\widehat{u_{N_1}}(t,p)}
		\widehat{u_{N_2}}(t,l)
		\widehat{n_{\pm_1,N_3}}(t,p-r)
		\widehat{n_{\pm_2,N_4}}(t,r-z)
		\widehat{n_{\pm_3,N_5}}(t,z-l)\,dt
		\right|
	\end{equation*}
	\begin{equation*}
		\lesssim N_3^{-2}\int_J
		\left\|\mathcal F_x^{-1}|\widehat{u_{N_1}}(t,\cdot)|\right\|_{L_x^\infty}
		\left\|\mathcal F_x^{-1}|\widehat{u_{N_2}}(t,\cdot)|\right\|_{L_x^\infty}
		\|n_{\pm_1,N_3}(t)\|_{L_x^2}
		\|n_{\pm_2,N_4}(t)\|_{L_x^2}
		\left\|\mathcal F_x^{-1}|\widehat{n_{\pm_3,N_5}}(t,\cdot)|\right\|_{L_x^\infty}\,dt
	\end{equation*}
	\begin{equation*}
		\lesssim
		hN_3^{-2}N_5
		A_{N_1}A_{N_2}B_{N_3}B_{N_4}B_{N_5},
	\end{equation*}
	where we used \eqref{eq:spatialu} and
	\eqref{eq:spatialwave}.
	
	By Cauchy--Schwarz,
	\begin{equation*}
		\sum_{N_1,N_2\leq N/8}A_{N_1}A_{N_2}
		\lesssim
		(1+\log(2+N))\sum_{N_1}A_{N_1}^2
		\lesssim_\eps N^\eps A^2.
	\end{equation*}
	Also,
	\begin{equation*}
		B_{N_5}\lesssim A,
		\qquad
		\sum_{N_4\sim N_3}B_{N_3}B_{N_4}\lesssim A^2,
		\qquad
		\sum_{N_3\gtrsim N_5^2}N_3^{-2}\lesssim N_5^{-4}.
	\end{equation*}
	Summing these estimates and including the symmetric region,
	we obtain
	\begin{equation*}
		\left|
		\int_J
		\sum_{\substack{p,r,z,l\\\text{\eqref{eq:R3-region-one}}}}
		\sum_{\pm_1,\pm_2,\pm_3}
		R_{3,\pm_1,\pm_2,\pm_3}(p,r,z,l)
		\overline{\widehat u(t,p)}\widehat u(t,l)
		\widehat{n_{\pm_1}}(t,p-r)
		\widehat{n_{\pm_2}}(t,r-z)
		\widehat{n_{\pm_3}}(t,z-l)\,dt
		\right|
	\end{equation*}
	\begin{equation*}
		\lesssim_{\eps}
		hN^\eps A^5
		\sum_{N_5\geq N/8}N_5
		\sum_{N_3\gtrsim N_5^2}N_3^{-2}
		\lesssim
		hN^\eps A^5\sum_{N_5\geq N/8}N_5^{-3}
		\lesssim hN^{-3+\eps}A^5.
	\end{equation*}
	
	\medskip
	\noindent\textbf{The region \eqref{eq:R3-region-two}.}
	
	The restrictions imply
	\begin{equation*}
		|r|\sim|z|,
		\qquad
		|p-r|\sim|r|,
		\qquad
		|z-l|\sim|z|.
	\end{equation*}
	Consequently,
	\begin{equation*}
		N_1,N_2\leq N/8,\qquad
		N_3\geq N/8,\qquad
		N_5\sim N_3,\qquad
		N_4\lesssim N_3.
	\end{equation*}
	By \eqref{eq:quinticcentralpointwise},
	\begin{equation*}
		|R_{3,\pm_1,\pm_2,\pm_3}(p,r,z,l)|
		\lesssim N_4N_3^{-4}.
	\end{equation*}
	For fixed signs and dyadic indices,
	\begin{equation*}
		\left|
		\int_J
		\sum_{\substack{p,r,z,l\\\text{\eqref{eq:R3-region-two}}}}
		R_{3,\pm_1,\pm_2,\pm_3}(p,r,z,l)
		\overline{\widehat{u_{N_1}}(t,p)}
		\widehat{u_{N_2}}(t,l)
		\widehat{n_{\pm_1,N_3}}(t,p-r)
		\widehat{n_{\pm_2,N_4}}(t,r-z)
		\widehat{n_{\pm_3,N_5}}(t,z-l)\,dt
		\right|
	\end{equation*}
	\begin{equation*}
		\lesssim N_4N_3^{-4}\int_J
		\left\|\mathcal F_x^{-1}|\widehat{u_{N_1}}(t,\cdot)|\right\|_{L_x^\infty}
		\left\|\mathcal F_x^{-1}|\widehat{u_{N_2}}(t,\cdot)|\right\|_{L_x^\infty}
		\|n_{\pm_1,N_3}(t)\|_{L_x^2}
		\left\|\mathcal F_x^{-1}|\widehat{n_{\pm_2,N_4}}(t,\cdot)|\right\|_{L_x^\infty}
		\|n_{\pm_3,N_5}(t)\|_{L_x^2}\,dt
	\end{equation*}
	\begin{equation*}
		\lesssim
		hN_4^2N_3^{-4}
		A_{N_1}A_{N_2}B_{N_3}B_{N_4}B_{N_5}.
	\end{equation*}
	Using
	\begin{equation*}
		\sum_{N_1,N_2\leq N/8}A_{N_1}A_{N_2}
		\lesssim_\eps N^\eps A^2,
		\qquad
		B_{N_4}\lesssim A,
		\qquad
		\sum_{N_5\sim N_3}B_{N_3}B_{N_5}\lesssim A^2,
	\end{equation*}
	we obtain
	\begin{equation*}
		\left|
		\int_J
		\sum_{\substack{p,r,z,l\\\text{\eqref{eq:R3-region-two}}}}
		\sum_{\pm_1,\pm_2,\pm_3}
		R_{3,\pm_1,\pm_2,\pm_3}(p,r,z,l)
		\overline{\widehat u(t,p)}\widehat u(t,l)
		\widehat{n_{\pm_1}}(t,p-r)
		\widehat{n_{\pm_2}}(t,r-z)
		\widehat{n_{\pm_3}}(t,z-l)\,dt
		\right|
	\end{equation*}
	\begin{equation*}
		\lesssim_{\eps}
		hN^\eps A^5
		\sum_{N_3\geq N/8}N_3^{-4}
		\sum_{N_4\lesssim N_3}N_4^2\lesssim hN^{-2+\eps}A^5.
	\end{equation*}

	\medskip
	\noindent\textbf{The region \eqref{eq:R3-region-three}.}
	
	Exchanging $p,l$ and $r,z$ and reversing the signs,
	we may assume
	\begin{equation*}
		|p|>N/16,\qquad |l|\leq N/64,\qquad
		\max(|p|,|r|,|z|,|l|)>c|p|^2.
	\end{equation*}
	Since $c|p|>4$, we have
	\begin{equation*}
		\max(|r|,|z|)>c|p|^2>4|p|.
	\end{equation*}
    We first restrict the sum below to
     \begin{equation*}
	\max(|r|,|z|)\sim N_6,
     \end{equation*}
     and then sum over dyadic $N_6$.
	Then
	\begin{equation*}
		N_1\geq N/32,\qquad
		N_2\leq N/32,\qquad
		N_6\gtrsim N_1^2,\qquad
		N_3,N_4,N_5\lesssim N_6,
	\end{equation*}
	and $m(N_2)=1$.
	Since both triples $(p,r,z)$ and $(r,z,l)$ have maximum
	magnitude comparable to $N_6$, \eqref{eq:R3-outside}
	and \eqref{eq:b2bound} give
	\begin{equation*}
		|R_{3,\pm_1,\pm_2,\pm_3}(p,r,z,l)|
		\lesssim N_6^{-2}.
	\end{equation*}
	
	For fixed signs and dyadic indices, first suppose
	$N_3=\min(N_3,N_4,N_5)$.
	Then
	\begin{equation*}
		\left|
		\int_J
		\sum_{\substack{p,r,z,l\\
				|p|>N/16,\ |l|\leq N/64\\
				\max(|r|,|z|)>c|p|^2\\
				\max(|r|,|z|)\sim N_6}}
		R_{3,\pm_1,\pm_2,\pm_3}(p,r,z,l)
		\overline{\widehat{u_{N_1}}(t,p)}
		\widehat{u_{N_2}}(t,l)
		\widehat{n_{\pm_1,N_3}}(t,p-r)
		\widehat{n_{\pm_2,N_4}}(t,r-z)
		\widehat{n_{\pm_3,N_5}}(t,z-l)\,dt
		\right|
	\end{equation*}
	\begin{equation*}
		\lesssim N_6^{-2}\int_J
		\left\|\mathcal F_x^{-1}|\widehat{u_{N_1}}(t,\cdot)|\right\|_{L_x^\infty}
		\left\|\mathcal F_x^{-1}|\widehat{u_{N_2}}(t,\cdot)|\right\|_{L_x^\infty}
		\left\|\mathcal F_x^{-1}|\widehat{n_{\pm_1,N_3}}(t,\cdot)|\right\|_{L_x^\infty}
		\|n_{\pm_2,N_4}(t)\|_{L_x^2}
		\|n_{\pm_3,N_5}(t)\|_{L_x^2}\,dt
	\end{equation*}
	\begin{equation*}
		\lesssim
		hN_6^{-2}N_3m(N_1)^{-1}
		A_{N_1}A_{N_2}B_{N_3}B_{N_4}B_{N_5}.
	\end{equation*}
	If $N_4$ or $N_5$ is the smallest index, use
	$L_x^\infty$ for the corresponding $n$ instead.
	Therefore, for every ordering of $N_3,N_4,N_5$,
	\begin{equation*}
		\left|
		\int_J
		\sum_{\substack{p,r,z,l\\
				|p|>N/16,\ |l|\leq N/64\\
				\max(|r|,|z|)>c|p|^2\\
				\max(|r|,|z|)\sim N_6}}
		R_{3,\pm_1,\pm_2,\pm_3}(p,r,z,l)
		\overline{\widehat{u_{N_1}}(t,p)}
		\widehat{u_{N_2}}(t,l)
		\widehat{n_{\pm_1,N_3}}(t,p-r)
		\widehat{n_{\pm_2,N_4}}(t,r-z)
		\widehat{n_{\pm_3,N_5}}(t,z-l)\,dt
		\right|
	\end{equation*}
	\begin{equation*}
		\lesssim
		hN_6^{-2}\min(N_3,N_4,N_5)m(N_1)^{-1}
		A_{N_1}A_{N_2}B_{N_3}B_{N_4}B_{N_5}.
	\end{equation*}
	
	Using $B_{N_a}\lesssim A$ and counting the dyadic choices,
	\begin{equation*}
		\sum_{N_3,N_4,N_5\lesssim N_6}
		\min(N_3,N_4,N_5)B_{N_3}B_{N_4}B_{N_5}
		\lesssim
		N_6\bigl(1+\log(2+N_6)\bigr)^3A^3
		\lesssim_\eps N_6^{1+\eps/8}A^3.
	\end{equation*}
	Also,
	\begin{equation*}
		\sum_{N_2\leq N/32}A_{N_2}
		\lesssim
		(1+\log(2+N))^{1/2}
		\left(\sum_{N_2}A_{N_2}^2\right)^{1/2}
		\lesssim_\eps N^{\eps/8}A.
	\end{equation*}
	After summing and including the symmetric region, we obtain
	\begin{equation*}
		\left|
		\int_J
		\sum_{\substack{p,r,z,l\\\text{\eqref{eq:R3-region-three}}}}
		\sum_{\pm_1,\pm_2,\pm_3}
		R_{3,\pm_1,\pm_2,\pm_3}(p,r,z,l)
		\overline{\widehat u(t,p)}\widehat u(t,l)
		\widehat{n_{\pm_1}}(t,p-r)
		\widehat{n_{\pm_2}}(t,r-z)
		\widehat{n_{\pm_3}}(t,z-l)\,dt
		\right|
	\end{equation*}
	\begin{equation*}
		\lesssim_{\eps}
		hN^{\eps/8}A^5
		\sum_{N_1\geq N/32}m(N_1)^{-1}
		\sum_{N_6\gtrsim N_1^2}N_6^{-1+\eps/8}\lesssim hN^{-2+\eps}A^5.
	\end{equation*}

	\medskip
	\noindent\textbf{The region \eqref{eq:R3-region-four}.}
	
	Here $b_3=0$, and hence
	\begin{equation*}
		R_{3,\pm_1,\pm_2,\pm_3}(p,r,z,l)
		=
		b_{2,\pm_2,\pm_3}(r,z,l)
		-b_{2,\pm_1,\pm_2}(p,r,z).
	\end{equation*}
	We first estimate the integral containing
	$b_{2,\pm_2,\pm_3}(r,z,l)$.
	Partition its sum according to
	\begin{equation*}
		\max(|r|,|z|,|l|)\sim N_6.
	\end{equation*}
	By Lemma~\ref{lem:b2}, a nonzero coefficient $b_2$ satisfies
	\begin{equation*}
		N_6\gtrsim N,
		\qquad
		|b_{2,\pm_2,\pm_3}(r,z,l)|\lesssim N_6^{-2}.
	\end{equation*}
	The support conditions give
	\begin{equation*}
		N_1,N_2\geq N/128,\qquad
		N_4,N_5\lesssim N_6,\qquad
		N_3\lesssim\max(N_1,N_6).
	\end{equation*}
	
	Insert $\psi_h$ on both occurrences of $u$ and on
	$n_{\pm_3,N_5}$, and insert $1_J$ on
	$n_{\pm_1,N_3}$ and $n_{\pm_2,N_4}$.
	The integral is unchanged because $\psi_h^3 1_J^2=1_J.$

	By Remark~\ref{absolute} and H\"older's inequality,
	\begin{equation*}
		\left|
		\int_J
		\sum_{\substack{p,r,z,l\\
				\text{\eqref{eq:R3-region-four}}\\
				\max(|r|,|z|,|l|)\sim N_6}}
		b_{2,\pm_2,\pm_3}(r,z,l)
		\overline{\widehat{u_{N_1}}(t,p)}
		\widehat{u_{N_2}}(t,l)
		\widehat{n_{\pm_1,N_3}}(t,p-r)
		\widehat{n_{\pm_2,N_4}}(t,r-z)
		\widehat{n_{\pm_3,N_5}}(t,z-l)\,dt
		\right|
	\end{equation*}
	\begin{equation*}
		\lesssim N_6^{-2}
		\left\|
		\left(\mathcal F_{t,x}^{-1}
		\left|\widetilde{\psi_h\overline{u_{N_1}}}\right|\right)
		\left(\mathcal F_{t,x}^{-1}
		\left|\widetilde{1_Jn_{\pm_1,N_3}}\right|\right)
		\right\|_{L^2_{t,x}}
		\left\|
		\left(\mathcal F_{t,x}^{-1}
		\left|\widetilde{\psi_hu_{N_2}}\right|\right)
		\left(\mathcal F_{t,x}^{-1}
		\left|\widetilde{1_Jn_{\pm_2,N_4}}\right|\right)
		\right\|_{L^2_{t,x}}
		\left\|
		\mathcal F_{t,x}^{-1}
		\left|\widetilde{\psi_hn_{\pm_3,N_5}}\right|
		\right\|_{L^\infty_{t,x}}.
	\end{equation*}
	Bernstein's inequality and
	\eqref{eq:smoothcutoff} give
	\begin{equation*}
		\left\|
		\mathcal F_{t,x}^{-1}
		\left|\widetilde{\psi_hn_{\pm_3,N_5}}\right|
		\right\|_{L^\infty_{t,x}}
		\lesssim N_5
		\left\|
		\mathcal F_{t,x}^{-1}
		\left|\widetilde{\psi_hn_{\pm_3,N_5}}\right|
		\right\|_{L^\infty_tL_x^2}
	\end{equation*}
	\begin{equation*}
		\lesssim N_5
		\left\|
		\mathcal F_{t,x}^{-1}
		\left|\widetilde{\psi_hn_{\pm_3,N_5}}\right|
		\right\|_{X_{W_{\pm_3}}^{0,1/2,1}}
		=
		N_5\|\psi_hn_{\pm_3,N_5}\|_{X_{W_{\pm_3}}^{0,1/2,1}}
		\lesssim
		N_5\|n_{\pm_3,N_5}\|_{X_{W_{\pm_3}}^{0,1/2,1}}
		=N_5B_{N_5}.
	\end{equation*}
	Applying \eqref{eq:pair} with $\eps/8$ to the two
	$L^2_{t,x}$ norms, and using $N_5\lesssim N_6$, yields
	\begin{equation*}
		\left|
		\int_J
		\sum_{\substack{p,r,z,l\\
				\text{\eqref{eq:R3-region-four}}\\
				\max(|r|,|z|,|l|)\sim N_6}}
		b_{2,\pm_2,\pm_3}(r,z,l)
		\overline{\widehat{u_{N_1}}(t,p)}
		\widehat{u_{N_2}}(t,l)
		\widehat{n_{\pm_1,N_3}}(t,p-r)
		\widehat{n_{\pm_2,N_4}}(t,r-z)
		\widehat{n_{\pm_3,N_5}}(t,z-l)\,dt
		\right|
	\end{equation*}
	\begin{equation*}
		\lesssim_{\eps}
		N_6^{-1}m(N_1)^{-1}m(N_2)^{-1}
		A_{N_1}A_{N_2}B_{N_3}B_{N_4}B_{N_5}
	\end{equation*}
	\begin{equation*}
		\qquad{}\times
		\bigl(N_1^{-1+\eps/8}
		+h^{1/2-\eps/8}N_1^{-1/2+\eps/8}\bigr)
		\bigl(N_2^{-1+\eps/8}
		+h^{1/2-\eps/8}N_2^{-1/2+\eps/8}\bigr).
	\end{equation*}
	
	We first sum in $N_3,N_4,N_5$.
	Since $B_{N_a}\lesssim A$,
	\begin{equation*}
		\sum_{\substack{
				N_3\lesssim\max(N_1,N_6)\\
				N_4,N_5\lesssim N_6}}
		B_{N_3}B_{N_4}B_{N_5}
		\lesssim_\eps (N_1N_6)^{\eps/8}A^3.
	\end{equation*}
	Using also $A_{N_1},A_{N_2}\lesssim A$, we obtain
	\begin{equation}\label{eq:R3-fourth-region-bound}
		\left|
		\int_J
		\sum_{\substack{p,r,z,l\\\text{\eqref{eq:R3-region-four}}}}
		\sum_{\pm_1,\pm_2,\pm_3}
		b_{2,\pm_2,\pm_3}(r,z,l)
		\overline{\widehat u(t,p)}\widehat u(t,l)
		\widehat{n_{\pm_1}}(t,p-r)
		\widehat{n_{\pm_2}}(t,r-z)
		\widehat{n_{\pm_3}}(t,z-l)\,dt
		\right|
	\end{equation}
	\begin{equation*}
		\lesssim_{\eps}
		A^5\sum_{N_6\gtrsim N}N_6^{-1+\eps/8}
		\left[
		\sum_{N_1\geq N/128}
		m(N_1)^{-1}
		\bigl(
		N_1^{-1+\eps/4}
		+h^{1/2-\eps/8}N_1^{-1/2+\eps/4}
		\bigr)
		\right]^2.
	\end{equation*}
	
	For $N_1\geq N/128$,
	\begin{equation*}
		m(N_1)^{-1}\sim N^{s-1}N_1^{1-s}.
	\end{equation*}
	Choose $0<\eps<\min(1,4(s-1/2))$.
	Then
	\begin{equation*}
		\sum_{N_1\geq N/128}
		m(N_1)^{-1}N_1^{-1+\eps/4}
		\lesssim
		N^{s-1}\sum_{N_1\geq N/128}N_1^{-s+\eps/4}
		\lesssim_{\eps}N^{-1+\eps/4},
	\end{equation*}
	and
	\begin{equation*}
		\sum_{N_1\geq N/128}
		m(N_1)^{-1}N_1^{-1/2+\eps/4}
		\lesssim
		N^{s-1}\sum_{N_1\geq N/128}N_1^{1/2-s+\eps/4}
		\lesssim_{\eps}N^{-1/2+\eps/4}.
	\end{equation*}
	Also,
	\begin{equation*}
		\sum_{N_6\gtrsim N}N_6^{-1+\eps/8}
		\lesssim_\eps N^{-1+\eps/8}.
	\end{equation*}
	Substitution into \eqref{eq:R3-fourth-region-bound}
	gives
	\begin{equation*}
		\left|
		\int_J
		\sum_{\substack{p,r,z,l\\\text{\eqref{eq:R3-region-four}}}}
		\sum_{\pm_1,\pm_2,\pm_3}
		b_{2,\pm_2,\pm_3}(r,z,l)
		\overline{\widehat u(t,p)}\widehat u(t,l)
		\widehat{n_{\pm_1}}(t,p-r)
		\widehat{n_{\pm_2}}(t,r-z)
		\widehat{n_{\pm_3}}(t,z-l)\,dt
		\right|
	\end{equation*}
	\begin{equation*}
		\lesssim_{\eps}
		N^{-1+\eps/8}
		\bigl(N^{-1+\eps/4}
		+h^{1/2-\eps/8}N^{-1/2+\eps/4}\bigr)^2A^5
		\lesssim
		\bigl(N^{-3+\eps}+N^{-2+\eps}h^{1-\eps}\bigr)A^5.
	\end{equation*}
	Exchanging $p,l$ and $r,z$ and reversing the signs
	gives the same bound for the integral containing
	$b_{2,\pm_1,\pm_2}(p,r,z)$.
	Therefore
	\begin{equation*}
		\left|
		\int_J
		\sum_{\substack{p,r,z,l\\\text{\eqref{eq:R3-region-four}}}}
		\sum_{\pm_1,\pm_2,\pm_3}
		R_{3,\pm_1,\pm_2,\pm_3}(p,r,z,l)
		\overline{\widehat u(t,p)}\widehat u(t,l)
		\widehat{n_{\pm_1}}(t,p-r)
		\widehat{n_{\pm_2}}(t,r-z)
		\widehat{n_{\pm_3}}(t,z-l)\,dt
		\right|
	\end{equation*}
	\begin{equation*}
		\lesssim_{\eps}
		\bigl(N^{-3+\eps}+N^{-2+\eps}h^{1-\eps}\bigr)A^5.
	\end{equation*}
	
	Since $R_3$ vanishes outside the four regions listed above,
	the definition of $\mathcal R_3$ gives
	\begin{equation*}
		\bigg|\int_J\mathcal R_3(t)\,dt\bigg|
		\lesssim_{\eps}
		\bigl(
		hN^{-3+\eps}+hN^{-2+\eps}
		+N^{-3+\eps}+N^{-2+\eps}h^{1-\eps}
		\bigr)A^5
		\lesssim
		\bigl(N^{-3+\eps}+N^{-2+\eps}h^{1-\eps}\bigr)A^5,
	\end{equation*}
	which proves \eqref{eq:quinticremainder}.
	
\end{proof}

\vspace{10pt}
\section{Fixed-time bounds for the additional corrections}
We now complete the proof of the desired fixed-time estimate.

\begin{proof}[Proof of Proposition \ref{prop:comparison}]
	
	We first estimate $\mathcal C_4$ and $\mathcal C_5$ at a fixed
	time and then compare $\mathcal E_N(u,n_+)$ with
	$\mathcal H(Iu,n_+)$.
	
	\textbf{The bound for $\mathcal C_4$.}
	
	By \eqref{eq:form},
	\begin{equation*}
		\mathcal C_4(u,n_+)
		=
		\frac{1}{4}
		\sum_{p,r,l}\sum_{\pm_1,\pm_2}
		b_{2,\pm_1,\pm_2}(p,r,l)
		\overline{\widehat u(p)}\widehat u(l)
		\widehat{n_{\pm_1}}(p-r)
		\widehat{n_{\pm_2}}(r-l).
	\end{equation*}
	Decompose
	\begin{equation*}
		u=\sum_{N_1}u_{N_1},
		\qquad
		u=\sum_{N_2}u_{N_2},
		\qquad
		n_{\pm_1}=\sum_{N_3}n_{\pm_1,N_3},
		\qquad
		n_{\pm_2}=\sum_{N_4}n_{\pm_2,N_4}.
	\end{equation*}
	Then the spatial Fourier supports satisfy
	\begin{equation*}
		|p|\sim N_1,\qquad |l|\sim N_2,\qquad
		|p-r|\sim N_3,\qquad |r-l|\sim N_4.
	\end{equation*}
	
	We first restrict the sum to
	\begin{equation*}
		\max(|p|,|r|,|l|)\sim N_5,
	\end{equation*}
	and then sum over dyadic $N_5$.
	On the nonzero support of $b_2$, Lemma~\ref{lem:b2} gives
	\begin{equation*}
		N_5\ge N/2,
		\qquad
		|b_{2,\pm_1,\pm_2}(p,r,l)|\lesssim N_5^{-2},
		\qquad
		\min(|p|,|l|)\leq N/4.
	\end{equation*}
	The triangle inequality also gives
	\begin{equation*}
		N_1,N_2,N_3,N_4\lesssim N_5.
	\end{equation*}
	Since at least one of $|p|,|l|$ is at most $N/4$,
	at least one of $N_1,N_2$ is at most $N/2$. Hence
	\begin{equation}\label{eq:fixed-u-weights}
		m(N_1)^{-1}m(N_2)^{-1}
		\lesssim m(N_5)^{-1}.
	\end{equation}
	
	 Bernstein's inequality and Plancherel's identity give
	\begin{equation}\label{eq:fixed-u-infinity}
		\left\|
		\mathcal F_x^{-1}|\widehat{u_{N_a}}|
		\right\|_{L_x^\infty}
		\lesssim N_a\|u_{N_a}\|_{L_x^2}
		\lesssim m(N_a)^{-1}\|Iu_{N_a}\|_{H^1}.
	\end{equation}
	Similarly,
	\begin{equation}\label{eq:fixed-n-norms}
		\left\|
		\mathcal F_x^{-1}|\widehat{n_{\pm,N_a}}|
		\right\|_{L_x^2}
		=\|n_{\pm,N_a}\|_{L_x^2},
		\qquad
		\left\|
		\mathcal F_x^{-1}|\widehat{n_{\pm,N_a}}|
		\right\|_{L_x^\infty}
		\lesssim N_a\|n_{\pm,N_a}\|_{L_x^2}.
	\end{equation}
	
	Then we can obtain
	\begin{equation*}
		\left|
		\sum_{\substack{p,r,l\\
				\max(|p|,|r|,|l|)\sim N_5}}
		b_{2,\pm_1,\pm_2}(p,r,l)
		\overline{\widehat{u_{N_1}}(p)}
		\widehat{u_{N_2}}(l)
		\widehat{n_{\pm_1,N_3}}(p-r)
		\widehat{n_{\pm_2,N_4}}(r-l)
		\right|
	\end{equation*}
	\begin{equation*}
		\lesssim N_5^{-2}
		\left\|\mathcal F_x^{-1}|\widehat{u_{N_1}}|\right\|_{L_x^\infty}
		\left\|\mathcal F_x^{-1}|\widehat{u_{N_2}}|\right\|_{L_x^\infty}
		\|n_{\pm_1,N_3}\|_{L_x^2}
		\|n_{\pm_2,N_4}\|_{L_x^2}
	\end{equation*}
	\begin{equation*}
		\lesssim N_5^{-2}m(N_5)^{-1}
		\|Iu_{N_1}\|_{H^1}\|Iu_{N_2}\|_{H^1}
		\|n_{\pm_1,N_3}\|_{L_x^2}
		\|n_{\pm_2,N_4}\|_{L_x^2}
		\lesssim N_5^{-2}m(N_5)^{-1}a^4.
	\end{equation*}
	Here we used \eqref{eq:fixed-u-weights}--%
	\eqref{eq:fixed-n-norms}.
	
	For fixed $N_5$, the number of possible dyadic
	$(N_1,N_2,N_3,N_4)$ is at most $N_5^\eps$.
	Therefore
	\begin{equation*}
		|\mathcal C_4(u,n_+)|
		\lesssim_\eps
		a^4\sum_{N_5\geq N/2}
		N_5^{-2+\eps}m(N_5)^{-1}\lesssim N^{-2+\eps}a^4.
	\end{equation*}

	\medskip
	\noindent\textbf{The bound for $\mathcal C_5$.}
	
	By \eqref{eq:form},
	\begin{equation*}
		\mathcal C_5(u,n_+)
		=
		\frac{1}{8}
		\sum_{p,r,z,l}\sum_{\pm_1,\pm_2,\pm_3}
		b_{3,\pm_1,\pm_2,\pm_3}(p,r,z,l)
		\overline{\widehat u(p)}\widehat u(l)
		\widehat{n_{\pm_1}}(p-r)
		\widehat{n_{\pm_2}}(r-z)
		\widehat{n_{\pm_3}}(z-l).
	\end{equation*}
	Dyadically decomposing the spatial frequency as before, we have
	\begin{equation*}
		|p|\sim N_1,\qquad |l|\sim N_2,\qquad
		|p-r|\sim N_3,\qquad |r-z|\sim N_4,\qquad
		|z-l|\sim N_5.
	\end{equation*}
	And one first restricts the sum to
	\begin{equation*}
		\max(|p|,|r|,|z|,|l|)\sim N_6,
	\end{equation*}
	and then sums over dyadic $N_6$.
	
	By \eqref{eq:b3-support} and \eqref{eq:b3bound},
	every nonzero summand satisfies
	\begin{equation*}
		N_6\geq N/2,
		\qquad
		|b_{3,\pm_1,\pm_2,\pm_3}(p,r,z,l)|
		\lesssim N^{-2}N_6^{-2},
		\qquad
		\min(|p|,|l|)\leq N/16.
	\end{equation*}
	Furthermore,
	\begin{equation*}
		N_1,N_2,N_3,N_4,N_5\lesssim N_6.
	\end{equation*}
	At least one of $N_1,N_2$ is at most $N/8$, so
	\begin{equation*}
		m(N_1)^{-1}m(N_2)^{-1}
		\lesssim m(N_6)^{-1}.
	\end{equation*}
	Then we can derive
	\begin{equation*}
		\left|
		\sum_{\substack{p,r,z,l\\
				\max(|p|,|r|,|z|,|l|)\sim N_6}}
		b_{3,\pm_1,\pm_2,\pm_3}(p,r,z,l)
		\overline{\widehat{u_{N_1}}(p)}
		\widehat{u_{N_2}}(l)
		\widehat{n_{\pm_1,N_3}}(p-r)
		\widehat{n_{\pm_2,N_4}}(r-z)
		\widehat{n_{\pm_3,N_5}}(z-l)
		\right|
	\end{equation*}
	\begin{equation*}
		\lesssim N^{-2}N_6^{-2}
		\left\|\mathcal F_x^{-1}|\widehat{u_{N_1}}|\right\|_{L_x^\infty}
		\left\|\mathcal F_x^{-1}|\widehat{u_{N_2}}|\right\|_{L_x^\infty}
		\|n_{\pm_1,N_3}\|_{L_x^2}
		\|n_{\pm_2,N_4}\|_{L_x^2}
		\left\|\mathcal F_x^{-1}|\widehat{n_{\pm_3,N_5}}|\right\|_{L_x^\infty}
	\end{equation*}
	\begin{equation*}
		\lesssim
		N^{-2}N_6^{-2}N_5m(N_6)^{-1}
		\|Iu_{N_1}\|_{H^1}\|Iu_{N_2}\|_{H^1}
		\|n_{\pm_1,N_3}\|_{L_x^2}
		\|n_{\pm_2,N_4}\|_{L_x^2}
		\|n_{\pm_3,N_5}\|_{L_x^2}
	\end{equation*}
	\begin{equation*}
		\lesssim
		N^{-2}N_6^{-1}m(N_6)^{-1}a^5.
	\end{equation*}
	Here we used \eqref{eq:fixed-u-infinity},
	\eqref{eq:fixed-n-norms}, and $N_5\lesssim N_6$.
	
	For fixed $N_6$, the number of possible dyadic
	$(N_1,N_2,N_3,N_4,N_5)$ is at most $N_6^\eps$.
	Thus
	\begin{equation*}
		|\mathcal C_5(u,n_+)|
		\lesssim_\eps
		N^{-2}a^5
		\sum_{N_6\geq N/2}N_6^{-1+\eps}m(N_6)^{-1}\lesssim N^{-3+\eps}a^5.
	\end{equation*}

	\medskip
	\noindent\textbf{The comparison of the energies.}
	
	By \eqref{eq:energydef},
	\begin{equation*}
		\mathcal E_N(u,n_+)-\mathcal H(Iu,n_+)
		=
		\widetilde{\mathcal H}(u,n_+)-\mathcal H(Iu,n_+)
		-\mathcal C_4(u,n_+)-\mathcal C_5(u,n_+).
	\end{equation*}
	Using \eqref{eq:cubiccomparison} and the estimates above,
	\begin{equation*}
		|\mathcal E_N(u,n_+)-\mathcal H(Iu,n_+)|
		\leq
		|\widetilde{\mathcal H}(u,n_+)-\mathcal H(Iu,n_+)|
		+|\mathcal C_4(u,n_+)|+|\mathcal C_5(u,n_+)|
	\end{equation*}
	\begin{equation*}
		\lesssim_\eps
		N^{-1+\eps}a^3+N^{-2+\eps}a^4+N^{-3+\eps}a^5.
	\end{equation*}
	This proves \eqref{eq:comparison}.

\end{proof}

\vspace{10pt}
\section{Estimates for $\mathcal F_2$ and $\mathcal F_3$}

We estimate the expressions $\mathcal F_2$ and $\mathcal F_3$
in \eqref{eq:energyidentity}. Recall that
\begin{equation}\label{eq:Fj-explicit}
	\mathcal F_j
	=
	-\frac{i}{2^j}
	\sum_{p_0,\ldots,p_j}
	\sum_{\pm_1,\ldots,\pm_j}
	b_{j,\pm_1,\ldots,\pm_j}(p_0,\ldots,p_j)
	\overline{\widehat u(p_0)}\widehat u(p_j)
	\sum_{a=1}^j
	\pm_a
	\widehat{D(|u|^2)}(p_{a-1}-p_a)
	\prod_{\substack{1\leq b\leq j\\b\neq a}}
	\widehat{n_{\pm_b}}(p_{b-1}-p_b).
\end{equation}
In particular, $\mathcal F_2$ contains four occurrences
of $u$ and one of $n$, while $\mathcal F_3$ contains
four occurrences of $u$ and two of $n$.

\begin{lemma}\label{lem:waveforcing}
	Let $1/2<s<1$. For smooth solutions on $J$, with $h=|J|$,
	\begin{equation}\label{eq:F2}
		\bigg|\int_J\mathcal F_2(t)\,dt\bigg|
		\lesssim_\eps hN^{-2+\eps}A^5,
	\end{equation}
	and
	\begin{equation}\label{eq:F3}
		\bigg|\int_J\mathcal F_3(t)\,dt\bigg|
		\lesssim_\eps hN^{-3+\eps}A^6.
	\end{equation}
\end{lemma}

\begin{proof}
	\textbf{The estimate for $D(|u|^2)$.}
	
	We first prove that, for every dyadic $N_3\geq1$,
	\begin{equation}\label{eq:newpair}
		\left\|
		\mathcal F_x^{-1}
		\left|\widehat{P_{N_3}D(|u|^2)}\right|
		\right\|_{L_t^\infty L_x^2}
		\lesssim_\eps
		N_3^{\eps/8}m(N_3)^{-2}A^2.
	\end{equation}

 Decomposing the two occurrences of $u$ gives
\begin{equation*}
	P_{N_3}(|u|^2)
	=
	\sum_{N_1,N_2}
	P_{N_3}(u_{N_1}\overline{u_{N_2}}).
\end{equation*}
The Fourier coefficient of the product is
\begin{equation*}
	\widehat{u_{N_1}\overline{u_{N_2}}}(q)
	=
	\sum_k
	\widehat{u_{N_1}}(k)
	\overline{\widehat{u_{N_2}}(k-q)}.
\end{equation*}
For a nonzero summand after applying $P_{N_3}$,
\begin{equation*}
	|k|\sim N_1,\qquad
	|k-q|\sim N_2,\qquad
	|q|\sim N_3.
\end{equation*}
Consequently, $N_3\lesssim\max(N_1,N_2)$.
By symmetry, we assume $N_1\geq N_2$. By Plancherel, H\"older's inequality and \eqref{eq:spatialu},
	\begin{equation*}
		\left\|
		\mathcal F_x^{-1}
		\left|
		\widehat{P_{N_3}D
			(u_{N_1}\overline{u_{N_2}})}
		\right|
		\right\|_{L_x^2}
		\lesssim N_3
		\left\|
		u_{N_1}
		\overline{u_{N_2}}
		\right\|_{L_x^2}
		\lesssim N_3
		\|u_{N_1}\|_{L_x^2}
		\left\|
		u_{N_2}
		\right\|_{L_x^\infty}.
	\end{equation*}
	\begin{equation*}
		\lesssim
		N_3N_1^{-1}m(N_1)^{-1}m(N_2)^{-1}
		A_{N_1}A_{N_2}\lesssim N_3N_1^{-1}m(N_1)^{-2}A^2.
	\end{equation*}
  Summing all dyadic indices, we obtain
	\begin{equation*}
		\left\|
		\mathcal F_x^{-1}
		\left|\widehat{P_{N_3}D(|u|^2)}\right|
		\right\|_{L_t^\infty L_x^2}
		\lesssim
		N_3A^2
		\sum_{N_1\gtrsim N_3}
		N_1^{-1}m(N_1)^{-2}
		\bigl(1+\log(2+N_1)\bigr)
	\end{equation*}
	\begin{equation*}
		\lesssim_\eps
		N_3A^2
		\sum_{N_1\gtrsim N_3}
		N_1^{-1+\eps/8}m(N_1)^{-2}.
	\end{equation*}
	
	If $N_3\leq N$, then
	\begin{equation*}
		N_3\sum_{\substack{N_1\gtrsim N_3\\N_1\leq2N}}
		N_1^{-1+\eps/8}m(N_1)^{-2}
		\lesssim_\eps
		N_3\sum_{N_1\gtrsim N_3}N_1^{-1+\eps/8}
		\lesssim_\eps N_3^{\eps/8}.
	\end{equation*}
	and 
	\begin{equation*}
		N_3\sum_{N_1>2N}
		N_1^{-1+\eps/8}m(N_1)^{-2}
		=
		N_3N^{2s-2}
		\sum_{N_1>2N}N_1^{1-2s+\eps/8}
		\lesssim_\eps N_3N^{-1+\eps/8}
		\leq N_3^{\eps/8}.
	\end{equation*}
	This proves \eqref{eq:newpair} when $N_3\leq N$.
	
	If $N_3>N$, then 
	\begin{equation*}
		N_3\sum_{N_1\gtrsim N_3}
		N_1^{-1+\eps/8}m(N_1)^{-2}
		\lesssim
		N_3N^{2s-2}
		\sum_{N_1\gtrsim N_3}N_1^{1-2s+\eps/8}
		\lesssim_\eps
		N^{2s-2}N_3^{2-2s+\eps/8}
		\lesssim
		N_3^{\eps/8}m(N_3)^{-2}.
	\end{equation*}
	This completes the proof of \eqref{eq:newpair}.
	
	\medskip
	\noindent\textbf{The estimate for $\mathcal F_2$.}
	
	We first estimate the summand in \eqref{eq:Fj-explicit}
	with $j=2$ and $a=1$. Apart from the fixed constant and sign, its integral is
	\begin{equation*}
		\int_J\sum_{p,r,l}
		b_{2,\pm_1,\pm_2}(p,r,l)
		\overline{\widehat u(t,p)}\widehat u(t,l)
		\widehat{D(|u|^2)}(t,p-r)
		\widehat{n_{\pm_2}}(t,r-l)\,dt.
	\end{equation*}
	Decompose the two displayed occurrences of $u$,
	$D(|u|^2)$, and $n_{\pm_2}$ so that
	\begin{equation*}
		|p|\sim N_1,\qquad
		|l|\sim N_2,\qquad
		|p-r|\sim N_3,\qquad
		|r-l|\sim N_4.
	\end{equation*}
	We first restrict the sum to
	\begin{equation*}
		\max(|p|,|r|,|l|)\sim N_5,
	\end{equation*}
	and then sum over dyadic $N_5$.
	
	By Lemma~\ref{lem:b2}, every nonzero summand satisfies
	\begin{equation*}
		N_5\geq N/2,
		\qquad
		|b_{2,\pm_1,\pm_2}(p,r,l)|\lesssim N_5^{-2},
		\qquad
		\min(|p|,|l|)\leq N/4.
	\end{equation*}
	Therefore
	\begin{equation*}
		N_1,N_2,N_3,N_4\lesssim N_5,
		\qquad
		m(N_1)^{-1}m(N_2)^{-1}\lesssim m(N_5)^{-1}.
	\end{equation*}
	
     Applying (\ref{eq:newpair}), one can obtain
	\begin{equation*}
		\left|
		\int_J
		\sum_{\substack{p,r,l\\
				N_5\leq\max(|p|,|r|,|l|)<2N_5}}
		b_{2,\pm_1,\pm_2}(p,r,l)
		\overline{\widehat{u_{N_1}}(t,p)}
		\widehat{u_{N_2}}(t,l)
		\widehat{P_{N_3}D(|u|^2)}(t,p-r)
		\widehat{n_{\pm_2,N_4}}(t,r-l)\,dt
		\right|
	\end{equation*}
	\begin{equation*}
		\lesssim N_5^{-2}\int_J
		\left\|\mathcal F_x^{-1}
		|\widehat{u_{N_1}}(t,\cdot)|\right\|_{L_x^\infty}
		\left\|\mathcal F_x^{-1}
		|\widehat{u_{N_2}}(t,\cdot)|\right\|_{L_x^\infty}
		\left\|\mathcal F_x^{-1}
		\left|\widehat{P_{N_3}D(|u|^2)}(t,\cdot)\right|
		\right\|_{L_x^2}
		\|n_{\pm_2,N_4}(t)\|_{L_x^2}\,dt
	\end{equation*}
	\begin{equation*}
		\lesssim_\eps
		hN_5^{-2}
		m(N_1)^{-1}m(N_2)^{-1}
		N_3^{\eps/8}m(N_3)^{-2}
		A_{N_1}A_{N_2}B_{N_4}A^2
	\end{equation*}
	\begin{equation*}
		\lesssim_\eps
		hN_5^{-2+\eps/8}m(N_5)^{-3}A^5.
	\end{equation*}

	Similarly, we can obtain the same bound for the summand with $a=2$. Summing both summands and all signs yields
	\begin{equation*}
		\bigg|\int_J\mathcal F_2(t)\,dt\bigg|
		\lesssim_\eps
		hA^5\sum_{N_5\geq N/2}
		N_5^{-2+\eps/4}m(N_5)^{-3}\lesssim hN^{-2+\eps}A^5,
	\end{equation*}
	which proves \eqref{eq:F2}.
	
	\medskip
	\noindent\textbf{The coefficient estimates used for $\mathcal F_3$.}
	
	We use the bounds obtained in the proof of
	Lemma~\ref{lem:b3geometry}, before they were weakened
	to \eqref{eq:b3bound}, and record the corresponding
	frequency restrictions.
	
	For $j=3$ in \eqref{eq:Fj-explicit}, write
	\begin{equation*}
		p_0=p,\qquad p_1=r,\qquad p_2=z,\qquad p_3=l.
	\end{equation*}
    Apply the dyadic decomposition so that
	\begin{equation*}
		|p|\sim N_1,\qquad |l|\sim N_2,\qquad
		|p-r|\sim N_3,\qquad |r-z|\sim N_4,\qquad
		|z-l|\sim N_5.
	\end{equation*}
	We now assign two additional dyadic parameters $N_6,N_7$
	according to the defining region of $b_3$.
	
	In \eqref{eq:b3-single-r}, take
	$N_6=N_7\sim|r|$; in \eqref{eq:b3-single-z}, take
	$N_6=N_7\sim|z|$.
	The denominator estimates give, respectively,
	\begin{equation*}
		|b_{3,\pm_1,\pm_2,\pm_3}(p,r,z,l)|
		\lesssim |r|^{-4}
		\quad\text{and}\quad
		|b_{3,\pm_1,\pm_2,\pm_3}(p,r,z,l)|
		\lesssim |z|^{-4}.
	\end{equation*}
	Thus, in both cases,
	\begin{equation*}
		|b_{3,\pm_1,\pm_2,\pm_3}(p,r,z,l)|
		\lesssim N_6^{-2}N_7^{-2},
		\qquad
		|p|,|l|\lesssim N_6,
	\end{equation*}
	and
	\begin{equation*}
		N_3,N_4,N_5\lesssim N_7.
	\end{equation*}
	
	In \eqref{eq:b3separated}, take
	\begin{equation*}
		N_6\sim\min(|r|,|z|),
		\qquad
		N_7\sim\max(|r|,|z|).
	\end{equation*}
	The defining formula and the denominator estimates give
	\begin{equation*}
		|b_{3,\pm_1,\pm_2,\pm_3}(p,r,z,l)|
		\lesssim
		\min(|r|,|z|)^{-2}\max(|r|,|z|)^{-2}
		\lesssim N_6^{-2}N_7^{-2}.
	\end{equation*}
	Moreover,
	\begin{equation*}
		|p|,|l|\leq N/16\lesssim N_6,
		\qquad
		N_3,N_4,N_5\lesssim N_7.
	\end{equation*}
	
	In \eqref{eq:b3comparable}, take
	$N_6=N_7\sim\max(|r|,|z|)$.
	The defining formula gives
	\begin{equation*}
		|b_{3,\pm_1,\pm_2,\pm_3}(p,r,z,l)|
		\lesssim N_6^{-4}=N_6^{-2}N_7^{-2},
	\end{equation*}
	and again
	\begin{equation*}
		|p|,|l|\lesssim N_6,
		\qquad
		N_3,N_4,N_5\lesssim N_7.
	\end{equation*}
	
	Finally, in \eqref{eq:b3highlow}, estimate its two
	fractions separately.
	For the fraction containing $b_{2,\pm_1,\pm_2}(p,r,z)$,
	take
	\begin{equation*}
		N_6\sim\max(|p|,|l|),
		\qquad
		N_7\sim\max(|p|,|r|,|z|).
	\end{equation*}
	By \eqref{eq:b2bound} and \eqref{eq:b3-exterior-phase},
	\begin{equation*}
		\left|
		\frac{b_{2,\pm_1,\pm_2}(p,r,z)}
		{\Omega_{3,\pm_1,\pm_2,\pm_3}(p,r,z,l)}
		\right|
		\lesssim N_6^{-2}N_7^{-2}.
	\end{equation*}
	The corresponding spatial restrictions are
	\begin{equation*}
		|p|,|l|\lesssim N_6,\qquad
		N_3,N_4\lesssim N_7,\qquad
		N_5\lesssim\max(N_6,N_7).
	\end{equation*}
	
	For the fraction containing $b_{2,\pm_2,\pm_3}(r,z,l)$,
	take
	\begin{equation*}
		N_6\sim\max(|p|,|l|),
		\qquad
		N_7\sim\max(|r|,|z|,|l|).
	\end{equation*}
	Then
	\begin{equation*}
		\left|
		\frac{b_{2,\pm_2,\pm_3}(r,z,l)}
		{\Omega_{3,\pm_1,\pm_2,\pm_3}(p,r,z,l)}
		\right|
		\lesssim N_6^{-2}N_7^{-2},
	\end{equation*}
	and
	\begin{equation*}
		|p|,|l|\lesssim N_6,\qquad
		N_4,N_5\lesssim N_7,\qquad
		N_3\lesssim\max(N_6,N_7).
	\end{equation*}
	
	In every case above, a nonzero coefficient satisfies
	\begin{equation*}
		N_6,N_7\gtrsim N,
		\qquad
		\min(|p|,|l|)\leq N/16, \qquad \max(|p|,|l|)\lesssim N_6.
	\end{equation*}
	Therefore
	\begin{equation}\label{eq:F3-u-weights}
		m(N_1)^{-1}m(N_2)^{-1}
		\lesssim m(N_6)^{-1}.
	\end{equation}
	Also, at least two of $N_3,N_4,N_5$ are bounded by
	$N_7$. 
	
	\begin{remark}
		The bound from Lemma~\ref{lem:b3geometry}
		\begin{equation*}
			|b_{3,\pm_1,\pm_2,\pm_3}(p,r,z,l)|
			\lesssim
			N^{-2}\max(|p|,|r|,|z|,|l|)^{-2}
		\end{equation*}
		is sufficient for the fixed-time estimate of $\mathcal C_5$.
		For $\mathcal F_3$, however, $D(|u|^2)$ introduces two
		additional occurrences of $u$.
		Using only this bound and a single dyadic scale
		\begin{equation*}
			N_0\sim\max(|p|,|r|,|z|,|l|)
		\end{equation*}
		would lead to
		\begin{equation*}
			hA^6N^{-2}
			\sum_{N_0\gtrsim N}
			N_0^{-1+\eps}m(N_0)^{-3}.
		\end{equation*}
		For $N_0\geq2N$,
		\begin{equation*}
			N_0^{-1+\eps}m(N_0)^{-3}
			=N^{3s-3}N_0^{2-3s+\eps},
		\end{equation*}
		so this argument requires $s>2/3$.
		To obtain the estimate for $s>1/2$, we retain
		the stronger bound $N_6^{-2}N_7^{-2}$ in each defining
		region, together with the fact that at least two of
		$N_3,N_4,N_5$ are at most $CN_7$.
	\end{remark}
	
	\medskip
	\noindent\textbf{The estimate for $\mathcal F_3$.}
	
	We first treat the case $a=1$ in \eqref{eq:Fj-explicit}.
	For each of the defining regions considered above,
	the coefficient is bounded by $N_6^{-2}N_7^{-2}$.
	In the region of \eqref{eq:b3highlow}, the following
	calculation is applied separately to its two fractions.
	
	Suppose first that $N_4\leq N_5$. We apply \eqref{eq:spatialu}, \eqref{eq:spatialwave},
	\eqref{eq:newpair} and \eqref{eq:F3-u-weights},
	\begin{equation*}
		\left|
		\int_J\sum_{p,r,z,l}
		b_{3,\pm_1,\pm_2,\pm_3}(p,r,z,l)
		\overline{\widehat{u_{N_1}}(t,p)}
		\widehat{u_{N_2}}(t,l)
		\widehat{P_{N_3}D(|u|^2)}(t,p-r)
		\widehat{n_{\pm_2,N_4}}(t,r-z)
		\widehat{n_{\pm_3,N_5}}(t,z-l)\,dt
		\right|
	\end{equation*}
	\begin{equation*}
		\lesssim N_6^{-2}N_7^{-2}\int_J
		\left\|\mathcal F_x^{-1}
		|\widehat{u_{N_1}}(t,\cdot)|\right\|_{L_x^\infty}
		\left\|\mathcal F_x^{-1}
		|\widehat{u_{N_2}}(t,\cdot)|\right\|_{L_x^\infty}
	\end{equation*}
	\begin{equation*}
		\qquad{}\times
		\left\|\mathcal F_x^{-1}
		\left|\widehat{P_{N_3}D(|u|^2)}(t,\cdot)\right|
		\right\|_{L_x^2}
		\left\|\mathcal F_x^{-1}
		|\widehat{n_{\pm_2,N_4}}(t,\cdot)|\right\|_{L_x^\infty}
		\|n_{\pm_3,N_5}(t)\|_{L_x^2}\,dt.
	\end{equation*}
	\begin{equation*}
		\lesssim_\eps
		hN_6^{-2}N_7^{-2}
		\bigl(m(N_1)^{-1}A_{N_1}\bigr)
		\bigl(m(N_2)^{-1}A_{N_2}\bigr)
		\bigl(N_3^{\eps/8}m(N_3)^{-2}A^2\bigr)
		\bigl(N_4B_{N_4}\bigr)B_{N_5}
	\end{equation*}
	\begin{equation*}
		\lesssim_\eps
		hN_6^{-2}N_7^{-2}N_4
		m(N_6)^{-1}
		N_3^{\eps/8}m(N_3)^{-2}A^6.
	\end{equation*}
	If $N_5<N_4$, interchange the two spatial exponents
	assigned to $n_{\pm_2,N_4}$ and $n_{\pm_3,N_5}$.
	Thus $N_4$ in the last bound can be replaced by
	$\min(N_4,N_5)$.
	
	The spatial restrictions proved above imply
	\begin{equation*}
		\min(N_4,N_5)\lesssim N_7,
		\qquad
		N_3\lesssim\max(N_6,N_7),
	\end{equation*}
	and therefore
	\begin{equation*}
		N_3^{\eps/8}m(N_3)^{-2}
		\lesssim
		\max(N_6,N_7)^{\eps/8}
		m(\max(N_6,N_7))^{-2}.
	\end{equation*}
	The same argument applies to the cases $a=2,3$.
	
	For fixed $N_6,N_7$, all five spatial indices
	$N_1,\ldots,N_5$ are at most $C\max(N_6,N_7)$.
	Hence their number of possible dyadic choices is bounded by $(N_6N_7)^{\eps/8}$.
	
	Thus, summing all defining regions, signs,
	and dyadic indices gives
	\begin{equation}\label{eq:F3sum}
		\bigg|\int_J\mathcal F_3(t)\,dt\bigg|
		\lesssim_\eps
		hA^6
		\sum_{N_6,N_7\gtrsim N}
		N_6^{-2+\eps/4}N_7^{-1+\eps/4}
		m(N_6)^{-1}m(\max(N_6,N_7))^{-2}.
	\end{equation}
	
	We first sum over $N_7\geq N_6$.
	Since $N_6,N_7\gtrsim N$,
	\begin{equation*}
		N_6^{-2+\eps/4}N_7^{-1+\eps/4}
		m(N_6)^{-1}m(N_7)^{-2}
		\lesssim
		N^{3s-3}N_6^{-1-s+\eps/4}N_7^{1-2s+\eps/4}.
	\end{equation*}
	Our choice of $\eps$ gives $1-2s+\eps/4<0$, so
	\begin{equation*}
		\sum_{N_6\gtrsim N}\sum_{N_7\geq N_6}
		N_6^{-2+\eps/4}N_7^{-1+\eps/4}
		m(N_6)^{-1}m(N_7)^{-2}
	\end{equation*}
	\begin{equation*}
		\lesssim_\eps
		N^{3s-3}
		\sum_{N_6\gtrsim N}
		N_6^{-1-s+\eps/4}N_6^{1-2s+\eps/4}
		=
		N^{3s-3}\sum_{N_6\gtrsim N}N_6^{-3s+\eps/2}
		\lesssim_\eps N^{-3+\eps/2}.
	\end{equation*}
	
	For $N_6>N_7$, we have
	\begin{equation*}
		N_6^{-2+\eps/4}N_7^{-1+\eps/4}m(N_6)^{-3}
		\lesssim
		N^{3s-3}N_6^{1-3s+\eps/4}N_7^{-1+\eps/4}.
	\end{equation*}
	Since $1-3s+\eps/4<0$,
	\begin{equation*}
		\sum_{N_7\gtrsim N}\sum_{N_6>N_7}
		N_6^{-2+\eps/4}N_7^{-1+\eps/4}m(N_6)^{-3}
	\end{equation*}
	\begin{equation*}
		\lesssim_\eps
		N^{3s-3}
		\sum_{N_7\gtrsim N}
		N_7^{-1+\eps/4}N_7^{1-3s+\eps/4}
		=
		N^{3s-3}\sum_{N_7\gtrsim N}N_7^{-3s+\eps/2}
		\lesssim_\eps N^{-3+\eps/2}.
	\end{equation*}
	Substituting these two bounds into \eqref{eq:F3sum},
	we conclude that
	\begin{equation*}
		\bigg|\int_J\mathcal F_3(t)\,dt\bigg|
		\lesssim_\eps hN^{-3+\eps/2}A^6
		\lesssim_\eps hN^{-3+\eps}A^6.
	\end{equation*}
	This proves \eqref{eq:F3}.
	
\end{proof}

\vspace{10pt}
\section{The estimate for $\mathcal R_4$}

Since $b_4=0$, definition \eqref{eq:Rj} gives
\begin{equation*}
	R_{4,\pm_1,\pm_2,\pm_3,\pm_4}(p_0,p_1,p_2,p_3,p_4)
	=
	b_{3,\pm_2,\pm_3,\pm_4}(p_1,p_2,p_3,p_4)
	-b_{3,\pm_1,\pm_2,\pm_3}(p_0,p_1,p_2,p_3).
\end{equation*}
 We estimate the above two terms separately.

\begin{lemma}\label{lem:six}
	Let $1/2<s<1$. For smooth solutions on $J$, with $h=|J|$, we have
	\begin{equation}\label{eq:six}
		\abs{\int_J\mathcal R_4(t)\,dt}
		\lesssim_\eps
		N^{-5/2+\eps}h^{1/2-\eps}A^6
	\end{equation}
	for every sufficiently small $\eps>0$.
\end{lemma}

\begin{proof}
	
	After relabelling the frequencies and fixing the signs  $\pm_0,\pm_1,\pm_2,\pm_3$, we estimate the first term in
	$\mathcal R_4$:
	\begin{equation*}
		\sum_{p,r,z,l,k}
		b_{3,\pm_1,\pm_2,\pm_3}(p,r,z,l)
		\overline{\widehat u(k)}\widehat u(l)
		\widehat{n_{\pm_0}}(k-p)
		\widehat{n_{\pm_1}}(p-r)
		\widehat{n_{\pm_2}}(r-z)
		\widehat{n_{\pm_3}}(z-l).
	\end{equation*}
	Exchanging the order of the five frequencies and reversing
	all signs reduces the second term to the complex conjugate
	of a sum of the same form, by \eqref{eq:b3-symmetry}.
	It therefore suffices to estimate the displayed sum.
	
	Apply the spatial dyadic decomposition so that
	\begin{equation*}
         |l|\sim N_2,\qquad
		|p-r|\sim N_3,\qquad |r-z|\sim N_4,\qquad
		|z-l|\sim N_5,
	\end{equation*}
	and
	\begin{equation*}
		|k|\sim N_8,\qquad |k-p|\sim N_9.
	\end{equation*}
	We also divide the sum according to $|p|\sim N_1$.
	This restriction is imposed on the summation variable
	$p$; there is no occurrence of $u_{N_1}$ in the sum.
	
	Use the frequency restrictions indexed by $N_6,N_7$
	from the preceding proof. On each such region,
	\begin{equation*}
		N_6,N_7\gtrsim N,\qquad
		|b_{3,\pm_1,\pm_2,\pm_3}(p,r,z,l)|
		\lesssim N_6^{-2}N_7^{-2},
	\end{equation*}
	\begin{equation*}
		|p|,|l|\lesssim N_6,\qquad
		N_3,N_4,N_5\lesssim\max(N_6,N_7),
	\end{equation*}
	and at least two of $N_3,N_4,N_5$ are at most $CN_7$.
	Consequently,
	\begin{equation}\label{eq:six-old-wave-product}
		N_3N_4N_5
		\lesssim N_7^2\max(N_6,N_7).
	\end{equation}
	In the region of \eqref{eq:b3highlow}, all calculations
	below are applied separately to the two fractions defining
	$b_3$, with the corresponding restrictions on $N_6,N_7$.
	
	Fix a sufficiently large absolute constant $C_0$.
	We consider separately
	\begin{equation*}
		N_8\leq C_0\max(N_6,N_7)
		\qquad\text{and}\qquad
		N_8>C_0\max(N_6,N_7).
	\end{equation*}
	
	\medskip
	\noindent\textbf{ The region
		$N_8\leq C_0\max(N_6,N_7)$.}
	
	Since $|p|\lesssim N_6$,
	\begin{equation*}
		|k-p|\leq |k|+|p|
		\lesssim \max(N_6,N_7),
	\end{equation*}
	and hence
	\begin{equation*}
		N_9\lesssim\max(N_6,N_7).
	\end{equation*}
	Combining this with \eqref{eq:six-old-wave-product} gives
	\begin{equation}\label{eq:six-four-wave-product}
		N_3N_4N_5N_9
		\lesssim N_7^2\max(N_6,N_7)^2.
	\end{equation}

	Recall that Bernstein's inequality
	and \eqref{eq:smoothcutoff} give
	\begin{equation*}
		\left\|
		\mathcal F_{t,x}^{-1}
		|\widetilde{\psi_hu_{N_a}}|
		\right\|_{L^\infty_{t,x}}
		\lesssim
		m(N_a)^{-1}
		\|\psi_h Iu_{N_a}\|_{X_S^{1,1/2,1}}
		\lesssim m(N_a)^{-1}A_{N_a}.
	\end{equation*}
	The same estimate holds for $\overline{u_{N_a}}$.
	Also, by \eqref{eq:wavefour} and Remark~\ref{absolute},
	\begin{equation*}
		\left\|
		\mathcal F_{t,x}^{-1}
		|\widetilde{1_Jn_{\pm,N_a}}|
		\right\|_{L^4_{t,x}}
		\lesssim N_a^{3/8}h^{1/8}B_{N_a}.
	\end{equation*}
	
	Insert $\psi_h$ on both occurrences of $u$ and $1_J$
	on all four occurrences of $n$. For fixed signs and dyadic indices, H\"older's inequality
	with exponents $\infty,\infty,4,4,4,4$ therefore gives
	\begin{equation*}
		\bigg|
		\int_J\sum_{p,r,z,l,k}
		b_{3,\pm_1,\pm_2,\pm_3}(p,r,z,l)
		\overline{\widehat{u_{N_8}}(t,k)}
		\widehat{u_{N_2}}(t,l)
		\widehat{n_{\pm_0,N_9}}(t,k-p)
		\widehat{n_{\pm_1,N_3}}(t,p-r)
		\widehat{n_{\pm_2,N_4}}(t,r-z)
		\widehat{n_{\pm_3,N_5}}(t,z-l)\,dt
		\bigg|
	\end{equation*}
	\begin{equation*}
		\lesssim
		N_6^{-2}N_7^{-2}
		\bigl(m(N_8)^{-1}A_{N_8}\bigr)
		\bigl(m(N_2)^{-1}A_{N_2}\bigr)
		\prod_{a\in\{3,4,5,9\}}
		\bigl(N_a^{3/8}h^{1/8}B_{N_a}\bigr)
	\end{equation*}
	\begin{equation*}
		\lesssim
		h^{1/2}N_6^{-2}N_7^{-5/4}
		\max(N_6,N_7)^{3/4}
		m(N_6)^{-1}m(\max(N_6,N_7))^{-1}A^6.
	\end{equation*}
	In the last step we used \eqref{eq:six-four-wave-product},
	$N_2\lesssim N_6$, and
	$N_8\lesssim\max(N_6,N_7)$.
	
	For fixed $N_6,N_7$, all seven remaining dyadic indices
	are bounded by $C\max(N_6,N_7)$. Their number of choices
	is at most $(N_6N_7)^{\eps/16}$.
	
	Thus it suffices to estimate
	\begin{equation}\label{eq:sixsum}
		\sum_{N_6,N_7\gtrsim N}
		N_6^{-2+\eps/16}N_7^{-5/4+\eps/16}
		\max(N_6,N_7)^{3/4}
		m(N_6)^{-1}m(\max(N_6,N_7))^{-1}.
	\end{equation}
	
	Choose $0<\eps<\min(1,s-1/2)$.
	For $N_7\geq N_6$, the corresponding part of
	\eqref{eq:sixsum} satisfies
	\begin{equation*}
		\sum_{N_6\gtrsim N}\sum_{N_7\geq N_6}
		N_6^{-2+\eps/16}N_7^{-1/2+\eps/16}
		m(N_6)^{-1}m(N_7)^{-1}
	\end{equation*}
	\begin{equation*}
		\lesssim
		N^{2s-2}
		\sum_{N_6\gtrsim N}N_6^{-1-s+\eps/16}
		\sum_{N_7\geq N_6}N_7^{1/2-s+\eps/16}
		\lesssim_\eps
		N^{2s-2}
		\sum_{N_6\gtrsim N}N_6^{-1/2-2s+\eps/8}
		\lesssim_\eps N^{-5/2+\eps/8}.
	\end{equation*}
	For $N_6>N_7$, we instead have
	\begin{equation*}
		\sum_{N_7\gtrsim N}\sum_{N_6>N_7}
		N_6^{-5/4+\eps/16}N_7^{-5/4+\eps/16}
		m(N_6)^{-2}
	\end{equation*}
	\begin{equation*}
		\lesssim
		N^{2s-2}
		\sum_{N_7\gtrsim N}N_7^{-5/4+\eps/16}
		\sum_{N_6>N_7}N_6^{3/4-2s+\eps/16}
		\lesssim_\eps
		N^{2s-2}
		\sum_{N_7\gtrsim N}N_7^{-1/2-2s+\eps/8}
		\lesssim_\eps N^{-5/2+\eps/8}.
	\end{equation*}
	Therefore the sum of the absolute values of the
	integrals in this region is bounded by
	\begin{equation*}
		 h^{1/2}N^{-5/2+\eps/8}A^6.
	\end{equation*}
	
	\medskip
	\noindent\textbf{ The region
		$N_8>C_0\max(N_6,N_7)$.}
	
	For $C_0$ sufficiently large, $|p|\leq |k|/2$.
	Hence $N_9\sim N_8$. First consider the ordering $N_3\geq N_4\geq N_5$.
	Since at least two of $N_3,N_4,N_5$ are at most $CN_7$, we see $N_4,N_5\lesssim N_7$.
	
	Insert $1_J$ on $n_{\pm_0,N_9}$ and $n_{\pm_1,N_3}$,
	and $\psi_h$ on the other four functions. 
	By \eqref{eq:pair} with $\eps/16$ in place of $\eps$,
	\begin{equation*}
		\left\|
		\left(\mathcal F_{t,x}^{-1}
		|\widetilde{\psi_h\overline{u_{N_8}}}|\right)
		\left(\mathcal F_{t,x}^{-1}
		|\widetilde{1_Jn_{\pm_0,N_9}}|\right)
		\right\|_{L^2_{t,x}}
	\end{equation*}
	\begin{equation*}
		\lesssim_\eps
		m(N_8)^{-1}
		\left(
		N_8^{-1+\eps/16}
		+h^{1/2-\eps/16}N_8^{-1/2+\eps/16}
		\right)A_{N_8}B_{N_9}.
	\end{equation*}
	Plancherel's identity and \eqref{eq:spatialwave} give
	\begin{equation*}
		\left\|
		\mathcal F_{t,x}^{-1}
		|\widetilde{1_Jn_{\pm_1,N_3}}|
		\right\|_{L^2_{t,x}}
		=
		\|1_Jn_{\pm_1,N_3}\|_{L^2_{t,x}}
		\lesssim h^{1/2}B_{N_3}.
	\end{equation*}
	Similarly, Bernstein's inequality and
	\eqref{eq:smoothcutoff} give
	\begin{equation*}
		\left\|
		\mathcal F_{t,x}^{-1}
		|\widetilde{\psi_hn_{\pm,N_a}}|
		\right\|_{L^\infty_{t,x}}
		\lesssim
		N_a\|\psi_hn_{\pm,N_a}\|_{X_{W_\pm}^{0,1/2,1}}
		\lesssim N_aB_{N_a},
		\qquad a=4,5.
	\end{equation*}
	and 
	\begin{equation*}
		\left\|
		\mathcal F_{t,x}^{-1}
		|\widetilde{\psi_hu_{N_2}}|
		\right\|_{L^\infty_{t,x}}
		\lesssim m(N_2)^{-1}A_{N_2}.
	\end{equation*}
	
	Applying Cauchy--Schwarz to the two $L^2_{t,x}$ quantities
	and using $L^\infty_{t,x}$ for the other three functions,
	we obtain
	\begin{equation*}
		\bigg|
		\int_J\sum_{p,r,z,l,k}
		b_{3,\pm_1,\pm_2,\pm_3}(p,r,z,l)
		\overline{\widehat{u_{N_8}}(t,k)}
		\widehat{u_{N_2}}(t,l)
		\widehat{n_{\pm_0,N_9}}(t,k-p)
		\widehat{n_{\pm_1,N_3}}(t,p-r)
		\widehat{n_{\pm_2,N_4}}(t,r-z)
		\widehat{n_{\pm_3,N_5}}(t,z-l)\,dt
		\bigg|
	\end{equation*}
	\begin{equation*}
		\lesssim_\eps
		N_6^{-2}N_7^{-2}
		\bigl(m(N_2)^{-1}A_{N_2}\bigr)
		\bigl(h^{1/2}B_{N_3}\bigr)
		\bigl(N_4B_{N_4}\bigr)
		\bigl(N_5B_{N_5}\bigr)
	\end{equation*}
	\begin{equation*}
		\qquad{}\times
		m(N_8)^{-1}
		\left(
		N_8^{-1+\eps/16}
		+h^{1/2-\eps/16}N_8^{-1/2+\eps/16}
		\right)A_{N_8}B_{N_9}
	\end{equation*}
	\begin{equation*}
		\lesssim_\eps
		N_6^{-2}m(N_6)^{-1}m(N_8)^{-1}
		\left(
		h^{1/2}N_8^{-1+\eps/16}
		+h^{1-\eps/16}N_8^{-1/2+\eps/16}
		\right)A^6
	\end{equation*}
	\begin{equation}\label{eq:sixlargeJ}
		\lesssim_\eps
		N^{2s-2}N_6^{-1-s}
		\left(
		h^{1/2}N_8^{-s+\eps/16}
		+h^{1-\eps/16}N_8^{1/2-s+\eps/16}
		\right)A^6.
	\end{equation}
	The other orderings of $N_3,N_4,N_5$ give the same
	bound: put the largest of the three in $L^2_{t,x}$
	and the other two in $L^\infty_{t,x}$.
	
	For fixed $N_6,N_8$, the restrictions imply
	\begin{equation*}
		N\lesssim N_7\lesssim N_8,\qquad
		N_1,N_2,N_3,N_4,N_5\lesssim N_8,\qquad
		N_9\sim N_8.
	\end{equation*}
	Thus the number of remaining dyadic choices is at most $N_8^{\eps/16}$.
	
	After summing these choices, the first term in
	\eqref{eq:sixlargeJ} is bounded by
	\begin{equation*}
		 h^{1/2}A^6N^{2s-2}
		\sum_{N_6\gtrsim N}N_6^{-1-s}
		\sum_{N_8\gtrsim N_6}N_8^{-s+\eps/8}
	\end{equation*}
	\begin{equation*}
		\lesssim_\eps
		h^{1/2}A^6N^{2s-2}
		\sum_{N_6\gtrsim N}N_6^{-1-2s+\eps/8}
		\lesssim_\eps h^{1/2}N^{-3+\eps/8}A^6.
	\end{equation*}
	The second term is bounded by
	\begin{equation*}
		 h^{1-\eps/16}A^6N^{2s-2}
		\sum_{N_6\gtrsim N}N_6^{-1-s}
		\sum_{N_8\gtrsim N_6}N_8^{1/2-s+\eps/8}
	\end{equation*}
	\begin{equation*}
		\lesssim_\eps
		h^{1-\eps/16}A^6N^{2s-2}
		\sum_{N_6\gtrsim N}N_6^{-1/2-2s+\eps/8}
		\lesssim_\eps
		h^{1-\eps/16}N^{-5/2+\eps/8}A^6.
	\end{equation*}
	
	Combining the two frequency regions and the two terms
	in $R_4$, we obtain
	\begin{equation*}
		\bigg|\int_J\mathcal R_4(t)\,dt\bigg|
		\lesssim_\eps
		\left(
		h^{1/2}N^{-5/2+\eps/8}
		+h^{1/2}N^{-3+\eps/8}
		+h^{1-\eps/16}N^{-5/2+\eps/8}
		\right)A^6
	\end{equation*}
	\begin{equation*}
		\lesssim_\eps
		h^{1/2-\eps}N^{-5/2+\eps}A^6,
	\end{equation*}
	since $N\geq1$ and $0<h\leq1$.
	This proves \eqref{eq:six}.
	
\end{proof}

\vspace{10pt}
\section{Completion of the almost conservation estimate}

\begin{proof}[Proof of Proposition~\ref{prop:increment}]
	Integrating \eqref{eq:energyidentity}, we obtain
	\begin{equation*}
		|\mathcal E_N(t_1)-\mathcal E_N(t_0)|
		\leq
		\sum_{j=1}^4\left|\int_J\mathcal R_j(t)\,dt\right|
		+\sum_{j=1}^3\left|\int_J\mathcal F_j(t)\,dt\right|.
	\end{equation*}
	Apply \eqref{eq:oldcubic}, \eqref{eq:quarticremainder},
	\eqref{eq:quinticremainder}, and \eqref{eq:six}
	to $\mathcal R_1,\mathcal R_2,\mathcal R_3,\mathcal R_4$,
	respectively, and \eqref{eq:oldquartic}, \eqref{eq:F2},
	and \eqref{eq:F3} to $\mathcal F_1,\mathcal F_2,\mathcal F_3$,
	respectively. Since $N\geq1$ and $0<h\leq1$,
	definitions \eqref{eq:d3}--\eqref{eq:d6} give
	\eqref{eq:increment}.
\end{proof}

\vspace{10pt}
\section{Proof of Theorem \ref{thm:main}}
\subsection{Modified local well-posedness}
The final ingredient in the proof of Theorem~\ref{thm:main}
is the following local well-posedness result from
\cite[Proposition 4.5]{KishimotoI}.
\begin{lemma}[Modified local theory]\label{lem:local}
	For $1/2\le s<1$, $u_0\in H^s$ and $n_{+0}\in L^2$, set
	$A_0=1+\norm{Iu_0}{H^1}+\norm{n_{+0}}2$. For every sufficiently small
	$\eps>0$ there is a time interval of length
	$h\sim_\eps A_0^{-2-\eps}$ on which the unique local solution satisfies
	\begin{equation*}
		\norm{Iu}{X_S^{1,1/2,1}([0,h])}
		+\norm{n_+}{X_{W_+}^{0,1/2,1}([0,h])}\lesssim_\eps A_0.
	\end{equation*}
	The constants are independent of $N$. The solution belongs to
	$C_t(H^s\times L^2)$, depends continuously on its data, and is
	unique in $X_S^{s,1/2,1}\times X_{W_+}^{0,1/2,1}$ locally in time.
\end{lemma}

\begin{remark}\label{rem:extension}
	By approximation with smooth solutions and continuous
	dependence in the local well-posedness spaces, the estimates
	for $\mathcal E_N$ in Propositions~\ref{prop:comparison}
	and~\ref{prop:increment} extend to the solutions of
	Lemma~\ref{lem:local}. The same approximation also gives
	conservation of the Schr\"odinger mass.
\end{remark}

\vspace{8pt}
\subsection{Coercivity under the strict mass condition}

We use the torus Gagliardo--Nirenberg inequality in the following
form, recorded in the introduction of \cite{KishimotoI}:
for every $f\in H^1(\T^2)$,
\begin{equation}\label{eq:gn}
	\norm f4^4\le\frac{2}{M(Q)}\norm f2^2\norm{\nabla f}2^2
	+C_\gamma\norm f2^4.
\end{equation}
The ground state and the normalization of $M(Q)$ are those in
\eqref{eq:groundstate}, and the torus norms use the measure fixed
in \eqref{eq:fourier}. We now prove explicitly the coercive
estimate needed below.

\begin{lemma}\label{lem:coercivity}
	Suppose $M(u)=M(u_0)<M(Q)$, and put
	$\theta=M(u_0)/M(Q)$. Then
	\begin{equation}\label{eq:coercivity}
		\mathcal H(Iu,n_+)+C_{\theta,\gamma}M(u_0)^2
		\ge\frac{1-\theta}{1+\theta}\norm{\nabla Iu}2^2
		+\frac{1-\theta}{4}\norm n2^2
		+\frac12\norm{D^{-1}n_t}2^2.
	\end{equation}
\end{lemma}

\begin{proof}
	Since $0<m\le1$, Parseval's identity gives
	$\norm{Iu}2^2\le M(u_0)$. Apply \eqref{eq:gn} to $Iu$:
	\begin{equation*}
		\norm{Iu}4^4\le2\theta\norm{\nabla Iu}2^2
		+C_\gamma M(u_0)^2.
	\end{equation*}
	Set $\delta=(1-\theta)/2$, so $1-\delta=(1+\theta)/2>0$.
	Cauchy--Schwarz and Young's inequality give
	\begin{equation*}
		\abs{\int_{\T^2}n|Iu|^2\,dx}
		\le\frac{1-\delta}{2}\norm n2^2
		+\frac{1}{2(1-\delta)}\norm{Iu}4^4.
	\end{equation*}
	Insert these two inequalities into \eqref{eq:hamiltonian}, using
	$\norm{n_+}2^2=\norm n2^2+\norm{D^{-1}n_t}2^2$.
	The coefficient of $\norm{\nabla Iu}2^2$ is
	$1-\theta/(1-\delta)=(1-\theta)/(1+\theta)$, and that of
	$\norm n2^2$ is $\delta/2=(1-\theta)/4$. Both are positive.
	This proves \eqref{eq:coercivity}.
	
\end{proof}

\vspace{8pt}
\subsection{Iteration of the local theory}

\begin{proof}[Proof of Theorem~\ref{thm:main}]
	Let $(u_0,n_{+0})\in H^s\times L^2$ satisfy the
	assumptions of the theorem. We follow the iteration
	in \cite[Section 4.3]{KishimotoI}, using
	Propositions~\ref{prop:comparison} and~\ref{prop:increment}.
	
	The definition of $m$ gives
	\begin{equation*}
		\langle k\rangle m(k)
		\lesssim N^{1-s}\langle k\rangle^s,
		\qquad 0<m(k)\leq1.
	\end{equation*}
	Consequently,
	\begin{equation*}
		1+\|Iu_0\|_{H^1}+\|n_{+0}\|_2
		\lesssim_{u_0,n_0,n_1}N^{1-s},
	\end{equation*}
	and
	\begin{equation*}
		\|Iu_0\|_2\leq\|u_0\|_2<\|Q\|_{L^2(\mathbb R^2)}.
	\end{equation*}
	By \eqref{eq:gn},
	\begin{equation*}
		\left|\int_{\T^2}n_0|Iu_0|^2\,dx\right|
		\leq\|n_0\|_2\|Iu_0\|_4^2
		\lesssim
		\|n_0\|_2
		\left(
		\|Iu_0\|_2\|\nabla Iu_0\|_2+\|Iu_0\|_2^2
		\right)
		\lesssim_{u_0,n_0,n_1}N^{1-s}.
	\end{equation*}
	Thus, we can obtain
	\begin{equation*}
		1+|\mathcal H(Iu_0,n_{+0})|
		\lesssim_{u_0,n_0,n_1} N^{2(1-s)}.
	\end{equation*}
	
	Mass conservation gives
	\begin{equation*}
		\|Iu(t)\|_2\leq\|u(t)\|_2=\|u_0\|_2.
	\end{equation*}
	Therefore \eqref{eq:coercivity} implies
	\begin{equation*}
		\left(1+\|Iu(t)\|_{H^1}+\|n_+(t)\|_2\right)^2
		\lesssim_{u_0,n_0,n_1}
		\mathcal H(Iu(t),n_+(t))+1.
	\end{equation*}
	In particular, whenever
	\begin{equation}\label{eq:iteration-condition}
		\mathcal H(Iu(t_0),n_+(t_0))
		\lesssim N^{2(1-s)},
	\end{equation}
	we have
	\begin{equation*}
		1+\|Iu(t_0)\|_{H^1}+\|n_+(t_0)\|_2
		\lesssim_{u_0,n_0,n_1}N^{1-s}.
	\end{equation*}
	Lemma~\ref{lem:local} then extends the solution from
	$t_0$ to $t_0+h$, with a uniform choice
	\begin{equation}\label{eq:step}
		h=cN^{-(1-s)(2+\eps)}\leq1,
	\end{equation}
	and gives
	\begin{equation}\label{eq:iteration-local-bound}
		1+\|Iu\|_{X_S^{1,1/2,1}([t_0,t_0+h])}
		+\|n_+\|_{X_{W_+}^{0,1/2,1}([t_0,t_0+h])}
		\lesssim N^{1-s}.
	\end{equation}
	Here and below, constants may depend on the initial data
	and the fixed sufficiently small $\eps>0$, but not on
	$N$ or the number of iterations.
	The embedding $X_\ast^{s,1/2,1}\hookrightarrow C_t H_x^s$ also gives
	\begin{equation*}
		\sup_{t_0\leq t\leq t_0+h}
		\left(1+\|Iu(t)\|_{H^1}+\|n_+(t)\|_2\right)
		\lesssim N^{1-s}.
	\end{equation*}
	
	Suppose that the local theory has been applied $M$ times,
	with \eqref{eq:iteration-condition} satisfied at the
	restart times $0,h,\ldots,(M-1)h$.
	We compare the energy at $Mh$ with its initial
	value:
	\begin{equation*}
		|\mathcal H(Iu(Mh),n_+(Mh))
		-\mathcal H(Iu_0,n_{+0})|
		\leq
		|\mathcal H(Iu(Mh),n_+(Mh))-\mathcal E_N(Mh)|
	\end{equation*}
	\begin{equation*}
		\qquad{}
		+\sum_{j=0}^{M-1}
		|\mathcal E_N((j+1)h)-\mathcal E_N(jh)|
		+|\mathcal E_N(0)-\mathcal H(Iu_0,n_{+0})|.
	\end{equation*}
	By Proposition~\ref{prop:comparison}, the sum of the
	two endpoint errors is at most
	\begin{equation*}
		C\left(
		N^{-1+\eps}N^{3(1-s)}
		+N^{-2+\eps}N^{4(1-s)}
		+N^{-3+\eps}N^{5(1-s)}
		\right)
	\end{equation*}
	\begin{equation*}
		=
		CN^{2(1-s)}
		\left(N^{-s+\eps}+N^{-2s+\eps}+N^{-3s+\eps}\right).
	\end{equation*}
	By Proposition~\ref{prop:increment} and
	\eqref{eq:iteration-local-bound},
	\begin{equation*}
		\sum_{j=0}^{M-1}
		|\mathcal E_N((j+1)h)-\mathcal E_N(jh)|
		\lesssim
		M\Bigl[
		N^{-1+\eps}h^{1/2-\eps}N^{3(1-s)}
		+\bigl(N^{-2+\eps}+N^{-1+\eps}h^{1-\eps}\bigr)
		N^{4(1-s)}
	\end{equation*}
	\begin{equation*}
		+\bigl(N^{-3+\eps}+N^{-2+\eps}h^{1-\eps}\bigr)
		N^{5(1-s)}
		+N^{-5/2+\eps}h^{1/2-\eps}N^{6(1-s)}
		\Bigr].
	\end{equation*}
	
	For $0<\eps<1/6$, formula \eqref{eq:step} gives
	\begin{equation*}
		h^{1/2-\eps}
		=
		c^{1/2-\eps}
		N^{-(1-s)+(1-s)(3\eps/2+\eps^2)}
		\lesssim N^{-(1-s)+2\eps},
	\end{equation*}
	and
	\begin{equation*}
		h^{1-\eps}
		=
		c^{1-\eps}
		N^{-2(1-s)+(1-s)(\eps+\eps^2)}
		\lesssim N^{-2(1-s)+\eps}.
	\end{equation*}
	Substituting these bounds into the preceding estimate
	and dividing by $N^{2(1-s)},$ we obtain
	\begin{equation*}
		\frac{1}{N^{2(1-s)}}
		\sum_{j=0}^{M-1}
		|\mathcal E_N((j+1)h)-\mathcal E_N(jh)|
		\lesssim
		M\Bigl(
		N^{-1+3\eps}+N^{-2s+\eps}+N^{-1+2\eps}
	\end{equation*}
	\begin{equation*}
		\qquad{}
		+N^{-3s+\eps}+N^{-1-s+2\eps}
		+N^{1/2-3s+3\eps}
		\Bigr)
		\lesssim MN^{-1+3\eps}.
	\end{equation*}
	The last inequality follows from $s>1/2$ and $N\geq1$.
	Combining the endpoint and increment estimates yields
	\begin{equation}\label{eq:totalincrement}
		|\mathcal H(Iu(Mh),n_+(Mh))
		-\mathcal H(Iu_0,n_{+0})|
		\lesssim
		\left(N^{-s+\eps}+MN^{-1+3\eps}\right)
		N^{2(1-s)}.
	\end{equation}
	
	Choose $\eps>0$ so small that
	\begin{equation*}
		\eps<\min\left(\frac16,\frac{2s-1}{8}\right).
	\end{equation*}
	First take $N$ sufficiently large that the first term
	on the right of \eqref{eq:totalincrement} is at most
	$\lesssim N^{2(1-s)}$.
	Then choose 
	\begin{equation*}
		M\sim N^{1-3\eps},
	\end{equation*}
	the second term is also at most
	$\lesssim N^{2(1-s)}$. Consequently,
	\begin{equation*}
		\mathcal H(Iu(Mh),n_+(Mh))
		\lesssim N^{2(1-s)}.
	\end{equation*}
	This verifies \eqref{eq:iteration-condition} at the next
	restart time. Induction therefore permits $M\sim N^{1-3\eps}$
	iterations. The same estimates apply at intermediate
	times by shortening the final interval.
	
	For sufficiently large $N$, the resulting existence
	time satisfies
	\begin{equation*}
		Mh\gtrsim
		N^{1-3\eps}N^{-(1-s)(2+\eps)}
		=
		N^{2s-1-(4-s)\eps}.
	\end{equation*}
	Since $2s-1-(4-s)\eps>0$, this time tends to infinity
	with $N$. Given any $T>0$, choose $N$ so that
	\begin{equation*}
		N\sim C(1+T)^{1/(2s-1-(4-s)\eps)},
	\end{equation*}
	 for some $C$ sufficiently large. Then $Mh\geq T$,
	and the solution exists on $[0,T]$.
	Time reversal gives existence on $[-T,0]$ as well.
	
	Moreover, \eqref{eq:iteration-local-bound} and 
	 embedding  (\ref{time uniform}) imply
	\begin{equation*}
		\sup_{|t|\leq T}
		\left(\|u(t)\|_{H^s}+\|n_+(t)\|_2\right)
		\lesssim
		\sup_{|t|\leq T}
		\left(\|Iu(t)\|_{H^1}+\|n_+(t)\|_2\right)
		\lesssim N^{1-s}
	\end{equation*}
	\begin{equation*}
		\lesssim
		(1+T)^{(1-s)/(2s-1-(4-s)\eps)}.
	\end{equation*}
	For any prescribed $\eps_0>0$, we may choose $\eps$
	still smaller so that
	\begin{equation*}
		\frac{1-s}{2s-1-(4-s)\eps}
		\leq \frac{1-s}{2s-1}+\eps_0.
	\end{equation*}
	Returning to the original variables and using
	\begin{equation*}
		\|n_+(t)\|_2^2
		=\|n(t)\|_2^2+\|D^{-1}n_t(t)\|_2^2,
	\end{equation*}
	we obtain
	\begin{equation*}
		\sup_{|t|\leq T}
		\left(
		\|u(t)\|_{H^s}+\|n(t)\|_2
		+\|D^{-1}n_t(t)\|_2
		\right)
		\lesssim_{u_0,n_0,n_1,\eps_0}
		(1+T)^{(1-s)/(2s-1)+\eps_0}.
	\end{equation*}
	This is \eqref{eq:growth}.

\end{proof}

\vspace{15pt}
\part{Blow-up for $s<1/2$}
\section{Local estimates and approximation on a finite interval}
\label{sec:blowup:local}

We first prove that a smooth approximate solution with
a sufficiently small error is close to an exact solution
with the same initial data. All estimates in this section are used at the regularity $H^{1/2}\times L^2$.

For an interval $J$, define
\begin{equation*}
	\|(u,n_+)\|_{\mathcal X(J)}
	=
	\|u\|_{X_S^{1/2,1/2,1}(J)}
	+\|n_+\|_{X_{W_+}^{0,1/2,1}(J)}.
\end{equation*}
The spaces on the right-hand side are those defined in
\eqref{eq:resolution}. We use the sum norm on
$H^{1/2}\times L^2$.

\vspace{8pt}
\subsection{The local estimates}

For $J=[t_0,t_0+h]$, define
\begin{equation*}
	\mathcal I_{S,t_0}F(t)
	=
	-i\int_{t_0}^t e^{i(t-t')\Delta}F(t')\,dt',
\end{equation*}
and 
\begin{equation*}
	\mathcal I_{W_+,t_0}F(t)
	=
	-i\int_{t_0}^t e^{-i(t-t')D}F(t')\,dt'.
\end{equation*}

\begin{lemma}\label{lem:blowup:local-inputs}
	Let $0<h\leq1$ and $J=[t_0,t_0+h]$. Then
	\begin{equation}\label{eq:blowup:free}
		\left\|
		\left(
		e^{i(t-t_0)\Delta}u_0,
		e^{-i(t-t_0)D}n_{+0}
		\right)
		\right\|_{\mathcal X(J)}
		\lesssim
		\|u_0\|_{H^{1/2}}+\|n_{+0}\|_{L^2},
	\end{equation}
	and
	\begin{equation}\label{eq:blowup:embedding}
		\sup_{t\in J}
		\left(
		\|u(t)\|_{H^{1/2}}+\|n_+(t)\|_{L^2}
		\right)
		\lesssim \|(u,n_+)\|_{\mathcal X(J)}.
	\end{equation}
	
    	Moreover, for every sufficiently small $\eps>0$,
    \begin{equation}\label{eq:blowup:bilinear-S}
    	\|\mathcal I_{S,t_0}(un_+)\|_{X_S^{1/2,1/2,1}(J)}
    	+
    	\|\mathcal I_{S,t_0}(u\overline{n_+})\|_{X_S^{1/2,1/2,1}(J)}
    \end{equation}
    \begin{equation*}
    	\lesssim_\eps
    	h^{1/2-\eps}
    	\|u\|_{X_S^{1/2,1/2,1}(J)}
    	\|n_+\|_{X_{W_+}^{0,1/2,1}(J)},
    \end{equation*}
    and
    \begin{equation}\label{eq:blowup:bilinear-W}
    	\|\mathcal I_{W_+,t_0}D(u\overline v)\|_{X_{W_+}^{0,1/2,1}(J)}
    	\lesssim_\eps
    	h^{1/2-\eps}
    	\|u\|_{X_S^{1/2,1/2,1}(J)}
    	\|v\|_{X_S^{1/2,1/2,1}(J)}.
    \end{equation}
    The constants are independent of $t_0$ and $h$.
\end{lemma}

\begin{proof}
	The estimate \eqref{eq:blowup:free} follows from
\cite[Lemma 4.1]{KishimotoLWP}. And
\eqref{eq:blowup:embedding} follows from embedding (\ref{time uniform}) in Part~I. The estimates \eqref{eq:blowup:bilinear-S} and \eqref{eq:blowup:bilinear-W} follow from
\cite[Proposition 4.3 and Remark 4.4]{KishimotoLWP}.

\end{proof}

The error obtained by substituting our approximate solution
into the equations is a sum of terms of the form
$f(t)e^{ik\cdot x-i\omega t}$.
We next estimate their Duhamel integrals to show that this error
produces a small change in the solution.

\begin{lemma}\label{lem:blowup:oscillatory}
Fix $T>0$ and let $f\in C_c^\infty(\mathbb R)$.
	For $*=S,W_+$, write
	\begin{equation*}
		\phi_S(k)=|k|^2,
		\qquad
		\phi_{W_+}(k)=|k|.
	\end{equation*}
	Then, for $J=[t_0,t_1]\subset[0,T]$,
	\begin{equation*}
		\left\|
		\mathcal I_{*,t_0}
		\left(f(t)e^{ik\cdot x-i\omega t}\right)
		\right\|_{X_*^{r,1/2,1}(J)}
		\lesssim_{T,r}\|f\|_{C^1}
		\langle k\rangle^r
		\langle\phi_*(k)-\omega\rangle^{-1/2}.
	\end{equation*}
    The constant is independent of $f,k,\omega,J$.
\end{lemma}

\begin{proof}
	Choose $\chi\in C_c^\infty(\mathbb R)$ equal to one
	on $[0,T]$, depending only on $T$, and set
	\begin{equation*}
		F(t)
		=
		\chi(t)\int_{t_0}^t
		f(t')e^{i(\phi_*(k)-\omega)t'}\,dt'.
	\end{equation*}
	Then
	\begin{equation*}
		\mathcal I_{*,t_0}
		\left(f(t)e^{ik\cdot x-i\omega t}\right)
		=
		-i e^{ik\cdot x-i\phi_*(k)t}F(t)
		\qquad\text{on }J.
	\end{equation*}
	
	We first estimate $F$ in $L^2_t$ and $H^1_t$.
	If $|\phi_*(k)-\omega|>1$, integration by parts gives
	\begin{equation*}
		F(t)
		=
		\frac{\chi(t)}
		{i(\phi_*(k)-\omega)}
		\left(
		f(t)e^{i(\phi_*(k)-\omega)t}
		-f(t_0)e^{i(\phi_*(k)-\omega)t_0}
		\right)
	\end{equation*}
	\begin{equation*}
		\qquad{}
		-\frac{\chi(t)}
		{i(\phi_*(k)-\omega)}
		\int_{t_0}^t
		f'(t')e^{i(\phi_*(k)-\omega)t'}\,dt'.
	\end{equation*}
	For $|\phi_*(k)-\omega|\leq1$, use the defining
	integral directly. In both cases,
	\begin{equation}\label{121}
		\|F\|_{L^2_t}
		\lesssim_T
		\|f\|_{C^1}
		\langle\phi_*(k)-\omega\rangle^{-1}.
	\end{equation}
	Moreover,
	\begin{equation*}
		F'(t)
		=
		\chi'(t)\int_{t_0}^t
		f(t')e^{i(\phi_*(k)-\omega)t'}\,dt'
		+\chi(t)f(t)e^{i(\phi_*(k)-\omega)t},
	\end{equation*}
	so
	\begin{equation}\label{1212}
		\|F\|_{H^1_t}\lesssim_T\|f\|_{C^1}.
	\end{equation}
	
We have
\begin{equation*}
	\mathcal F_{t,x}
	\left(e^{ik\cdot x-i\phi_*(k)t}F(t)\right)(\tau,k')
	=
	\mathbf 1_{\{k'=k\}}
	\widehat F(\tau+\phi_*(k)).
\end{equation*}
For the modulation projections $Q_L^*$ already defined in
\eqref{eq:resolution}, 
Plancherel's identity gives
\begin{equation*}
	\left\|
	Q_L^*\left(e^{ik\cdot x-i\phi_*(k)t}F(t)\right)
	\right\|_{L^2_{t,x}}
	\lesssim
	\left(
	\int_{\langle\tau+\phi_*(k)\rangle\sim L}
	|\widehat F(\tau+\phi_*(k))|^2\,d\tau
	\right)^{1/2}
\end{equation*}
\begin{equation*}
	=
	\left(
	\int_{\langle\xi\rangle\sim L}
	|\widehat F(\xi)|^2\,d\xi
	\right)^{1/2}
	\leq 	\min\left(
	\|F\|_{L^2_t},L^{-1}\|F\|_{H^1_t}
	\right)
\end{equation*}

	Splitting the dyadic sum at
	$L=\|F\|_{H^1_t}/\|F\|_{L^2_t}\geq1$, we can obtain
	
	\begin{equation*}
		\sum_L L^{1/2}
		\min\left(\|F\|_{L^2_t},L^{-1}\|F\|_{H^1_t}\right)
		\lesssim
		\|F\|_{L^2_t}^{1/2}\|F\|_{H^1_t}^{1/2}.
	\end{equation*}
	 Consequently, from (\ref{121}) and (\ref{1212}),
	\begin{equation*}
		\left\|
		\mathcal I_{*,t_0}
		\left(f(t)e^{ik\cdot x-i\omega t}\right)
		\right\|_{X_*^{r,1/2,1}(J)}
		\lesssim_r
		\langle k\rangle^r
		\|F\|_{L^2_t}^{1/2}\|F\|_{H^1_t}^{1/2}
		\lesssim_{T,r}
		\|f\|_{C^1}
		\langle k\rangle^r
		\langle\phi_*(k)-\omega\rangle^{-1/2}.
	\end{equation*}

\end{proof}

\vspace{8pt}
\subsection{Approximation by an exact solution}
We now establish a stability statement for smooth approximate
solutions using Lemma~\ref{lem:blowup:local-inputs} and a standard contraction argument.
We give the proof to make explicit that the constants depend
only on the bound in $\mathcal X$, the time interval, and
the fixed parameters, rather than on higher Sobolev norms.

\begin{proposition}\label{prop:blowup:stability}
	Fix $T>0$, $M\geq1$, and a sufficiently small $\eps>0$.
	Let $(u^{\mathrm{app}},n_+^{\mathrm{app}})$ be smooth
	on $[0,T]\times\T^2$ and satisfy
	\begin{equation*}
		\|(u^{\mathrm{app}},n_+^{\mathrm{app}})\|_{\mathcal X(J)}
		\leq M
	\end{equation*}
	for every interval $J\subset[0,T]$ of length at most one.
	Define the errors
	\begin{equation}\label{eq:blowup:residual-S}
		F_S
		=
		(i\partial_t+\Delta)u^{\mathrm{app}}
		-\frac12
		\left(n_+^{\mathrm{app}}
		+\overline{n_+^{\mathrm{app}}}\right)
		u^{\mathrm{app}},
	\end{equation}
	and
	\begin{equation}\label{eq:blowup:residual-W}
		F_{W_+}
		=
		(i\partial_t-D)n_+^{\mathrm{app}}
		-D|u^{\mathrm{app}}|^2.
	\end{equation}
	
	There exist a partition
	\begin{equation*}
		0=t_0<t_1<\cdots<t_m=T,
		\qquad J_j=[t_j,t_{j+1}],
	\end{equation*}
	and $\rho_0,C>0$, all depending only on
	$M,T,\eps$, with the following property.
	If $0\leq\rho\leq\rho_0$ and
	\begin{equation*}
		\max_{0\leq j<m}
		\left\|
		\left(
		\mathcal I_{S,t_j}F_S,
		\mathcal I_{W_+,t_j}F_{W_+}
		\right)
		\right\|_{\mathcal X(J_j)}
		\leq\rho,
	\end{equation*}
	then there exists a smooth exact solution $(u,n_+)$
	on $[0,T]$ with
	\begin{equation*}
		(u,n_+)(0)
		=
		(u^{\mathrm{app}},n_+^{\mathrm{app}})(0)
	\end{equation*}
	and
	\begin{equation}\label{eq:blowup:stability}
		\max_{0\leq j<m}
		\|(u-u^{\mathrm{app}},
		n_+-n_+^{\mathrm{app}})\|_{\mathcal X(J_j)}
		\leq C\rho.
	\end{equation}
\end{proposition}
\begin{proof}
	To find the exact solution, it suffices to find $(e_u, e_n)$ satisfying
\begin{equation*}
	(i\partial_t+\Delta)e_u
	=
	-F_S
	+\frac12
	\left(n_+^{\mathrm{app}}
	+\overline{n_+^{\mathrm{app}}}\right)e_u
	+\frac12(e_n+\overline{e_n})u^{\mathrm{app}}
	+\frac12(e_n+\overline{e_n})e_u,
\end{equation*}
and
\begin{equation*}
	(i\partial_t-D)e_n
	=
	-F_{W_+}
	+D\left(
	u^{\mathrm{app}}\overline{e_u}
	+e_u\overline{u^{\mathrm{app}}}
	+|e_u|^2
	\right),
\end{equation*}
with initial data $e_u(0)=e_n(0)=0.$

	We solve these equations on successive short intervals.
	Choose $h\leq1$ sufficiently small (to be specified later)
	and partition $[0,T]$ into intervals
	$J_j=[t_j,t_{j+1}]$ of length at most $h$.
	The constant $C_\eps$ is chosen to cover the estimates
	below.
	
    	On $J_j=[t_j,t_{j+1}]$, the values
    $e_u(t_j),e_n(t_j)$ are given by the construction
    on the preceding interval; for $j=0$, both are zero.
    We define the solution map
    \begin{equation*}
    	\mathcal T_j(e_u,e_n)
    	=
    	\left(
    	\mathcal T_{j,u}(e_u,e_n),
    	\mathcal T_{j,n}(e_u,e_n)
    	\right),
    \end{equation*}
    where
    \begin{equation*}
    	\mathcal T_{j,u}(e_u,e_n)(t)
    	=
    	e^{i(t-t_j)\Delta}e_u(t_j)
    	-\mathcal I_{S,t_j}F_S(t)
    	+\frac12\mathcal I_{S,t_j}
    	\left(
    	(n_+^{\mathrm{app}}
    	+\overline{n_+^{\mathrm{app}}})e_u
    	+(e_n+\overline{e_n})u^{\mathrm{app}}
    	+(e_n+\overline{e_n})e_u
    	\right)(t),
    \end{equation*}
    and
    \begin{equation*}
    	\mathcal T_{j,n}(e_u,e_n)(t)
    	=
    	e^{-i(t-t_j)D}e_n(t_j)
    	-\mathcal I_{W_+,t_j}F_{W_+}(t)
    	+\mathcal I_{W_+,t_j}D
    	\left(
    	u^{\mathrm{app}}\overline{e_u}
    	+e_u\overline{u^{\mathrm{app}}}
    	+|e_u|^2
    	\right)(t).
    \end{equation*}
    
    We first show that $\mathcal T_j$ preserves a
    sufficiently small ball.
    Lemma~\ref{lem:blowup:local-inputs} gives
    \begin{equation*}
    	\|\mathcal T_{j,u}(e_u,e_n)\|_{X_S^{1/2,1/2,1}(J_j)}
    	\leq
    	C\|e_u(t_j)\|_{H^{1/2}}
    	+\|\mathcal I_{S,t_j}F_S\|_{X_S^{1/2,1/2,1}(J_j)}
    \end{equation*}
    \begin{equation*}
    	\qquad{}
    	+C_\eps h^{1/2-\eps}
    	\|n_+^{\mathrm{app}}\|_{X_{W_+}^{0,1/2,1}(J_j)}
    	\|e_u\|_{X_S^{1/2,1/2,1}(J_j)}
    \end{equation*}
    \begin{equation*}
    	\qquad{}
    	+C_\eps h^{1/2-\eps}
    	\|e_n\|_{X_{W_+}^{0,1/2,1}(J_j)}
    	\left(
    	\|u^{\mathrm{app}}\|_{X_S^{1/2,1/2,1}(J_j)}
    	+\|e_u\|_{X_S^{1/2,1/2,1}(J_j)}
    	\right),
    \end{equation*}
    and
    \begin{equation*}
    	\|\mathcal T_{j,n}(e_u,e_n)\|_{X_{W_+}^{0,1/2,1}(J_j)}
    	\leq
    	C\|e_n(t_j)\|_{L^2}
    	+\|\mathcal I_{W_+,t_j}F_{W_+}\|_{X_{W_+}^{0,1/2,1}(J_j)}
    \end{equation*}
    \begin{equation*}
    	\qquad{}
    	+C_\eps h^{1/2-\eps}
    	\left(
    	2\|u^{\mathrm{app}}\|_{X_S^{1/2,1/2,1}(J_j)}
    	\|e_u\|_{X_S^{1/2,1/2,1}(J_j)}
    	+\|e_u\|_{X_S^{1/2,1/2,1}(J_j)}^2
    	\right).
    \end{equation*}
    Adding these inequalities and using the assumptions
    on the approximate solution and its errors yields
    \begin{equation}\label{de}
    	\|\mathcal T_j(e_u,e_n)\|_{\mathcal X(J_j)}
    	\leq
    	C\|(e_u,e_n)(t_j)\|_{H^{1/2}\times L^2}
    	+\rho
    	+C_\eps h^{1/2-\eps}
    	\left(
    	M\|(e_u,e_n)\|_{\mathcal X(J_j)}
    	+\|(e_u,e_n)\|_{\mathcal X(J_j)}^2
    	\right).
    \end{equation}
    
    Similarly, we can obtain the estimate for the difference of two outputs:
    \begin{equation*}
    	\|\mathcal T_j(e_u,e_n)
    	-\mathcal T_j(\widetilde e_u,\widetilde e_n)
    	\|_{\mathcal X(J_j)}
    \end{equation*}
    \begin{equation}\label{ed}
    	\leq
    	C_\eps h^{1/2-\eps}
    	\left(
    	M+\|(e_u,e_n)\|_{\mathcal X(J_j)}
    	+\|(\widetilde e_u,\widetilde e_n)\|_{\mathcal X(J_j)}
    	\right)
    	\|(e_u-\widetilde e_u,e_n-\widetilde e_n)
    	\|_{\mathcal X(J_j)}.
    \end{equation}
    
    Choose $h$ such that
    \begin{equation*}
    	C_\eps h^{1/2-\eps}(M+2)\leq\frac14.
    \end{equation*}
    For this interval, set
    \begin{equation*}
    	R_j
    	=
    	2\left(
    	C\|(e_u,e_n)(t_j)\|_{H^{1/2}\times L^2}
    	+\rho
    	\right),
    \end{equation*}
    and suppose $R_j\leq1$.
    On the closed ball of radius $R_j$ in
    $\mathcal X(J_j)$, the first estimate (\ref{de}) gives
    \begin{equation*}
    	\|\mathcal T_j(e_u,e_n)\|_{\mathcal X(J_j)}
    	\leq
    	\frac{R_j}{2}
    	+C_\eps h^{1/2-\eps}(M+R_j)R_j
    	\leq\frac{3R_j}{4}
    	\leq R_j.
    \end{equation*}
    And the second estimate (\ref{ed}) gives
    \begin{equation*}
    	\|\mathcal T_j(e_u,e_n)
    	-\mathcal T_j(\widetilde e_u,\widetilde e_n)
    	\|_{\mathcal X(J_j)}
    	\leq
    	\frac14
    	\|(e_u-\widetilde e_u,e_n-\widetilde e_n)
    	\|_{\mathcal X(J_j)}.
    \end{equation*}
    Banach's fixed point theorem therefore produces an exact solution
    $(e_u,e_n)$ on $J_j$ with
    \begin{equation*}
    	\|(e_u,e_n)\|_{\mathcal X(J_j)}\leq R_j.
    \end{equation*}

    We now verify that $R_j\leq1$ at every step.
    Let $C_{\mathrm{emb}}$ be the constant in (\ref{eq:blowup:embedding}).
    Then
    \begin{equation*}
    	\|(e_u,e_n)(t_{j+1})\|_{H^{1/2}\times L^2}
    	\leq
    	C_{\mathrm{emb}}\|(e_u,e_n)\|_{\mathcal X(J_j)}
    	\leq C_{\mathrm{emb}}R_j
    \end{equation*}
    \begin{equation*}
    	=
    	2C_{\mathrm{emb}}C
    	\|(e_u,e_n)(t_j)\|_{H^{1/2}\times L^2}
    	+2C_{\mathrm{emb}}\rho.
    \end{equation*}
    Choose
    \begin{equation*}
    	C_1\geq
    	\max\left(2,2C_{\mathrm{emb}}C,2C_{\mathrm{emb}}\right).
    \end{equation*}
    Since the initial error is zero, the resulting recurrence
    implies
    \begin{equation*}
    	\|(e_u,e_n)(t_j)\|_{H^{1/2}\times L^2}
    	\leq \rho\sum_{a=1}^{j}C_1^a, \qquad \forall j
    \end{equation*}
    
    Choosing $\rho_0>0$ sufficiently small, we can thus ensure $R_j\le 1$ and (\ref{eq:blowup:stability}).
    It remains to prove that the constructed solution is smooth.
    
    Its initial data agree with those of the smooth pair
    $(u^{\mathrm{app}},n_+^{\mathrm{app}})$.
    Hence, for every integer $m\geq1$,
    \cite[Theorem 1.1]{KishimotoLWP} gives a local solution in
    $H^{m+1/2}\times H^m$, which agrees with $(u,n_+)$
    by local uniqueness. We show by induction on $m$ that
    this higher regularity persists through $T$.
    
    	The case $m=0$ has already been established. Suppose that
    all spatial derivatives of order at most $m-1$ have finite
    $\mathcal X$ norms on the intervals of a finite partition
    of $[0,T]$. On any sufficiently short subinterval $I=[a,b]$
    on which the solution belongs to $H^{m+1/2}\times H^m$,
    differentiate the integral equations $m$ times in space.
    For $|\alpha|=m$, the terms in $\partial^\alpha(n_+u)$
    with all $m$ derivatives on one factor are
    $n_+\partial^\alpha u$ and $(\partial^\alpha n_+)u$.
    Similarly, in $\partial^\alpha D(|u|^2)$, these terms are
    $D((\partial^\alpha u)\overline u)$ and
    $D(u\,\overline{\partial^\alpha u})$.
    In every remaining product, each factor carries at most
    $m-1$ spatial derivatives.
    
    Applying Lemma~\ref{lem:blowup:local-inputs} as before
    and absorbing the terms of order $m$ gives
    \begin{equation*}
    	\sum_{|\alpha|=m}
    	\|(\partial^\alpha u,\partial^\alpha n_+)\|_{\mathcal X(I)}
    	<\infty,
    \end{equation*}
    with a bound depending only on the $H^{m+1/2}\times H^m$
    norm at $a$ and the lower-order bounds. Iterating this
    estimate over a finite subdivision continues the solution at this regularity throughout $[0,T]$.

\end{proof}

\vspace{10pt}
\section{Choice of frequencies and construction of an approximate solution}
\label{sec:blowup:pulse}

We will construct an approximate solution whose $u$ component
is supported at two frequencies $p,l$ and whose $n_+$ component
is supported at $q=p-l$. We choose these frequencies so that
$|p|^2-|l|^2-|q|$ tends to zero. We then choose the three
Fourier coefficients so that the equations hold exactly at
$p,l,q$. At every other frequency where an error occurs,
Lemma~\ref{lem:blowup:oscillatory} will show that its Duhamel
integral tends to zero.

\vspace{8pt}
\subsection{Choice of frequencies on the rectangular torus}

\begin{lemma}\label{lem:blowup:frequencies}
	For every $\gamma_1,\gamma_2>0$, there exist lattice
	frequencies $p_\nu,l_\nu,q_\nu$, indexed by positive even
	integers $\nu$, and numbers $N_\nu\to\infty$ such that
	\begin{equation}\label{eq:blowup:frequency-size}
		q_\nu=p_\nu-l_\nu,
		\qquad
		|p_\nu|\sim|l_\nu|\sim|q_\nu|\sim N_\nu,
	\end{equation}
	and
	\begin{equation}\label{eq:blowup:detuning}
		\delta_\nu
		=
		|p_\nu|^2-|l_\nu|^2-|q_\nu|
		\longrightarrow0.
	\end{equation}
	The first coordinate of $q_\nu$ is positive.
	Moreover, for every fixed nonzero lattice frequency $r$,
	\begin{equation}\label{eq:blowup:fixed-frequency}
		|p_\nu\cdot r|\sim_r N_\nu,
		\qquad
		|l_\nu\cdot r|\sim_r N_\nu
	\end{equation}
	for all sufficiently large $\nu$.
\end{lemma}

\begin{proof}
	Choose an integer $h_\nu$ satisfying
	\begin{equation*}
		\left|2h_\nu-(\gamma_1-\sqrt2\,\nu)\right|\leq1.
	\end{equation*}
	For each fixed $\nu$, choose integers $K_\nu$
	and $v_\nu$ such that
	\begin{equation}\label{eq:blowup:integer-choice}
		\left|
		v_\nu-
		K_\nu\frac{\gamma_2^2}{\gamma_1^2}
		(\gamma_1-2h_\nu)
		\right|
		\leq \nu^{-2}.
	\end{equation}
	Indeed, if
	$\frac{\gamma_2^2}{\gamma_1^2}(\gamma_1-2h_\nu)$
	is rational, take $K_\nu$ to be any sufficiently large
	multiple of its denominator. 
	
	If $\frac{\gamma_2^2}{\gamma_1^2}(\gamma_1-2h_\nu)$ is
	irrational, Dirichlet's approximation theorem gives arbitrarily
	large integers $K_\nu$ and corresponding integers $v_\nu$ such that
	\begin{equation*}
		\left|
		K_\nu\frac{\gamma_2^2}{\gamma_1^2}(\gamma_1-2h_\nu)-v_\nu
		\right|<K_\nu^{-1}.
	\end{equation*}
	We may therefore require $K_\nu\geq\nu^3$ to ensure
	\eqref{eq:blowup:integer-choice}.
	
	Set $N_\nu=\nu K_\nu$ and define
	\begin{equation*}
		p_\nu
		=
		\left(
		\frac{N_\nu+h_\nu}{\gamma_1},
		\frac{v_\nu+\nu/2}{\gamma_2}
		\right),
		\qquad
		l_\nu
		=
		\left(
		\frac{-N_\nu+h_\nu}{\gamma_1},
		\frac{v_\nu-\nu/2}{\gamma_2}
		\right).
	\end{equation*}
	Since $\nu$ is even, these are lattice frequencies, and
	\begin{equation*}
		q_\nu=p_\nu-l_\nu
		=
		\left(
		\frac{2N_\nu}{\gamma_1},
		\frac{\nu}{\gamma_2}
		\right).
	\end{equation*}
	To verify \eqref{eq:blowup:detuning}, compute
	\begin{equation*}
		|p_\nu|^2-|l_\nu|^2
		=
		\frac{4N_\nu h_\nu}{\gamma_1^2}
		+\frac{2\nu v_\nu}{\gamma_2^2},
		\qquad
		|q_\nu|-\frac{2N_\nu}{\gamma_1}
		=
		\frac{\nu^2}
		{\gamma_2^2(|q_\nu|+2N_\nu/\gamma_1)}.
	\end{equation*}
	Hence
	\begin{equation*}
		\delta_\nu
		=
		\frac{2\nu}{\gamma_2^2}
		\left(
		v_\nu-
		K_\nu\frac{\gamma_2^2}{\gamma_1^2}
		(\gamma_1-2h_\nu)
		\right)
		-
		\frac{\nu^2}
		{\gamma_2^2(|q_\nu|+2N_\nu/\gamma_1)}.
	\end{equation*}
	By \eqref{eq:blowup:integer-choice},
	\begin{equation*}
		|\delta_\nu|
		\leq
		\frac{2}{\gamma_2^2\nu}
		+\frac{\gamma_1\nu}{2\gamma_2^2K_\nu}
		\longrightarrow0.
	\end{equation*}
	
	It remains to verify \eqref{eq:blowup:frequency-size}
	and \eqref{eq:blowup:fixed-frequency}.
	The choices of $h_\nu$ and $v_\nu$ give
	\begin{equation*}
		\frac{h_\nu}{N_\nu}\longrightarrow0,
		\qquad
		\frac{v_\nu}{N_\nu}
		=
		\frac{\gamma_2^2}{\gamma_1^2}
		\frac{\gamma_1-2h_\nu}{\nu}
		+O\left(\frac{\nu^{-2}}{N_\nu}\right)
		\longrightarrow
		\frac{\sqrt2\,\gamma_2^2}{\gamma_1^2}.
	\end{equation*}
	Therefore
	\begin{equation*}
		\frac{p_\nu}{N_\nu}
		\longrightarrow
		\left(\frac1{\gamma_1},
		\frac{\sqrt2\,\gamma_2}{\gamma_1^2}\right),
		\qquad
		\frac{l_\nu}{N_\nu}
		\longrightarrow
		\left(-\frac1{\gamma_1},
		\frac{\sqrt2\,\gamma_2}{\gamma_1^2}\right),
	\end{equation*}
	and $q_\nu/N_\nu\to(2/\gamma_1,0)$.
	This proves \eqref{eq:blowup:frequency-size}.
	For $r=(r_1/\gamma_1,r_2/\gamma_2)\neq0$, we also have
	\begin{equation*}
		\frac{p_\nu\cdot r}{N_\nu}
		\longrightarrow
		\frac{r_1+\sqrt2\,r_2}{\gamma_1^2},
		\qquad
		\frac{l_\nu\cdot r}{N_\nu}
		\longrightarrow
		\frac{-r_1+\sqrt2\,r_2}{\gamma_1^2}.
	\end{equation*}
	Neither limit vanishes, since $r_1,r_2$ are integers,
	not both zero, and $\sqrt2$ is irrational. This proves
	\eqref{eq:blowup:fixed-frequency}.
\end{proof}

\vspace{8pt}
\subsection{The three Fourier coefficients}

Fix $\lambda>0$ independently of $\nu$.
For the calculations below, write $p=p_\nu$, $l=l_\nu$,
$q=q_\nu$, and $\delta=\delta_\nu$. We take $\nu$ sufficiently
large that $|\delta|\leq1$. Consider
\begin{equation}\label{eq:blowup:pulse-u}
	u_\nu^{\mathrm{app}}(t,x)
	=
	a(t)e^{ip\cdot x-i|p|^2t}
	+c(t)e^{il\cdot x-i|l|^2t}
\end{equation}
and
\begin{equation}\label{eq:blowup:pulse-n}
	n_{+,\nu}^{\mathrm{app}}(t,x)
	=
	b(t)e^{iq\cdot x-i(|q|+\delta)t}.
\end{equation}
We first explain why we choose the time factor in
\eqref{eq:blowup:pulse-n} as $|q|+\delta$.
Since $q=p-l$, we have
\begin{equation*}
	\widehat{D|u_\nu^{\mathrm{app}}|^2}(t,q)
	=
	|q|
	\left(a(t)e^{-i|p|^2t}\right)
	\overline{\left(c(t)e^{-i|l|^2t}\right)}
	=
	|q|a(t)\overline{c(t)}e^{-i(|p|^2-|l|^2)t}.
\end{equation*}
By the definition of $\delta$,
$|p|^2-|l|^2=|q|+\delta$. We therefore choose
$\widehat{n_{+,\nu}^{\mathrm{app}}}(t,q)
=b(t)e^{-i(|q|+\delta)t}$, so that the two sides of
the wave equation at $q$ have the same exponential factor:
\begin{equation*}
	(i\partial_t-|q|)
	\left(b(t)e^{-i(|q|+\delta)t}\right)
	=
	\left(ib'(t)+\delta b(t)\right)e^{-i(|q|+\delta)t},
\end{equation*}
\begin{equation*}
	\widehat{D|u_\nu^{\mathrm{app}}|^2}(t,q)
	=
	|q|a(t)\overline{c(t)}e^{-i(|q|+\delta)t}.
\end{equation*}

We next compute the two coefficients in the Schr\"odinger
equation. At $p=q+l$, only the product of
$n_{+,\nu}^{\mathrm{app}}$ at $q$ and $u_\nu^{\mathrm{app}}$
at $l$ occurs. Since $|q|+\delta+|l|^2=|p|^2$,
\begin{equation*}
	\frac12
	\widehat{
		\left(n_{+,\nu}^{\mathrm{app}}
		+\overline{n_{+,\nu}^{\mathrm{app}}}\right)
		u_\nu^{\mathrm{app}}}(t,p)
	=
	\frac12
	\left(b(t)e^{-i(|q|+\delta)t}\right)
	\left(c(t)e^{-i|l|^2t}\right)
	=
	\frac12 b(t)c(t)e^{-i|p|^2t}.
\end{equation*}
At $l=p-q$, only the product of
$\overline{n_{+,\nu}^{\mathrm{app}}}$ at $-q$ and
$u_\nu^{\mathrm{app}}$ at $p$ occurs. Since
$|p|^2-(|q|+\delta)=|l|^2$,
\begin{equation*}
	\frac12
	\widehat{
		\left(n_{+,\nu}^{\mathrm{app}}
		+\overline{n_{+,\nu}^{\mathrm{app}}}\right)
		u_\nu^{\mathrm{app}}}(t,l)
	=
	\frac12
	\overline{\left(b(t)e^{-i(|q|+\delta)t}\right)}
	\left(a(t)e^{-i|p|^2t}\right)
	=
	\frac12\overline{b(t)}a(t)e^{-i|l|^2t}.
\end{equation*}
The corresponding left-hand sides are
\begin{equation*}
	(i\partial_t-|p|^2)\left(a(t)e^{-i|p|^2t}\right)
	=ia'(t)e^{-i|p|^2t},
	\qquad
	(i\partial_t-|l|^2)\left(c(t)e^{-i|l|^2t}\right)
	=ic'(t)e^{-i|l|^2t}.
\end{equation*}
Cancelling the common exponential factor in each equation,
we see that the Schr\"odinger equation holds at $p,l$ and
the wave equation holds at $q$ if and only if
\begin{equation}\label{eq:blowup:amplitude-equations}
	ia'=\frac12bc,
	\qquad
	ic'=\frac12\overline b\,a,
	\qquad
	ib'=-\delta b+|q|a\overline c.
\end{equation}

With the parameter $\lambda$, we choose
\begin{equation}\label{eq:blowup:amplitude-a}
	a(t)
	=
	-\sqrt{\frac2{|q|}}\,\lambda
	\left(
	\tanh(\lambda(t-T))+\frac{i\delta}{2\lambda}
	\right),
\end{equation}
\begin{equation}\label{eq:blowup:amplitude-c}
	c(t)
	=
	\sqrt{\frac2{|q|}}\,\lambda
	\operatorname{sech}(\lambda(t-T))
	e^{-i\delta(t-T)/2},
\end{equation}
and
\begin{equation*}
	b(t)
	=
	-2i\lambda
	\operatorname{sech}(\lambda(t-T))
	e^{i\delta(t-T)/2}.
\end{equation*}
Indeed, direct differentiation gives
\begin{equation*}
	ia'
	=
	-i\sqrt{\frac2{|q|}}\,\lambda^2
	\operatorname{sech}^2(\lambda(t-T))
	=\frac12bc,
\end{equation*}
\begin{equation*}
	ic'
	=
	c(t)
	\left(
	-i\lambda\tanh(\lambda(t-T))+\frac{\delta}{2}
	\right)
	=\frac12\overline b\,a,
\end{equation*}
and
\begin{equation*}
	ib'
	=
	\operatorname{sech}(\lambda(t-T))e^{i\delta(t-T)/2}
	\left(
	-2\lambda^2\tanh(\lambda(t-T))+i\lambda\delta
	\right)
	=-\delta b+|q|a\overline c.
\end{equation*}

We record the bounds needed for the construction.
By \eqref{eq:blowup:amplitude-a} and
\eqref{eq:blowup:amplitude-c},
\begin{equation}\label{asdd}
	|a(t)|^2+|c(t)|^2
	=
	\frac{2\lambda^2}{|q|}
	\left(
	\tanh^2(\lambda(t-T))
	+\operatorname{sech}^2(\lambda(t-T))
	+\frac{\delta^2}{4\lambda^2}
	\right)
	=
	\frac{2\lambda^2+\delta^2/2}{|q|}.
\end{equation}
Since $|p|,|l|,|q|\sim N_\nu$, Parseval's identity and (\ref{asdd}) give
\begin{equation}\label{eq:blowup:pulse-small-u}
	\sup_{0\leq t\leq T}
	\|u_\nu^{\mathrm{app}}(t)\|_{H^{s_0}}
	\lesssim_{\lambda}
	N_\nu^{s_0-1/2}
	\longrightarrow0,
\end{equation}
where $\max(s,0)<s_0<1/2$ is the one in (\ref{eq:blowup:s0}).

In contrast, the wave coefficient has modulus
\begin{equation*}
	\left|
	\widehat{n_{+,\nu}^{\mathrm{app}}}(t,q_\nu)
	\right|
	=
	2\lambda\operatorname{sech}(\lambda(t-T)),
\end{equation*}
independently of $\nu$.

Next, we claim that
\begin{equation}\label{eq:blowup:pulse-local-bound}
	\sup_{\substack{J\subset[0,T]}}
	\|(u_\nu^{\mathrm{app}},n_{+,\nu}^{\mathrm{app}})
	\|_{\mathcal X(J)}
	\lesssim_{\lambda,T}1.
\end{equation}
We first verify the estimate for a single spatial frequency.
By the estimate proved in the proof of
Lemma~\ref{lem:blowup:oscillatory}, for $*=S,W_+$ and
$F\in H^1(\mathbb R)$,
\begin{equation}\label{asd}
	\|e^{ik\cdot x-i\phi_*(k)t}F(t)\|_{X_*^{r,1/2,1}}
	\lesssim_r \langle k\rangle^r\|F\|_{H^1_t}.
\end{equation}

From the formulas for $a,c,b$, using $|\delta|\leq1$ and
the boundedness of $\tanh$, $\operatorname{sech}$, and their
first derivatives, we obtain
\begin{equation}\label{C1}
	|a(t)|+|a'(t)|+|c(t)|+|c'(t)|
	\lesssim_\lambda |q|^{-1/2},
	\qquad
	|b(t)|+|b'(t)|\lesssim_\lambda 1.
\end{equation}
Take $\chi\in C_c^\infty$ with $\chi= 1$ on $[0,T]$. The product rule gives
\begin{equation*}
	\|\chi a\|_{H^1_t}+\|\chi c\|_{H^1_t}
	\lesssim_{\lambda,T}|q|^{-1/2},
	\qquad
	\|\chi b e^{-i\delta t}\|_{H^1_t}
	\lesssim_{\lambda,T}1.
\end{equation*}

Applying (\ref{asd})
with $r=1/2$ at $k=p,l$ and with $r=0$ at $k=q$, we find
\begin{equation*}
	\|(u_\nu^{\mathrm{app}},n_{+,\nu}^{\mathrm{app}})\|_{\mathcal X(J)}
	\lesssim
	\langle p\rangle^{1/2}\|\chi a\|_{H^1_t}
	+\langle l\rangle^{1/2}\|\chi c\|_{H^1_t}
	+\|\chi b e^{-i\delta t}\|_{H^1_t}
\end{equation*}
\begin{equation*}
	\lesssim_{\lambda,T}
	\frac{\langle p\rangle^{1/2}+\langle l\rangle^{1/2}}
	{|q|^{1/2}}+1
	\lesssim_{\lambda,T}1.
\end{equation*}
The bound is independent
of $J$ and $\nu$, which proves \eqref{eq:blowup:pulse-local-bound}.

\vspace{8pt}
\subsection{The error in the two equations}

\begin{lemma}\label{lem:blowup:pulse-error}
	Let $F_{S,\nu},F_{W_+,\nu}$ be the errors obtained by
	substituting \eqref{eq:blowup:pulse-u} and
	\eqref{eq:blowup:pulse-n} into
	\eqref{eq:blowup:residual-S} and
	\eqref{eq:blowup:residual-W}.
	For every $J=[t_0,t_1]\subset[0,T]$ of length at most one,
	\begin{equation*}
		\|\mathcal I_{S,t_0}F_{S,\nu}\|_{X_S^{1/2,1/2,1}(J)}
		\lesssim_{\lambda,T}N_\nu^{-1},
	\end{equation*}
	and
	\begin{equation*}
		\|\mathcal I_{W_+,t_0}F_{W_+,\nu}\|_{X_{W_+}^{0,1/2,1}(J)}
		\lesssim_{\lambda,T}N_\nu^{-1/2}.
	\end{equation*}
\end{lemma}

\begin{proof}
	By \eqref{eq:blowup:amplitude-equations}, the errors at
	$p,l$ in the Schr\"odinger equation and at $q$ in the wave
	equation vanish. We compute all remaining terms.
	
	After simple calculations, we can obtain
	\begin{equation*}
		F_{S,\nu}
		=
		-\frac12 ba\,
		e^{i(p+q)\cdot x-i(|p|^2+|q|+\delta)t}
		-\frac12\overline b\,c\,
		e^{i(l-q)\cdot x-i(|l|^2-|q|-\delta)t}.
	\end{equation*}
Since $\chi=1$ on $J$, we may multiply both terms of
$F_{S,\nu}$ by $\chi$ before applying
Lemma~\ref{lem:blowup:oscillatory}. With $r=1/2$, the lemma gives
\begin{equation*}
	\|\mathcal I_{S,t_0}F_{S,\nu}\|_{X_S^{1/2,1/2,1}(J)}
	\lesssim_T
	\frac{\langle p+q\rangle^{1/2}\|\chi ba\|_{C^1}}
	{\langle\phi_S(p+q)-(|p|^2+|q|+\delta)\rangle^{1/2}}
	+\frac{\langle l-q\rangle^{1/2}
		\|\chi\overline b\,c\|_{C^1}}
	{\langle\phi_S(l-q)-(|l|^2-|q|-\delta)\rangle^{1/2}}.
\end{equation*}
We now evaluate the two denominators. Since $\phi_S(k)=|k|^2$,
$p-q=l$, and $|q|+\delta=|p|^2-|l|^2$, we have
\begin{equation*}
	\phi_S(p+q)-(|p|^2+|q|+\delta)
	=
	|p+q|^2-2|p|^2+|p-q|^2
	=
	2|q|^2,
\end{equation*}
and
\begin{equation*}
	\phi_S(l-q)-(|l|^2-|q|-\delta)
	=
	|l-q|^2-2|l|^2+|l+q|^2
	=
	2|q|^2.
\end{equation*}
By the product rule and the preceding bounds for $a,b,c$
and their first derivatives,
\begin{equation*}
	\|\chi ba\|_{C^1}+\|\chi\overline b\,c\|_{C^1}
	\lesssim_T
	\|b\|_{C^1}\bigl(\|a\|_{C^1}+\|c\|_{C^1}\bigr)
	\lesssim_{\lambda,T}|q|^{-1/2}.
\end{equation*}
Substituting these bounds and using
$|p|\sim|l|\sim|q|\sim N_\nu$, we obtain
\begin{equation*}
	\|\mathcal I_{S,t_0}F_{S,\nu}\|_{X_S^{1/2,1/2,1}(J)}
	\lesssim_{\lambda,T}
	\frac{\bigl(\langle p+q\rangle^{1/2}
		+\langle l-q\rangle^{1/2}\bigr)|q|^{-1/2}}
	{\langle 2|q|^2\rangle^{1/2}}
	\lesssim N_\nu^{-1}.
\end{equation*}

Similarly, one can obtain
\begin{equation*}
	F_{W_+,\nu}
	=
	-|q|c\overline a\,
	e^{-iq\cdot x+i(|q|+\delta)t}.
\end{equation*}
and 
\begin{equation*}
	\|\mathcal I_{W_+,t_0}F_{W_+,\nu}\|_{X_{W_+}^{0,1/2,1}(J)}
	\lesssim_{\lambda,T}N_\nu^{-1/2}.
\end{equation*}

\end{proof}

\vspace{10pt}
\section{Adding the approximate solution to a smooth solution}
\label{sec:blowup:addition}

Fix a smooth solution $(u,n_+)$ on $[0,T]$ and fix $\lambda>0$.
We will construct a second exact solution whose initial data
are obtained by adding
$(u_\nu^{\mathrm{app}},n_{+,\nu}^{\mathrm{app}})(0)$ to
$(u,n_+)(0)$. The point is to make the difference between
the new exact solution and the sum
$(u+u_\nu^{\mathrm{app}},n_++n_{+,\nu}^{\mathrm{app}})$
arbitrarily small on the whole interval $[0,T]$ by taking
$\nu$ sufficiently large.

\begin{proposition}\label{prop:blowup:addition}
	Let $(u,n_+)$ be a smooth solution on $[0,T]$ with
	$\widehat{n_+}(0,0)=0$, and fix $\lambda>0$.
	For all sufficiently large $\nu$, there exists a smooth
	solution $(u_\nu,n_{+,\nu})$ on $[0,T]$ with
	\begin{equation}\label{eq:blowup:addition-data}
		(u_\nu,n_{+,\nu})(0)
		=
		(u,n_+)(0)
		+(u_\nu^{\mathrm{app}},n_{+,\nu}^{\mathrm{app}})(0).
	\end{equation}
	Moreover,
	\begin{equation}\label{eq:blowup:addition-error}
		\sup_{0\leq t\leq T}
		\bigl\|
		(u_\nu-u-u_\nu^{\mathrm{app}},
		n_{+,\nu}-n_+-n_{+,\nu}^{\mathrm{app}})(t)
		\bigr\|_{H^{1/2}\times L^2}
		\longrightarrow0
	\end{equation}
	as $\nu\to\infty$. The spatial mean of $n_{+,\nu}$ is zero
	on $[0,T]$.
\end{proposition}

\begin{proof}
	We apply Proposition~\ref{prop:blowup:stability} to
	\begin{equation*}
		\widetilde u_\nu=u+u_\nu^{\mathrm{app}},
		\qquad
		\widetilde n_{+,\nu}=n_++n_{+,\nu}^{\mathrm{app}}.
	\end{equation*}
	Since $(u,n_+)$ is fixed and smooth, its $\mathcal X(J)$
	norm is bounded uniformly over $J\subset[0,T]$ with
	$|J|\leq1$. The same bound holds for
	$(u_\nu^{\mathrm{app}},n_{+,\nu}^{\mathrm{app}})$,
	uniformly in $\nu$, by (\ref{eq:blowup:pulse-local-bound}).
	Consequently, there exists $M\geq1$, independent of $\nu$,
	such that
	\begin{equation*}
		\sup_{\substack{J\subset[0,T]\\ |J|\leq1}}
		\|(\widetilde u_\nu,\widetilde n_{+,\nu})\|_{\mathcal X(J)}
		\leq M.
	\end{equation*}
	Fix the finite partition given by
	Proposition~\ref{prop:blowup:stability} for this $M$.
	We will prove that the Duhamel integrals of the two
	errors tend to zero on each interval of this partition.
	
	\medskip
	\noindent\textbf{1. The errors for the sum.}
	Since $(u,n_+)$ satisfies the equations exactly,
	\begin{equation*}
		(i\partial_t+\Delta)\widetilde u_\nu
		-\frac12
		(\widetilde n_{+,\nu}+\overline{\widetilde n_{+,\nu}})
		\widetilde u_\nu
		=
		F_{S,\nu}
		-\frac12(n_++\overline{n_+})u_\nu^{\mathrm{app}}
		-\frac12
		(n_{+,\nu}^{\mathrm{app}}+\overline{n_{+,\nu}^{\mathrm{app}}})u,
	\end{equation*}
	and
	\begin{equation*}
		(i\partial_t-D)\widetilde n_{+,\nu}
		-D|\widetilde u_\nu|^2
		=
		F_{W_+,\nu}
		-D\left(
		u_\nu^{\mathrm{app}}\overline u
		+u\,\overline{u_\nu^{\mathrm{app}}}
		\right).
	\end{equation*}
	The Duhamel integrals of $F_{S,\nu}$ and $F_{W_+,\nu}$
	already tend to zero from Lemma \ref{lem:blowup:pulse-error}. It remains to estimate the displayed
	products containing one factor from each solution.
	
    	\medskip
    \noindent\textbf{2. Reduction to finitely many frequencies
    	of the fixed solution.}
    Fix an interval $J_j=[t_j,t_{j+1}]$ of the partition chosen above. To estimate the remaining products, write
    \begin{equation*}
    	\mathcal B_{\nu,j}(v,m)
    	=
    	\left(
    	\frac12\mathcal I_{S,t_j}
    	\left((m+\overline m)u_\nu^{\mathrm{app}}
    	+(n_{+,\nu}^{\mathrm{app}}
    	+\overline{n_{+,\nu}^{\mathrm{app}}})v\right),
    	\;
    	\mathcal I_{W_+,t_j}
    	D\left(u_\nu^{\mathrm{app}}\overline v
    	+v\overline{u_\nu^{\mathrm{app}}}\right)
    	\right).
    \end{equation*}
    Our goal is to show that
    $\max_j\|\mathcal B_{\nu,j}(u,n_+)\|_{\mathcal X(J_j)}
    \to0$.
    Define $P_{\leq N_0}$ by
    \begin{equation*}
    	\widehat{P_{\leq N_0}f}(t,k)
    	=\mathbf 1_{\{|k|\leq N_0\}}\widehat f(t,k),
    	\qquad P_{>N_0}=1-P_{\leq N_0}.
    \end{equation*}
    Then
    \begin{equation*}
    	\mathcal B_{\nu,j}(u,n_+)
    	=
    	\mathcal B_{\nu,j}(P_{\leq N_0}u,P_{\leq N_0}n_+)
    	+\mathcal B_{\nu,j}(P_{>N_0}u,P_{>N_0}n_+).
    \end{equation*}
    
    We first estimate the second term uniformly in $\nu$.
    Applying Lemma~\ref{lem:blowup:local-inputs}, one can derive 
    \begin{equation*}
    	\|\mathcal B_{\nu,j}(P_{>N_0}u,P_{>N_0}n_+)\|_{\mathcal X(J_j)}
    	\lesssim_\eps |J_j|^{1/2-\eps}
    	\|P_{>N_0}n_+\|_{X_{W_+}^{0,1/2,1}(J_j)}
    	\|u_\nu^{\mathrm{app}}\|_{X_S^{1/2,1/2,1}(J_j)}
    \end{equation*}
    \begin{equation*}
    	\qquad{}
    	+ |J_j|^{1/2-\eps}
    	\|P_{>N_0}u\|_{X_S^{1/2,1/2,1}(J_j)}
    	\left(
    	\|n_{+,\nu}^{\mathrm{app}}\|_{X_{W_+}^{0,1/2,1}(J_j)}
    	+\|u_\nu^{\mathrm{app}}\|_{X_S^{1/2,1/2,1}(J_j)}
    	\right).
    \end{equation*}
    Thus $|J_j|\leq1$ and the uniform bound for the approximate
    solution imply
    \begin{equation*}
    	\sup_\nu
    	\|\mathcal B_{\nu,j}(P_{>N_0}u,P_{>N_0}n_+)\|_{\mathcal X(J_j)}
    	\lesssim_{\lambda,T,\eps}
    	\|(P_{>N_0}u,P_{>N_0}n_+)\|_{\mathcal X(J_j)}.
    \end{equation*}
    
    Since $u,n_+$ are smooth,
    \begin{equation*}
    	\|u\|_{X_S^{3/2,1/2,1}([0,T])}
    	+\|n_+\|_{X_{W_+}^{1,1/2,1}([0,T])}<\infty.
    \end{equation*}
    On $|k|>N_0$, we have $1\leq N_0^{-1}\langle k\rangle$.
    Thus, by the definition of the norms on $J_j$,
    \begin{equation}\label{asdf}
    	\|(P_{>N_0}u,P_{>N_0}n_+)\|_{\mathcal X(J_j)}
    	\lesssim N_0^{-1}
    	\left(
    	\|u\|_{X_S^{3/2,1/2,1}([0,T])}
    	+\|n_+\|_{X_{W_+}^{1,1/2,1}([0,T])}
    	\right).
    \end{equation}

    We now fix $N_0\gg 1$  and estimate
    $\mathcal B_{\nu,j}(P_{\leq N_0}u,P_{\leq N_0}n_+)$.
    In the following arguments, constants may depend on the fixed
    functions $u,n_+$, on $\lambda,T,N_0$,
    but not on $\nu$ or $j$. As before, write
    $p=p_\nu$, $l=l_\nu$, $q=q_\nu$, and $\delta=\delta_\nu$.
    
    \medskip
    \noindent\textbf{3. The products in the Schr\"odinger equation.}
    We first estimate
    \begin{equation*}
    	\mathcal I_{S,t_j}
    	\left((P_{\leq N_0}n_+
    	+\overline{P_{\leq N_0}n_+})u_\nu^{\mathrm{app}}\right).
    \end{equation*}
    The wave equation at $k=0$ and the hypothesis on the
    initial wave mean give
    \begin{equation*}
    	i\partial_t\widehat{n_+}(t,0)=0,
    	\qquad
    	\widehat{n_+}(t,0)=\widehat{n_+}(0,0)=0
    	\qquad (0\leq t\leq T).
    \end{equation*}
    Since $\widehat{n_+}(t,0)=0$, we have
    \begin{equation*}
    	P_{\leq N_0}n_++\overline{P_{\leq N_0}n_+}
    	=
    	\sum_{0<|r|\leq N_0}
    	\widehat{n_++\overline{n_+}}(t,r)e^{ir\cdot x}.
    \end{equation*}
    Multiplying by $u_\nu^{\mathrm{app}}$ from
    \eqref{eq:blowup:pulse-u} gives
    \begin{equation*}
    	(P_{\leq N_0}n_++\overline{P_{\leq N_0}n_+})
    	u_\nu^{\mathrm{app}}
    	=
    	\sum_{0<|r|\leq N_0}
    	\widehat{n_++\overline{n_+}}(t,r)\,
    	a(t)e^{i(p+r)\cdot x-i|p|^2t}
    \end{equation*}
    \begin{equation*}
    	\qquad{}
    	+\sum_{0<|r|\leq N_0}
    	\widehat{n_++\overline{n_+}}(t,r)\,
    	c(t)e^{i(l+r)\cdot x-i|l|^2t}.
    \end{equation*}
   Applying Lemma~\ref{lem:blowup:oscillatory} with spatial
   exponent $1/2$ to the terms containing $a$ and $c$ separately,
   we obtain
   \begin{equation*}
   	\left\|
   	\mathcal I_{S,t_j}
   	\left((P_{\leq N_0}n_+
   	+\overline{P_{\leq N_0}n_+})u_\nu^{\mathrm{app}}\right)
   	\right\|_{X_S^{1/2,1/2,1}(J_j)}
   \end{equation*}
   \begin{equation*}
   	\lesssim_T
   	\sum_{0<|r|\leq N_0}\left(
   	\frac{\langle p+r\rangle^{1/2}
   		\|\widehat{n_++\overline{n_+}}(\cdot,r)\,a\|_{C^1}}
   	{\langle |p+r|^2-|p|^2\rangle^{1/2}}
   	+
   	\frac{\langle l+r\rangle^{1/2}
   		\|\widehat{n_++\overline{n_+}}(\cdot,r)\,c\|_{C^1}}
   	{\langle |l+r|^2-|l|^2\rangle^{1/2}}\right).
   \end{equation*}
   For each fixed $r\neq0$, the product rule and (\ref{C1}) give
   \begin{equation*}
   	\|\widehat{n_++\overline{n_+}}(\cdot,r)\,a\|_{C^1}
   	+\|\widehat{n_++\overline{n_+}}(\cdot,r)\,c\|_{C^1}
   	\lesssim
   	\|\widehat{n_++\overline{n_+}}(\cdot,r)\|_{C^1}
   	\left(\|a\|_{C^1}+\|c\|_{C^1}\right)
   	\lesssim_r N_\nu^{-1/2}.
   \end{equation*}
   Moreover, Lemma~\ref{lem:blowup:frequencies} gives
   $|p\cdot r|\sim_r N_\nu$ and $|l\cdot r|\sim_r N_\nu$.
   Since $|r|^2$ is fixed, it follows that
   \begin{equation*}
   	\bigl||p+r|^2-|p|^2\bigr|
   	=\bigl|2p\cdot r+|r|^2\bigr|\sim_r N_\nu,
   	\qquad
   	\bigl||l+r|^2-|l|^2\bigr|
   	=\bigl|2l\cdot r+|r|^2\bigr|\sim_r N_\nu.
   \end{equation*}
   There are only finitely many $r$ with $0<|r|\leq N_0$.
   Using $\langle p+r\rangle^{1/2}+\langle l+r\rangle^{1/2}
   \lesssim_{N_0}N_\nu^{1/2}$, we conclude that
   \begin{equation*}
   	\left\|
   	\mathcal I_{S,t_j}
   	\left((P_{\leq N_0}n_+
   	+\overline{P_{\leq N_0}n_+})u_\nu^{\mathrm{app}}\right)
   	\right\|_{X_S^{1/2,1/2,1}(J_j)}
   	\lesssim_{N_0}
   	\frac{N_\nu^{1/2}N_\nu^{-1/2}}{N_\nu^{1/2}}
   	=N_\nu^{-1/2}.
   \end{equation*}
   
   For the second product, we can similarly obtain
   \begin{equation*}
   	(n_{+,\nu}^{\mathrm{app}}
   	+\overline{n_{+,\nu}^{\mathrm{app}}})P_{\leq N_0}u
   	=
   	\sum_{|k|\leq N_0}
   	b(t)\widehat u(t,k)
   	e^{i(k+q)\cdot x-i(|q|+\delta)t}
   \end{equation*}
   \begin{equation*}
   	\qquad{}
   	+\sum_{|k|\leq N_0}
   	\overline{b(t)}\widehat u(t,k)
   	e^{i(k-q)\cdot x+i(|q|+\delta)t}.
   \end{equation*}
   Applying Lemma~\ref{lem:blowup:oscillatory} gives
   \begin{equation*}
   	\left\|\mathcal I_{S,t_j}
   	\left((n_{+,\nu}^{\mathrm{app}}
   	+\overline{n_{+,\nu}^{\mathrm{app}}})P_{\leq N_0}u\right)
   	\right\|_{X_S^{1/2,1/2,1}(J_j)}
   \end{equation*}
   \begin{equation*}
   	\lesssim_T
   	\sum_{|k|\leq N_0}\left(
   	\frac{\langle k+q\rangle^{1/2}
   		\|b\,\widehat u(\cdot,k)\|_{C^1}}
   	{\langle |k+q|^2-|q|-\delta\rangle^{1/2}}
   	+
   	\frac{\langle k-q\rangle^{1/2}
   		\|\overline b\,\widehat u(\cdot,k)\|_{C^1}}
   	{\langle |k-q|^2+|q|+\delta\rangle^{1/2}}\right).
   \end{equation*}
   For $|k|\leq N_0$ and sufficiently large $\nu$,
   $|k\pm q|\sim N_\nu$. Since $|\delta|\leq1$,
   \begin{equation*}
   	\bigl||k+q|^2-|q|-\delta\bigr|\sim N_\nu^2,
   	\qquad
   	\bigl||k-q|^2+|q|+\delta\bigr|\sim N_\nu^2.
   \end{equation*}
   Using $\|b\|_{C^1}\lesssim_\lambda1$ and summing over $k$ gives
   \begin{equation*}
   	\left\|\mathcal I_{S,t_j}
   	\left((n_{+,\nu}^{\mathrm{app}}
   	+\overline{n_{+,\nu}^{\mathrm{app}}})P_{\leq N_0}u\right)
   	\right\|_{X_S^{1/2,1/2,1}(J_j)}
   	\lesssim_{N_0}\frac{N_\nu^{1/2}}{N_\nu}
   	=N_\nu^{-1/2}.
   \end{equation*}
   
   \medskip
   \noindent\textbf{4. The products in the wave equation.}
   Expanding $u_\nu^{\mathrm{app}}\overline{P_{\leq N_0}u}$
   and applying $D$ gives
   \begin{equation*}
   	D\left(u_\nu^{\mathrm{app}}\overline{P_{\leq N_0}u}\right)
   	=
   	\sum_{|k|\leq N_0}
   	|p-k|\,a(t)\overline{\widehat u(t,k)}
   	e^{i(p-k)\cdot x-i|p|^2t}
   \end{equation*}
   \begin{equation*}
   	\qquad{}
   	+\sum_{|k|\leq N_0}
   	|l-k|\,c(t)\overline{\widehat u(t,k)}
   	e^{i(l-k)\cdot x-i|l|^2t}.
   \end{equation*}
   Since $D$ commutes with complex conjugation, the other
   product is
   \begin{equation*}
   	D\left((P_{\leq N_0}u)\overline{u_\nu^{\mathrm{app}}}\right)
   	=
   	\sum_{|k|\leq N_0}
   	|p-k|\,\overline{a(t)}\widehat u(t,k)
   	e^{i(k-p)\cdot x+i|p|^2t}
   \end{equation*}
   \begin{equation*}
   	\qquad{}
   	+\sum_{|k|\leq N_0}
   	|l-k|\,\overline{c(t)}\widehat u(t,k)
   	e^{i(k-l)\cdot x+i|l|^2t}.
   \end{equation*}
   Applying Lemma~\ref{lem:blowup:oscillatory} with spatial
   exponent zero to the four sums yields
   \begin{equation*}
   	\left\|
   	\mathcal I_{W_+,t_j}D
   	\left(
   	u_\nu^{\mathrm{app}}\overline{P_{\leq N_0}u}
   	+(P_{\leq N_0}u)\overline{u_\nu^{\mathrm{app}}}
   	\right)
   	\right\|_{X_{W_+}^{0,1/2,1}(J_j)}
   \end{equation*}
   \begin{equation*}
   	\lesssim_T
   	\sum_{|k|\leq N_0}|p-k|
   	\left(
   	\frac{\|a\,\overline{\widehat u(\cdot,k)}\|_{C^1}}
   	{\langle |p-k|-|p|^2\rangle^{1/2}}
   	+
   	\frac{\|\overline a\,\widehat u(\cdot,k)\|_{C^1}}
   	{\langle |p-k|+|p|^2\rangle^{1/2}}
   	\right)
   \end{equation*}
   \begin{equation*}
   	\qquad{}
   	+\sum_{|k|\leq N_0}|l-k|
   	\left(
   	\frac{\|c\,\overline{\widehat u(\cdot,k)}\|_{C^1}}
   	{\langle |l-k|-|l|^2\rangle^{1/2}}
   	+
   	\frac{\|\overline c\,\widehat u(\cdot,k)\|_{C^1}}
   	{\langle |l-k|+|l|^2\rangle^{1/2}}
   	\right).
   \end{equation*}
   Here $|p-k|+|l-k|\lesssim N_\nu$, whereas
   $|p|^2\sim |l|^2\sim N_\nu^2$. Thus each of the four
   denominators is comparable to $N_\nu$. Since
   $\|a\|_{C^1}+\|c\|_{C^1}\lesssim_\lambda N_\nu^{-1/2}$,
   we obtain
   \begin{equation*}
   	\left\|
   	\mathcal I_{W_+,t_j}D
   	\left(
   	u_\nu^{\mathrm{app}}\overline{P_{\leq N_0}u}
   	+(P_{\leq N_0}u)\overline{u_\nu^{\mathrm{app}}}
   	\right)
   	\right\|_{X_{W_+}^{0,1/2,1}(J_j)}
   	\lesssim_{N_0}
   	\frac{N_\nu N_\nu^{-1/2}}{N_\nu}
   	=N_\nu^{-1/2}.
   \end{equation*}
   
   \medskip
   \noindent\textbf{5. Application of stability.}
   The above three estimates give
   \begin{equation*}
   	\max_j
   	\|\mathcal B_{\nu,j}(P_{\leq N_0}u,P_{\leq N_0}n_+)
   	\|_{\mathcal X(J_j)}
   	\leq C_{N_0}N_\nu^{-1/2}.
   \end{equation*}
   Adding the estimate (\ref{asdf}) for $(P_{>N_0}u,P_{>N_0}n_+)$ yields
   \begin{equation*}
   	\max_j\|\mathcal B_{\nu,j}(u,n_+)\|_{\mathcal X(J_j)}
   	\leq C N_0^{-1}+C_{N_0}N_\nu^{-1/2},
   \end{equation*}
   where $C$ is independent of $N_0,\nu$, and $C_{N_0}$ is
   independent of $\nu$. For each fixed $N_0$, taking the
   upper limit as $\nu\to\infty$ gives
   \begin{equation*}
   	\limsup_{\nu\to\infty}
   	\max_j\|\mathcal B_{\nu,j}(u,n_+)\|_{\mathcal X(J_j)}
   	\leq C N_0^{-1}.
   \end{equation*}
   Letting $N_0\to\infty$ proves that the left-hand side is zero.

   Hence the Duhamel integrals of the errors for $(\widetilde u_\nu, \widetilde n_{+,\nu})$ tend
   to zero on the fixed partition. For all sufficiently large
   $\nu$, Proposition~\ref{prop:blowup:stability} gives a
   smooth exact solution $(u_\nu,n_{+,\nu})$ on $[0,T]$ with
   initial data \eqref{eq:blowup:addition-data} and
   \begin{equation*}
   	\max_j
   	\|(u_\nu-\widetilde u_\nu,
   	n_{+,\nu}-\widetilde n_{+,\nu})\|_{\mathcal X(J_j)}
   	\longrightarrow0.
   \end{equation*}
   The embedding gives \eqref{eq:blowup:addition-error}.
   Finally, $q_\nu\neq0$, so
   $\widehat{n_{+,\nu}}(0,0)=\widehat{n_+}(0,0)=0$.
   Since the symbol of $D$ vanishes at zero, the wave equation
   gives $\partial_t\widehat{n_{+,\nu}}(t,0)=0$.
   Thus $\widehat{n_{+,\nu}}(t,0)=0$ throughout $[0,T]$.
	
\end{proof}

\vspace{10pt}
\section{Proof of the blow-up theorem}

We prove Theorem~\ref{thm:blowup:main}. Fix $T>0$, $\eta>0$,
and $\max(s,0)<s_0<1/2$. We first construct a sequence of
smooth exact solutions on $[0,T]$. We then pass to the limit
on each $[0,T']$ with $T'<T$ and determine the behavior of
the wave component as $t\uparrow T$.

\vspace{8pt}
\subsection{Construction of the smooth solutions}

For $j\geq1$, let $\lambda_j=2^j\Lambda$, where $\Lambda>0$
will be chosen first. By the explicit formula for the wave
coefficient, the approximate solution corresponding to
$\lambda_j$ satisfies
\begin{equation*}
	\left|
	\widehat{n_{+,j}^{\mathrm{app}}}(0,q_j)
	\right|
	=
	2\lambda_j\operatorname{sech}(\lambda_jT).
\end{equation*}
We choose $\Lambda$ sufficiently large that
\begin{equation}\label{eq:blowup:initial-wave-budget}
	\left(
	4|\T^2|
	\sum_{j=1}^\infty
	\lambda_j^2\operatorname{sech}^2(\lambda_jT)
	\right)^{1/2}<\frac{\eta}{4}.
\end{equation}
Thus the sequence $(\lambda_j)$ is fixed before
any spatial frequencies are chosen.

Start with $(u^{(0)},n_+^{(0)})=(0,0)$. At the beginning of
step $j$, the indices $\nu_1,\ldots,\nu_{j-1}$, the frequencies
$p_i,l_i,q_i$ for $i<j$, and the smooth exact solution
$(u^{(j-1)},n_+^{(j-1)})$ on $[0,T]$ have already been chosen.
Assume that
\begin{equation*}
	\widehat{n_+^{(j-1)}}(t,0)=0,
	\qquad 0\leq t\leq T.
\end{equation*}
We now vary only $\nu$, using $\lambda=\lambda_j$ in the
formulas for $(u_\nu^{\mathrm{app}},n_{+,\nu}^{\mathrm{app}})$.

First, from (\ref{eq:blowup:pulse-small-u}),
\begin{equation*}
	\sup_{0\leq t\leq T}\|u_\nu^{\mathrm{app}}(t)\|_{H^{s_0}}
	\leq C_{\lambda_j}N_\nu^{s_0-1/2}
	\longrightarrow0, \qquad \text{if}\;\ \nu\to \infty.
\end{equation*}
Second, Proposition~\ref{prop:blowup:addition}, applied with
$(u,n_+)=(u^{(j-1)},n_+^{(j-1)})$ and $\lambda=\lambda_j$,
gives, for all sufficiently large $\nu$, a smooth exact
solution $(u_\nu,n_{+,\nu})$ satisfying
\begin{equation*}
	(u_\nu,n_{+,\nu})(0)
	=
	(u^{(j-1)},n_+^{(j-1)})(0)
	+(u_\nu^{\mathrm{app}},n_{+,\nu}^{\mathrm{app}})(0)
\end{equation*}
and
\begin{equation*}
	\sup_{0\leq t\leq T}
	\left\|
	(u_\nu,n_{+,\nu})(t)
	-(u^{(j-1)},n_+^{(j-1)})(t)
	-(u_\nu^{\mathrm{app}},n_{+,\nu}^{\mathrm{app}})(t)
	\right\|_{H^{1/2}\times L^2}
	\longrightarrow0,
	\qquad \text{if}\;\ \nu\to\infty.
\end{equation*}
Third, Lemma~\ref{lem:blowup:frequencies} gives
\begin{equation*}
	|\delta_\nu|\longrightarrow0,
	\qquad
	\min\{|p_\nu|,|l_\nu|,|q_\nu|\}\longrightarrow\infty, \qquad \text{if}\;\ \nu\to\infty
\end{equation*}

Choose an even index $\nu_j$ so large that the quantities
in the first two limits are both at most
$\frac{\eta}{8}2^{-j}$, and
\begin{equation*}
	|\delta_{\nu_j}|\leq1,
	\qquad
	\min\{|p_{\nu_j}|,|l_{\nu_j}|,|q_{\nu_j}|\}
	>
	2\max\bigl(
	\{|p_i|,|l_i|,|q_i|:1\leq i<j\}\bigr).
\end{equation*}

With this choice, set
\begin{equation*}
	(p_j,l_j,q_j,\delta_j)
	=(p_{\nu_j},l_{\nu_j},q_{\nu_j},\delta_{\nu_j}),
	\qquad
	(u^{(j)},n_+^{(j)})=(u_{\nu_j},n_{+,\nu_j}),
\end{equation*}
and write
$(u_j^{\mathrm{app}},n_{+,j}^{\mathrm{app}})
=(u_{\nu_j}^{\mathrm{app}},n_{+,\nu_j}^{\mathrm{app}})$.
Then
\begin{equation}\label{eq:blowup:small-app-j}
	\sup_{0\leq t\leq T}\|u_j^{\mathrm{app}}(t)\|_{H^{s_0}}
	\leq\frac{\eta}{8}\,2^{-j}.
\end{equation}
Defining
\begin{equation*}
	(e_{u,j},e_{n,j})
	=
	(u^{(j)},n_+^{(j)})
	-(u^{(j-1)},n_+^{(j-1)})
	-(u_j^{\mathrm{app}},n_{+,j}^{\mathrm{app}}),
\end{equation*}
we also have
\begin{equation}\label{eq:blowup:small-error-j}
	\sup_{0\leq t\leq T}
	\|(e_{u,j},e_{n,j})(t)\|_{H^{1/2}\times L^2}
	\leq\frac{\eta}{8}\,2^{-j}.
\end{equation}
The equality of the initial data gives
\begin{equation}\label{eq:blowup:error-at-zero}
	e_{u,j}(0)=e_{n,j}(0)=0.
\end{equation}

Since the selected frequencies belong to the sequence in
Lemma~\ref{lem:blowup:frequencies},
\begin{equation}\label{eq:blowup:chosen-frequencies}
	|p_j|\sim|l_j|\sim|q_j|,
	\qquad (q_j)_1>0,
\end{equation}
with comparison constants depending only on the periods $\gamma_1, \gamma_2$.
Within each step, $p_j-l_j=q_j\neq0$; between different
steps, the frequency magnitudes are separated by our choice
of $\nu_j$. Hence all the $p_j,l_j$ are distinct, and the
$q_j$ are distinct. Finally, Proposition~\ref{prop:blowup:addition}
gives $\widehat{n_+^{(j)}}(t,0)=0$ on $[0,T]$, as required
for step $j+1$.

\vspace{8pt}
\subsection{Passage to the limit and the initial data}

Summing the definitions of $e_{u,j}$ and $e_{n,j}$ gives
the exact identities
\begin{equation}\label{eq:blowup:partial-sums}
	u^{(m)}
	=\sum_{j=1}^m u_j^{\mathrm{app}}
	+\sum_{j=1}^m e_{u,j},
	\qquad
	n_+^{(m)}
	=\sum_{j=1}^m n_{+,j}^{\mathrm{app}}
	+\sum_{j=1}^m e_{n,j}.
\end{equation}
By \eqref{eq:blowup:small-error-j}, the series
\begin{equation*}
	(e_u,e_n):=\sum_{j=1}^\infty(e_{u,j},e_{n,j})
\end{equation*}
converges in $C([0,T];H^{1/2}\times L^2)$ and has norm
at most $\eta/8$. By \eqref{eq:blowup:small-app-j} and
$H^{1/2}\hookrightarrow H^{s_0}$, the first identity in
\eqref{eq:blowup:partial-sums} therefore gives
\begin{equation}\label{eq:blowup:u-limit}
	u^{(m)}\longrightarrow
	u:=\sum_{j=1}^\infty u_j^{\mathrm{app}}+e_u
	\quad\text{in }C([0,T];H^{s_0}).
\end{equation}

For the wave component, fix $0<T'<T$. By orthogonality
of the distinct frequencies $q_j$, for $m'>m$,
\begin{equation*}
	\sup_{0\leq t\leq T'}
	\left\|\sum_{j=m+1}^{m'}n_{+,j}^{\mathrm{app}}(t)
	\right\|_{L^2}^2
	\leq
	4|\T^2|
	\sum_{j>m}\lambda_j^2
	\operatorname{sech}^2(\lambda_j(T-T')).
\end{equation*}
Since $T-T'>0$ and $\lambda_j=2^j\Lambda$, the series on
the right converges by exponential decay. Its tail tends
to zero as $m\to\infty$. The second identity in
\eqref{eq:blowup:partial-sums} now gives
\begin{equation}\label{eq:blowup:wave-limit}
	n_+^{(m)}\longrightarrow
	n_+:=\sum_{j=1}^\infty n_{+,j}^{\mathrm{app}}+e_n
	\quad\text{in }C([0,T'];L^2).
\end{equation}
In particular, $n_+$ has zero spatial mean.

We next verify the equations for $(u,n_+)$. On $[0,T']$,
both convergences \eqref{eq:blowup:u-limit} and
\eqref{eq:blowup:wave-limit} hold in $C_tL_x^2$.
Thus, uniformly in $t\in[0,T']$,
\begin{equation*}
	\|n_+^{(m)}u^{(m)}-n_+u\|_{L^1}
	\leq
	\|n_+^{(m)}-n_+\|_{L^2}\|u^{(m)}\|_{L^2}
	+\|n_+\|_{L^2}\|u^{(m)}-u\|_{L^2}
	\longrightarrow0,
\end{equation*}
and
\begin{equation*}
	\bigl\||u^{(m)}|^2-|u|^2\bigr\|_{L^1}
	\leq
	\left(\|u^{(m)}\|_{L^2}+\|u\|_{L^2}\right)
	\|u^{(m)}-u\|_{L^2}
	\longrightarrow0.
\end{equation*}
The first inequality also holds with the wave components
conjugated. In the wave equation, we apply $D$ to the
smooth test function when passing to the limit in
$D|u^{(m)}|^2$. Since each $(u^{(m)},n_+^{(m)})$ is an
exact solution, we obtain, in distributions on $[0,T)$,
\begin{equation*}
	i\partial_t u+\Delta u
	=\frac12(n_++\overline{n_+})u,
	\qquad
	i\partial_t n_+-Dn_+=D|u|^2.
\end{equation*}

Set $n=\operatorname{Re}n_+$. Taking real and imaginary
parts of the wave equation gives
\begin{equation*}
	\partial_t n=D\operatorname{Im}n_+,
	\qquad
	\partial_t\operatorname{Im}n_+=-Dn-D|u|^2.
\end{equation*}
By \eqref{eq:blowup:wave-limit},
$n=\operatorname{Re}n_+\in C([0,T);L^2)$. Since
$D:L^2\to H^{-1}$ is bounded, the first identity also gives
$n\in C^1([0,T);H^{-1})$. Together with
$u\in C([0,T];H^{s_0})$ from \eqref{eq:blowup:u-limit},
this proves \eqref{eq:blowup:regularity}.
Applying $D$ to the second identity yields
\begin{equation*}
	\partial_t^2n=-D^2n-D^2|u|^2
	=\Delta n+\Delta|u|^2.
\end{equation*}
Thus $(u,n)$ solves the original Zakharov system in
distributions. Each smooth solution conserves its
Schr\"odinger mass. For every $0\leq t<T$, taking the
limit in
\begin{equation*}
	\|u^{(m)}(t)\|_{L^2}=\|u^{(m)}(0)\|_{L^2}
\end{equation*}
by \eqref{eq:blowup:u-limit} proves \eqref{eq:blowup:mass}.

We next verify the assertions about the initial data.
By \eqref{eq:blowup:error-at-zero},
\begin{equation*}
	u_0=u(0)=\sum_{j=1}^\infty u_j^{\mathrm{app}}(0),
	\qquad
	n_+(0)=\sum_{j=1}^\infty n_{+,j}^{\mathrm{app}}(0).
\end{equation*}
It follows from \eqref{eq:blowup:small-app-j},
orthogonality of the distinct frequencies $q_j$, and
\eqref{eq:blowup:initial-wave-budget} that
\begin{equation*}
	\|u_0\|_{H^{s_0}}\leq\frac{\eta}{8},
	\qquad
	\|n_+(0)\|_{L^2}<\frac{\eta}{4}.
\end{equation*}
Writing $n_0=n(0)$ and $n_1=\partial_t n(0)$, we have
\begin{equation*}
	\|n_0\|_{L^2}\leq\|n_+(0)\|_{L^2},
	\qquad
	\|n_1\|_{H^{-1}}
	=\|D\operatorname{Im}n_+(0)\|_{H^{-1}}
	\leq\|n_+(0)\|_{L^2}.
\end{equation*}
Both wave data have zero mean. Since $s<s_0$ and $s_0>0$,
we obtain the upper bounds in \eqref{eq:blowup:smallness}:
\begin{equation*}
	\|u_0\|_{H^s}+\|n_0\|_{L^2}+\|n_1\|_{H^{-1}}<\eta,
	\qquad
	\|u_0\|_{L^2}<\eta.
\end{equation*}

We now show that $u_0\notin H^{1/2}$. Since all the
$p_j,l_j$ are distinct, only $u_j^{\mathrm{app}}(0)$
contributes to $\widehat u_0(p_j)$. The explicit formula
for $a(0)$ at step $j$ gives
\begin{equation*}
	\langle p_j\rangle|\widehat u_0(p_j)|^2
	=
	\frac{2\lambda_j^2\langle p_j\rangle}{|q_j|}
	\left(
	\tanh^2(\lambda_jT)
	+\frac{\delta_j^2}{4\lambda_j^2}
	\right).
\end{equation*}
By \eqref{eq:blowup:chosen-frequencies},
$\langle p_j\rangle/|q_j|$ is bounded below uniformly in
$j$, and $\tanh(\lambda_jT)\to1$. Therefore
\begin{equation*}
	\|u_0\|_{H^{1/2}}^2
	\geq
	|\T^2|
	\sum_{j=1}^\infty\langle p_j\rangle|\widehat u_0(p_j)|^2
	\gtrsim\sum_{j\ \mathrm{sufficiently\ large}}\lambda_j^2
	=\infty.
\end{equation*}
In particular, $u_0\neq0$, so $0<\|u_0\|_{L^2}<\eta$.
This completes the proof of \eqref{eq:blowup:smallness}.

\vspace{8pt}
\subsection{Blow-up of the wave component}

It remains to prove \eqref{eq:blowup:n-rate} and
\eqref{eq:blowup:nt-rate}. For $\tau=T-t>0$,
orthogonality gives
\begin{equation}\label{eq:blowup:wave-series}
	\left\|\sum_{j=1}^\infty n_{+,j}^{\mathrm{app}}(t)
	\right\|_{L^2}^2
	=
	4|\T^2|
	\sum_{j=1}^\infty
	\lambda_j^2\operatorname{sech}^2(\lambda_j\tau).
\end{equation}
For $0<\tau<\lambda_1^{-1}$, let $j_0$ be the first index
such that $\lambda_{j_0}\tau\geq1$. Since
$\lambda_{j_0}=2\lambda_{j_0-1}$,
\begin{equation*}
	1\leq\lambda_{j_0}\tau<2.
\end{equation*}
The single term with $j=j_0$ gives the lower bound
\begin{equation*}
	\lambda_{j_0}^2\operatorname{sech}^2(\lambda_{j_0}\tau)
	\geq \tau^{-2}\operatorname{sech}^2(2).
\end{equation*}
For the upper bound, first sum over $j<j_0$:
\begin{equation*}
	\sum_{j<j_0}\lambda_j^2\operatorname{sech}^2(\lambda_j\tau)
	\leq\sum_{j<j_0}\lambda_j^2
	\lesssim\lambda_{j_0}^2
	\lesssim\tau^{-2}.
\end{equation*}
For $j=j_0+k$, $k\geq0$, we have
$2^k\leq\lambda_j\tau<2^{k+1}$. Consequently,
\begin{equation*}
	\sum_{j\geq j_0}\lambda_j^2\operatorname{sech}^2(\lambda_j\tau)
	\lesssim
	\tau^{-2}\sum_{k=0}^\infty 4^k e^{-2^{k+1}}
	\lesssim\tau^{-2}.
\end{equation*}
Substituting these bounds into \eqref{eq:blowup:wave-series}
gives
\begin{equation}\label{eq:blowup:wave-sum-rate}
	\left\|\sum_{j=1}^\infty n_{+,j}^{\mathrm{app}}(t)
	\right\|_{L^2}
	\sim (T-t)^{-1}
	\qquad (t\uparrow T).
\end{equation}

By \eqref{eq:blowup:chosen-frequencies}, the Fourier support
of $\sum_j n_{+,j}^{\mathrm{app}}(t)$ lies in $\{k_1>0\}$,
whereas the support of its complex conjugate lies in
$\{k_1<0\}$. These two functions are orthogonal in $L^2$.
Thus
\begin{equation*}
	\left\|\operatorname{Re}\sum_j n_{+,j}^{\mathrm{app}}(t)
	\right\|_{L^2}^2
	=
	\left\|\operatorname{Im}\sum_j n_{+,j}^{\mathrm{app}}(t)
	\right\|_{L^2}^2
	=
	\frac12
	\left\|\sum_j n_{+,j}^{\mathrm{app}}(t)\right\|_{L^2}^2.
\end{equation*}
By \eqref{eq:blowup:wave-limit}, the difference between
$n_+(t)$ and this sum is $e_n(t)$, with
$\sup_{0\leq t\leq T}\|e_n(t)\|_{L^2}\leq\eta/8$.
The triangle inequality in both directions and
\eqref{eq:blowup:wave-sum-rate} therefore give
\begin{equation*}
	\|n(t)\|_{L^2}
	=\|\operatorname{Re}n_+(t)\|_{L^2}
	\sim (T-t)^{-1},
	\qquad
	\|\operatorname{Im}n_+(t)\|_{L^2}
	\sim (T-t)^{-1}
	\qquad (t\uparrow T).
\end{equation*}
The first estimate proves \eqref{eq:blowup:n-rate}.

Finally, $\operatorname{Im}n_+(t)$ has zero mean.
For every mean-zero $f\in L^2(\T^2)$,
\begin{equation*}
	\|Df\|_{H^{-1}}^2
	=
	|\T^2|
	\sum_{k\neq0}\frac{|k|^2}{1+|k|^2}|\widehat f(k)|^2
	\sim \|f\|_{L^2}^2,
\end{equation*}
because the nonzero frequencies of the fixed torus are
bounded away from zero. Taking $f=\operatorname{Im}n_+(t)$
and using $\partial_t n=D\operatorname{Im}n_+$, we obtain
\begin{equation*}
	\|\partial_t n(t)\|_{H^{-1}}
	\sim\|\operatorname{Im}n_+(t)\|_{L^2}
	\sim(T-t)^{-1}
	\qquad (t\uparrow T).
\end{equation*}
This proves \eqref{eq:blowup:nt-rate} and completes the
proof of Theorem~\ref{thm:blowup:main}.

\clearpage
\section*{Acknowledgement}
The author thanks Prof. Alex Cohen for providing financial support to cover the costs of ChatGPT usage during the preparation of this work.

\section*{Conflict of interest statement}
The author does not have any possible conflict of interest.

\section*{Data availability statement}
The manuscript has no associated data.

\bibliographystyle{plain}
\bibliography{zakharov_gwp}

\begin{thebibliography}{10}

\bibitem{BGHN2015}
Ioan Bejenaru, Zihua Guo, Sebastian Herr, and Kenji Nakanishi.
\newblock Well-posedness and scattering for the {Zakharov} system in four
  dimensions.
\newblock {\em Analysis \& PDE}, 8(8):2029--2055, 2015.

\bibitem{BejenaruHerr2011}
Ioan Bejenaru and Sebastian Herr.
\newblock Convolutions of singular measures and applications to the {Zakharov}
  system.
\newblock {\em Journal of Functional Analysis}, 261(2):478--506, 2011.

\bibitem{BHHT2009}
Ioan Bejenaru, Sebastian Herr, Justin Holmer, and Daniel Tataru.
\newblock On the {2D} {Zakharov} system with {$L^2$}-{Schr\"odinger} data.
\newblock {\em Nonlinearity}, 22(5):1063--1089, 2009.

\bibitem{Bourgain1994}
Jean Bourgain.
\newblock On the {Cauchy} and invariant measure problem for the periodic
  {Zakharov} system.
\newblock {\em Duke Mathematical Journal}, 76(1):175--202, 1994.

\bibitem{BourgainColliander1996}
Jean Bourgain and James Colliander.
\newblock On wellposedness of the {Zakharov} system.
\newblock {\em International Mathematics Research Notices}, 1996(11):515--546,
  1996.

\bibitem{Candy2024}
Timothy Candy.
\newblock Concentration compactness for the energy critical {Zakharov} system.
\newblock {\em Discrete and Continuous Dynamical Systems}, 44(5):1395--1445,
  2024.

\bibitem{CandyHerrNakanishi2021}
Timothy Candy, Sebastian Herr, and Kenji Nakanishi.
\newblock Global wellposedness for the energy-critical {Zakharov} system below
  the ground state.
\newblock {\em Advances in Mathematics}, 384:107746, 2021.

\bibitem{CandyHerrNakanishi2023}
Timothy Candy, Sebastian Herr, and Kenji Nakanishi.
\newblock The {Zakharov} system in dimension {$d\geq4$}.
\newblock {\em Journal of the European Mathematical Society}, 25(8):3177--3228,
  2023.

\bibitem{CandyNakanishi2026}
Timothy Candy and Kenji Nakanishi.
\newblock Scattering for the {4D} {Zakharov} system below the ground state,
  2026.
\newblock Preprint.

\bibitem{ChenWu2021}
Zijun Chen and Shengkun Wu.
\newblock Local well-posedness for the {Zakharov} system in dimension
  {$d=2,3$}.
\newblock {\em Communications on Pure and Applied Analysis}, 20(12):4307--4319,
  2021.

\bibitem{CollianderHolmerTzirakis2008}
James Colliander, Justin Holmer, and Nikolaos Tzirakis.
\newblock Low regularity global well-posedness for the {Zakharov} and
  {Klein--Gordon--Schr\"odinger} systems.
\newblock {\em Transactions of the American Mathematical Society},
  360(9):4619--4638, 2008.

\bibitem{CKSTT2002}
James Colliander, Markus Keel, Gigliola Staffilani, Hideo Takaoka, and Terence
  Tao.
\newblock Almost conservation laws and global rough solutions to a nonlinear
  {Schr\"odinger} equation.
\newblock {\em Mathematical Research Letters}, 9(5--6):659--682, 2002.

\bibitem{Compaan2019}
Erin Compaan.
\newblock A note on global existence for the {Zakharov} system on {$\mathbb
  T$}.
\newblock {\em Communications on Pure and Applied Analysis}, 18(5):2473--2489,
  2019.

\bibitem{evans2022partial}
Lawrence~C Evans.
\newblock {\em Partial differential equations}, volume~19.
\newblock American mathematical society, 2022.

\bibitem{FangPecherZhong2009}
Daoyuan Fang, Hartmut Pecher, and Sijia Zhong.
\newblock Low regularity global well-posedness for the two-dimensional
  {Zakharov} system.
\newblock {\em Analysis}, 29:265--281, 2009.

\bibitem{GinibreTsutsumiVelo1997}
Jean Ginibre, Yoshio Tsutsumi, and Giorgio Velo.
\newblock On the {Cauchy} problem for the {Zakharov} system.
\newblock {\em Journal of Functional Analysis}, 151(2):384--436, 1997.

\bibitem{GlangetasMerle1994II}
L.~Glangetas and Frank Merle.
\newblock Concentration properties of blow-up solutions and instability results
  for {Zakharov} equation in dimension two. {II}.
\newblock {\em Communications in Mathematical Physics}, 160(2):349--389, 1994.

\bibitem{GlangetasMerle1994I}
L.~Glangetas and Frank Merle.
\newblock Existence of self-similar blow-up solutions for {Zakharov} equation
  in dimension two. {I}.
\newblock {\em Communications in Mathematical Physics}, 160(1):173--215, 1994.

\bibitem{Grube2025}
Florian Grube.
\newblock Norm inflation for the {Zakharov} system.
\newblock {\em Communications on Pure and Applied Analysis}, 24(9):1804--1829,
  2025.

\bibitem{GuoNakanishi2021}
Zihua Guo and Kenji Nakanishi.
\newblock The {Zakharov} system in {4D} radial energy space below the ground
  state.
\newblock {\em American Journal of Mathematics}, 143(5):1527--1600, 2021.

\bibitem{GuoNakanishi2026}
Zihua Guo and Kenji Nakanishi.
\newblock Global dynamics above the ground state energy for the {3D} {Zakharov}
  system, 2026.
\newblock Preprint.

\bibitem{GuoNakanishiWang2013}
Zihua Guo, Kenji Nakanishi, and Shuxia Wang.
\newblock Global dynamics below the ground state energy for the {Zakharov}
  system in the {3D} radial case.
\newblock {\em Advances in Mathematics}, 238:412--441, 2013.

\bibitem{HaniPusateriShatah2013}
Zaher Hani, Fabio Pusateri, and Jalal Shatah.
\newblock Scattering for the {Zakharov} system in 3 dimensions.
\newblock {\em Communications in Mathematical Physics}, 322(3):731--753, 2013.

\bibitem{KatoKinoshita2023}
Isao Kato and Shinya Kinoshita.
\newblock Global well-posedness and scattering for the {Zakharov} system at the
  critical space in three spatial dimensions with small and radial initial
  data.
\newblock {\em Journal of Mathematical Analysis and Applications},
  518(1):126667, 2023.

\bibitem{KatoTsugawa2017}
Isao Kato and Kotaro Tsugawa.
\newblock Scattering and well-posedness for the {Zakharov} system at a critical
  space in four and more spatial dimensions.
\newblock {\em Differential and Integral Equations}, 30(9--10):763--794, 2017.

\bibitem{KNS2023}
Shinya Kinoshita, Shohei Nakamura, and Akansha Sanwal.
\newblock Decoupling inequality for paraboloid under shell type restriction and
  its application to the periodic {Zakharov} system.
\newblock {\em Communications in Mathematical Physics}, 404(2):883--926, 2023.

\bibitem{KishimotoLWP}
Nobu Kishimoto.
\newblock Local well-posedness for the {Zakharov} system on the
  multidimensional torus.
\newblock {\em Journal d'Analyse Math\'ematique}, 119:213--253, 2013.

\bibitem{KishimotoI}
Nobu Kishimoto.
\newblock Resonant decomposition and the {$I$}-method for the two-dimensional
  {Zakharov} system.
\newblock {\em Discrete and Continuous Dynamical Systems}, 33(9):4095--4122,
  2013.

\bibitem{Kishimoto2022}
Nobu Kishimoto.
\newblock Remarks on periodic {Zakharov} systems.
\newblock {\em Electronic Journal of Differential Equations}, 2022(20):1--19,
  2022.

\bibitem{KishimotoMaeda}
Nobu Kishimoto and Masaya Maeda.
\newblock Construction of blow-up solutions for {Zakharov} system on
  {$\mathbb{T}^2$}.
\newblock {\em Annales de l'Institut Henri Poincar\'e C, Analyse non
  lin\'eaire}, 30(5):791--824, 2013.

\bibitem{KriegerSchmid2024I}
Joachim Krieger and Tobias Schmid.
\newblock Finite time blow up for the energy critical {Zakharov} system {I}:
  approximate solutions, 2024.
\newblock Preprint.

\bibitem{KriegerSchmid2024II}
Joachim Krieger and Tobias Schmid.
\newblock Finite time blow up for the energy critical {Zakharov} system {II}:
  exact solutions, 2024.
\newblock Preprint.

\bibitem{OzawaTsutsumi1992}
Tohru Ozawa and Yoshio Tsutsumi.
\newblock Existence and smoothing effect of solutions for the {Zakharov}
  equations.
\newblock {\em Publications of the Research Institute for Mathematical
  Sciences}, 28(3):329--361, 1992.

\bibitem{Sanwal2022}
Akansha Sanwal.
\newblock Local well-posedness for the {Zakharov} system in dimension
  {$d\leq3$}.
\newblock {\em Discrete and Continuous Dynamical Systems}, 42(3):1067--1103,
  2022.

\bibitem{Takaoka1999}
Hideo Takaoka.
\newblock Well-posedness for the {Zakharov} system with the periodic boundary
  condition.
\newblock {\em Differential and Integral Equations}, 12(6):789--810, 1999.

\bibitem{Zakharov1972}
Vladimir~E. Zakharov.
\newblock Collapse of {Langmuir} waves.
\newblock {\em Soviet Physics JETP}, 35:908--914, 1972.

\end{thebibliography}
\bigskip
\noindent
{\small
	Courant Institute of Mathematical Sciences,
	New York University, New York, NY 10012, USA\par
	\noindent
	\textit{Email:} \texttt{jw9409@nyu.edu}\par
}
\end{document}